%% file: main.tex
\documentclass[12pt, reqno]{amsart}
\input{preamble}

\begin{document}

\title[Mordell--Lang. Selmer groups. superelliptic curves. global function field]{Mordell--Lang and disparate Selmer ranks of odd twists of some superelliptic curves over global function fields}
\author{Sun Woo Park}
\address{Max Planck Institute for Mathematics, Vivatsgasse 7, 53115 Bonn, Germany}
\email{} 
\email{\href{mailto:s.park@mpim-bonn.mpg.de}{s.park@mpim-bonn.mpg.de}}

\begin{abstract}
    Fix a prime number $\ell \geq 5$. Let $K = \mathbb{F}_q(t)$ be a global function field of characteristic $p$ coprime to $2,3$, and $q \equiv 1 \text{ mod } \ell$. Let $C:y^\ell = F(x)$ be a non-isotrivial superelliptic curve over $K$ such that $F$ is a degree $3$ polynomial over $\mathbb{F}_q(t)$. Denote by $C_f: fy^\ell = F(x)$ the twist of $C$ by a polynomial $f$ over $\mathbb{F}_q$. Assuming some conditions on $C$, we show that the expected number of $K$-rational points of $C_f$ is bounded, and at least $99\%$ of such curves $C_f$ have at most $(3p)^{4\ell} \cdot \ell!$ many $K$-rational points, as $f$ ranges over the set of polynomials of sufficiently large degree over $\mathbb{F}_q$. To achieve this, we compute the distribution of dimensions of $1-\zeta_\ell$ Selmer groups of Jacobians of such superelliptic curves. This is done by generalizing the technique of constructing a governing Markov operator, as developed from previous studies by Swinnerton-Dyer, Klagsbrun--Mazur--Rubin, Yu, and the author.
\end{abstract}

\maketitle
\tableofcontents

\input{Introduction}
\input{Alternating}

\input{Application}

\newpage

\bibliographystyle{alpha}
\bibliography{main}

\end{document}

%% file: preamble.tex
\usepackage[margin=0.75in]{geometry}
\usepackage[
  bookmarks=true]{hyperref}
\usepackage[OT2, T1]{fontenc}
\usepackage[english]{babel}
\usepackage[utf8]{inputenc}
\usepackage{csquotes}
\usepackage[final]{microtype}
\usepackage{lmodern}
\usepackage{amsthm}
\usepackage{amssymb}
\usepackage{mathrsfs}
\usepackage{enumerate}
\usepackage{tikz-cd} 
\usepackage{tikz}
\usepackage{multicol}
\usepackage{multirow}
\usepackage{tikz-qtree}
\usepackage{rotating}
\usepackage{graphicx}
\usepackage{comment}
\usepackage{enumitem}
\usepackage{manfnt}
\usetikzlibrary{arrows,calc,matrix,trees,arrows.meta,positioning,decorations.pathreplacing,bending}

\newtheorem{theorem}{Theorem}[section]

\newtheorem{proposition}[theorem]{Proposition}
\newtheorem{lemma}[theorem]{Lemma}
\newtheorem{corollary}[theorem]{Corollary}

\theoremstyle{definition}

\newtheorem{remark}[theorem]{Remark}
\newtheorem{definition}[theorem]{Definition}

\newtheorem{condition}[theorem]{Condition}

\newcommand{\F}{\mathbb{F}}

\newcommand{\GL}{\mathrm{GL}}

\newcommand{\Oh}{\mathcal{O}} 
\newcommand{\Sel}{\mathrm{Sel}} 
\newcommand{\Gal}{\mathrm{Gal}} 

\newcommand{\Frob}{\mathrm{Frob}}
\newcommand{\Hom}{\mathrm{Hom}}
\newcommand{\et}{\text{\'et}}

\newcommand{\genlegendre}[4]{%
  \genfrac{(}{)}{}{#1}{#3}{#4}%
  \if\relax\detokenize{#2}\relax\else_{\!#2}\fi
}

\DeclareSymbolFont{cyrletters}{OT2}{wncyr}{m}{n}
\DeclareMathSymbol{\Sha}{\mathalpha}{cyrletters}{"58}

%% file: Introduction.tex
\section{Introduction}

\subsection{Main result}
This paper obtains upper bounds on the moments of the number of rational solutions to twist families of some superelliptic curves over global function fields.

Given a global field $K$ and a smooth projective curve $C$ of genus at least two, the Mordell--Lang conjecture states that the number of $K$-rational solutions of $C$ is finite. The proof of the conjecture with non-effective bounds was obtained in groundbreaking studies by Faltings \cite{Faltings83} and Vojta \cite{Vojta91} when $K$ is a number field, and by Manin \cite{Manin63}, Grauert \cite{Grauert65}, and Vojta \cite{Vojta89} when $K$ is a global function field with finite characteristic. The uniform effective upper bounds for the number of $K$-rational solutions of $C$ can be found in a collection of seminal works, such as those of Buium--Voloch \cite{BV96}, Katz--Rabinoff--Zureick-Brown \cite{KRZB16}, and Dimitrov--Gao--Habegger \cite{DGH21}. We state below one of the effective versions of these results.
\begin{theorem}[Uniform Mordell--Lang for curves]
    Let $g \geq 2$ be an integer, and let $K$ be a global field. Given any smooth projective curve $C$ of genus $g$ over $K$ which is not defined over $K^p$, there exists an absolute constant $c := c(g,K) > 0$ such that
    \begin{equation*}
        \# C(K) \leq c^{1 + \mathrm{Rank}_\mathbb{Z} \mathrm{Jac}(C)(K)},
    \end{equation*}
    where $\mathrm{Jac}(C)$ is the Jacobian variety of the curve $C$ over $K$.
\end{theorem}
A key parameter which appears in the statement of the uniform Mordell--Lang conjecture is the Mordell-Weil rank of $\mathrm{Jac}(C)(K)$. It is a natural question to ask to which extent the uniform Mordell--Lang conjecture can be utilized to analyze the distribution of number of $K$-rational points of families of smooth projective curves $\mathcal{C}$ of fixed genus at least two.

Before we give a partial answer to this question, we fix a few notations. Fix a prime number $\ell$. Let $K = \mathbb{F}_q(t)$ be the global function field whose characteristic $p$ is coprime to $2, 3,$ and $\ell$. We denote by $\mu_\ell$ the set of primitive $\ell$-th roots of unity. Throughout this manuscript we assume that $\mu_\ell \subset K$. Given a non-negative integer $n$, we denote by $F_n(\mathbb{F}_q)$ the set of polynomials $f$ over $\mathbb{F}_q$ of degree $n$. We denote by $C$ a non-isotrivial superelliptic curve over $K$ whose Weierstrass equation is
\begin{equation}
    C: y^\ell = F(x),
\end{equation}
where $F$ is a degree $3$ polynomial with coefficients in $K$ but not in $K^p$. The family of superelliptic curves we consider is given by a family of twists of $C$:
\begin{equation}
    \mathcal{C}_n := \left\{ C_f : fy^\ell = F(x) \; | \; f \in F_n(\mathbb{F}_q) \right\}.
\end{equation}

The main result of this paper computes effective probabilistic stratification on the number of $K$-rational points of superelliptic curves $C_f \in \mathcal{C}_n$ as $n$ grows arbitrarily large. Provided below is a simplification of the main result of the manuscript.
\begin{theorem}[A simplification of Theorem \ref{thm:main}] \label{thm:main_simple}
    Suppose the following three conditions hold for the superelliptic curve $C: y^\ell = F(x)$ over $K$.
    \begin{itemize}
        \item $C$ is not defined over $K^p$.
        \item The splitting field of $F(x)$ is an $S_3$-Galois extension over $K$, whose constant field is $\mathbb{F}_q$.
        \item $\mathrm{Jac}(C)$ has a place $v$ of totally split multiplicative reduction.
    \end{itemize}
    Then for sufficiently large $n$, there exists a fixed constant $B(\ell,p)$ depending only on $\ell$ and $p$ such that
    \begin{equation}
        \limsup_{n \to \infty} \frac{\sum_{f \in F_n(\mathbb{F}_q)} \# C_f(K)}{\# F_n(\mathbb{F}_q)} \leq B(\ell,p).
    \end{equation}
    Furthermore, at least $99\%$ of polynomials $f \in F_n(\mathbb{F}_q)$ satisfy $\# C_f(K) \leq (3p)^{4\ell} \cdot \ell!$, and for every positive integer $k \geq 0$, the following inequality holds:
    \begin{equation}
        \limsup_{n \to \infty} \frac{\#\{f \in F_n(\mathbb{F}_q) \; : \; \# C_f(K) > (3p)^{(3+k) \cdot \ell} \cdot \ell!\}}{\# F_n(\mathbb{F}_q)} < 2 \cdot {\ell^{-\frac{(k+1)(k+2)}{2}}}.
    \end{equation}
\end{theorem}

\subsection{Selmer groups of Jacobians}
The main result follows from computing the dimensions of Selmer groups of Jacobians of superelliptic curves in $\mathcal{C}_n$. The computation of distribution of Selmer groups of families of abelian varieties is one of the widely used techniques to understand the distribution of Mordell-Weil ranks of their $K$-rational points. Given an abelian variety $A$ over a global field $K$ and a $K$-rational isogeny $\phi \in \mathrm{End}(A/K)$, the $\phi$-Selmer group of $A$, denoted as $\Sel_\phi(A/K)$, is a finite subset of the first cohomology group $H^1_\et(K, A[\phi])$. In case $\phi$ is a multiplication by $\ell$ map for some prime $\ell$, the group $\Sel_\ell(A/K)$ is a finite dimensional $\mathbb{F}_\ell$ vector space of $H^1_\et(K, A[\ell])$ whose dimension gives an upper bound on the Mordell-Weil rank of $A(K)$. Most notably, the Poonen--Rains heuristics \cite{PR12} predict the following distribution of $\Sel_\ell(A/K)$ for a ``nice'' family of abelian varieties $\mathcal{A}$ endowed with a height function $h: \mathcal{A} \to \mathbb{R}_{\geq 0}$ satisfying the Northcott property (i.e. there are only finitely many isomorphism classes of abelian varieties of bounded height $X$):
\begin{equation}
    \lim_{X \to \infty} \frac{\# \{A \in \mathcal{A} \; | \; h(A) \leq X, \dim_{\mathbb{F}_\ell} \Sel_\ell(A/K) = r\}}{\# \{A \in \mathcal{A} \; | \; h(A) \leq X\}} = \prod_{k=0}^\infty \frac{1}{1+\ell^{-k}} \cdot \prod_{k=1}^r \frac{\ell}{\ell^k-1}.
\end{equation}
As the upper bound on the height $X$ grows arbitrarily large, the average size of $\Sel_\ell(A/K)$ as $A$ ranges over elements of $\mathcal{A}$ converges to $\ell + 1$, and the probability that the dimension of $\Sel_\ell(A/K)$ is even (or odd, respectively) as $A$ ranges over elements of $\mathcal{A}$ converges to $\frac{1}{2}$. The formulation of this heuristic rests upon identifying Selmer groups as intersections of two random maximal isotropic subspaces of an infinite dimensional vector space over finite fields \cite[Theorem 4.14]{PR12}.

There is a wealth of previous studies which compute the moments and distributions of Selmer groups of families of Jacobians of higher genus curves over global fields $K$. When $K$ is a number field, the works by Bhargava--Gross\cite{BG13}, Shankar--Wang \cite{SW18}, and Shankar \cite{Sh19} compute moments of 2-Selmer groups of universal families of hyperelliptic curves. The two works by Yu computes that there are infinitely many twist families of hyperelliptic and superelliptic curves over number fields whose Selmer groups have a given dimension \cite{Yu16}, and that the distribution of 2-Selmer groups of quadratic twist families of hyperelliptic curves conforms to the Poonen-Rains heuristics assuming certain conditions on the local conditions defining the desired Selmer groups \cite{Yu19}. These results, combined with Chabauty's method, may be leveraged to prove strong results on the number of $\mathbb{Q}$-rational points of such hyperelliptic curves, as shown in the groundbreaking work by Poonen--Stoll \cite{PS14} and Smith \cite{Sm22_01, Sm22_02}. For cases where $K$ is a global function field, the work by Thinh computes the first moments of 2-Selmer groups of universal families of even hyperelliptic curves \cite{Thinh23}. The recent groundbreaking work by Ellenberg--Landesman computes the distribution of $n$-Selmer groups of quadratic twist families of principally polarized abelian varieties \cite{EL23}. Both results compute the upper and lower bounds of the distribution, the two values of which converge to the predicted moments conforming to the Poonen-Rains heuristics assuming that the size of the constant field of $K$ grows arbitrarily large.

Among a variety of techniques available to analyze Selmer groups of families of curves in $\mathcal{C}_n$, we focus on the approach of constructing a Markov operator over the state space of non-negative integers, as previously studied in the works of Swinnerton-Dyer \cite{SD08}, Klagsbrun--Mazur--Rubin \cite{KMR14}, Yu \cite{Yu16, Yu19}, and the author \cite{park2022prime}. All these works focus on computing the distribution of Selmer groups of twist families of abelian varieties that share a common $\mathrm{Gal}(\overline{K}/K)$-module that admits a non-degenerate alternating Weil pairing.

\subsection{New contribution}

In this manuscript, we focus on computing the Selmer groups of $\mathrm{Jac}(C_f)$ with respect to the isogeny $1-\zeta_\ell: \mathrm{Jac}(C_f) \to \mathrm{Jac}(C_f)$ induced from the map $\zeta_\ell: C_f \to C_f$ that sends $(x,y)$ to $(x, \zeta_\ell y)$. Unlike the previous works which utilize Markov operators to compute the distribution of Selmer groups, the twist families of Jacobians of superelliptic curves $\mathrm{Jac}(C_f)$ admit a subtle yet crucially different structure: as proven by Poonen and Schaefer \cite[Lemma 10.1]{PS97}, the Weil pairing 
\begin{equation}
    \mathrm{Jac}(C_f)[1-\zeta_\ell] \times \mathrm{Jac}(C_f)[1-\zeta_\ell] \to \mu_\ell
\end{equation}
is a non-degenerate symmetric pairing, not an alternating pairing. (This reference was brought up by Peter Koymans during a discussion with me and Cihan Sabuncu in late July 2026. As a result, several non-trivial corrections to the previous version of the manuscript are made.) This different structure prevents us from using a direct adaptation of the previous works that utilizes Markov machinery to compute the desired distribution of $1-\zeta_\ell$ Selmer groups. Nevertheless, we demonstrate in this mansucript that it is still possible to devise a new Markov model for $1-\zeta_\ell$ Selmer groups of $\mathrm{Jac}(C_f)$. 

\begin{figure}[ht]
    \centering
    \includegraphics[width=0.9\linewidth]{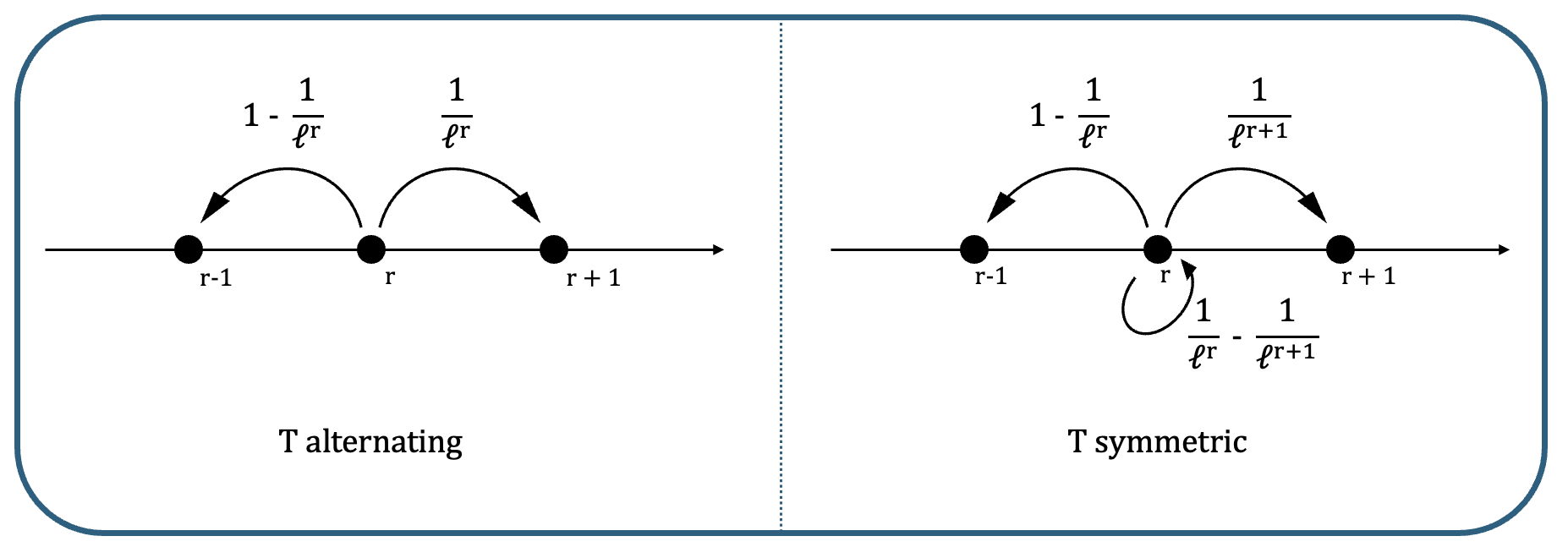}
    \caption{A comparison of Markov operators $M_\ell$ that govern the distribution of dimensions of Selmer groups of an alternating module $T$ and a symmetric module $T$.}
    \label{fig:Markov}
\end{figure}

Figure \ref{fig:Markov} illustrates two types of Markov operators defined over the countable state space of non-negative integers that govern the distribution of Selmer groups associated to a $\mathrm{Gal}(\overline{K}/K)$ module $T$ satisfying some technical conditions (see Condition \ref{condition:admissible} for further details). The left side of the figure denotes the Markov operator studied in aforementioned previous works \cite{SD08, KMR14, Yu19, park2022prime}. These Markov operators are associated to any $\mathrm{Gal}(\overline{K}/K)$-module $T$ satisfying some technical conditions whose Weil pairing is an alternating pairing. The right side of the figure, which is newly proposed in this manuscript, is associated to such modules $T$ whose Weil pairing is a symmetric pairing. The difference in transition probabilities of these two stochastic processes lead to the difference between their stationary distributions, which we will use to characterize the distribution of dimensions of Selmer groups associated to $T$.

There are two key merits for utilizing the Markov machinery: the ability to compute distribution of desired Selmer groups for any size of the constant field of $K$: and the capability to explicitly determine the rate of convergence in terms of $n$. By generalizing the previous work of the author \cite{park2022prime}, we compute the distribution of $1-\zeta_\ell$ Selmer groups of superelliptic curves $C_f \in \mathcal{C}$ with both merits fulfilled. A non-trivial byproduct we obtain is that the distribution of $1-\zeta_\ell$ Selmer groups of Jacobians of curves $C_f \in \mathcal{C}_n$ deviates from the prediction of Poonen-Rains heuristics, but adheres to the random matrix model proposed by Cais--Ellenberg--Zureick-Brown \cite[Proposition 3.5]{Cais_Ellenberg_Zureick-Brown_2013}. Most notably, the density of curves $C_f \in \mathcal{C}$ whose Jacobian has even dimensional $1-\zeta_\ell$ Selmer group does not converge to $1/2$. Furthermore, as $\ell$ grows arbitrarily large, the density of curves $C_f \in \mathcal{C}$ whose Jacobian has its $1-\zeta_\ell$ Selmer group to be of dimension $0$ converges to $1$.
\begin{theorem}[A simplification of Theorem \ref{thm:superelliptic}]
    Assume all the conditions and notations from Theorem \ref{thm:main_simple}. Then for any small enough $\delta > 0$, there exist sufficiently large $n$ and fixed constants $B(C,q) > 0$ and $0 < \alpha(\ell) < 1$ independent of $n$ such that
    \begin{equation}
        \left| \frac{\#\{f \in F_n(\mathbb{F}_q) \; | \; \dim_{\mathbb{F}_\ell} \Sel_{1-\zeta_\ell}(\mathrm{Jac}(C_f)/K) = r\}}{\# F_n(\mathbb{F}_q)} - \prod_{k=1}^\infty \frac{1}{1+\ell^{-k}} \cdot \prod_{i=1}^r \frac{1}{\ell^i-1} \right| < \frac{B(C,q)}{n^{\alpha(\ell) - \delta}}.
    \end{equation}
\end{theorem}
We refer to Theorem \ref{thm:superelliptic} for an explicit description of the constants $B(C,q)$ and $\alpha(\ell)$. To differentiate two cases where the Weil pairing on $T$ is alternating or symmetric (see Definition \ref{def:alternating-symmetric}), we obtain separate distributions of dimensions of Selmer groups of $T$ in Theorems \ref{thm:Galois_module-alternating} and \ref{thm:Galois_module-symmetric}. The first result also incorporates generalizations of remarkable observations made by Klagsbrun--Mazur--Rubin \cite[Corollary 7.10, Theorem 8.2]{KMR13}, Yu \cite[Theorem 1]{Yu16}, and Morgan \cite[Theorem 7.4]{Mor19} where they demonstrate that there are disparities in the dimensions of 2-Selmer groups of quadratic twist families of principally polarized abelian varieties. We then use the second result to obtain Theorem \ref{thm:superelliptic}.

\subsection{Organization}
Section \ref{sec:Selmer} generalizes the main theorem of \cite{park2022prime} to admissible $\Gal(\overline{K}/K)$ modules $T$, a $2$-dimensional vector space over the finite field $\mathbb{F}_\ell$ satisfying a number of mild technical conditions outlined in Condition \ref{condition:admissible}. The highlights of the section are Theorems \ref{thm:Galois_module-alternating} and \ref{thm:Galois_module-symmetric}, where we extend the main result of \cite{KMR14, park2022prime} to any such admissible $\Gal(\overline{K}/K)$ modules $T$ over global function fields $K$. Section \ref{sec:superelliptic} focuses on applying Theorem \ref{thm:Galois_module-symmetric} to obtain Theorem \ref{thm:main} and Theorem \ref{thm:superelliptic}.

\subsection*{Acknowledgements}
I would like to thank Valentin Blomer, Jordan Ellenberg, Daniel Keliher, Chris Keyes, Zev Klagsbrun, Peter Koymans, Cihan Sabuncu, and David Zureick-Brown for helpful discussions and comments. I would like to also thank the Max Planck Institute for Mathematics for its generous support. 

I would like to express my sincerest gratitude to Peter Koymans for pointing out the reference by Poonen and Schaefer \cite{PS97}, which shows that the Weil pairing over $\mathrm{Jac}(C_f)[1-\zeta_\ell]$ is symmetric, during a discussion with me and Cihan Sabuncu on late July 2026. The proof provided in the previous version of the manuscript regrettably was based on an error that the Weil pairing over $\mathrm{Jac}(C_f)[1-\zeta_\ell]$ would induce a non-trivial quadratic form over $H^1_\et(K_v, \mathrm{Jac}(C_f)[1-\zeta_\ell])$, which is not true. As a result, subsequent major revisions have been made, which led to a major expansion of Section \ref{sec:Selmer}.

%% file: Alternating.tex
\section{Distribution of twisted Selmer structures of rank 2 Galois modules}
\label{sec:Selmer}

In this section, we generalize the Markov machinery discussed in \cite{SD08, KMR14, Yu19, park2022prime} to 2-dimensional $\mathrm{Gal}(\overline{K}/K)$ modules $T$ that are endowed with a non-degenerate pairing (either alternating or symmetric) and satisfy some generic technical conditions. We compute the distribution of dimensions of twisted Selmer groups of $T$ using stationary distribution of governing Markov operators.

\subsection{Setup}

We fix a set $T$ equipped with a continuous action of $\Gal(\overline{K}/K)$.
\begin{condition} \label{condition:admissible}
    Fix a prime number $\ell$. We say that a $\Gal(\overline{K}/K)$-module $T$ is \textbf{an admissible} if $T$ satisfies all the six conditions provided below.
    \begin{enumerate}
        \item $T$ is a 2-dimensional $\mathbb{F}_\ell$ vector space.
        \item The constant field of $K(T)$ is equal to $\mathbb{F}_q$.
        \item $T$ is a simple $\Gal(\overline{K}/K)$-module.
        \item $\mathrm{Hom}_{\Gal(\overline{K}/K)}(T,T) = \mathbb{F}_\ell$.
        \item $H^1(\Gal(K(T)/K),T) = 0$.
        \item There exists a place $v$ of $K$ such that $T^{\Gal(\overline{K}_v/K_v\overline{\mathbb{F}}_q)} \cap \mu_{\ell^\infty} = \mu_\ell$.
    \end{enumerate}
\end{condition}

\begin{definition} \label{def:alternating-symmetric}
    We say that an admissible module $T$ is \textbf{alternating} if there exists a $\Gal(\overline{K}/K)$-equivariant non-degenerate alternating pairing $e_T: T \times T \to \mu_\ell$ such that the following is a symmetric pairing:
        \begin{equation*}
            H^1_{\et}(K_v, T) \times H^1_{\et}(K_v, T) \to H^2_{\et}(K_v, T \otimes T) \to \mathbb{F}_\ell,
        \end{equation*}
        where the first map is the cup product, and the second map is the map induced from $e_T$.
        
    Suppose $\ell \geq 3$. We say that an admissible module $T$ is \textbf{symmetric} if $e_T$ is a $\Gal(\overline{K}/K)$-equivariant non-degenerate symmetric pairing such that the analogously defined pairing is alternating:
        \begin{equation*}
            H^1_{\et}(K_v, T) \times H^1_{\et}(K_v, T) \to H^2_{\et}(K_v, T \otimes T) \to \mathbb{F}_\ell.
        \end{equation*}
\end{definition}
The goal of this section is to generalize the results by \cite{KMR14} and \cite{park2022prime} on the distribution of Selmer groups of cyclic twists of a fixed admissible $\Gal(\overline{K}/K)$-module $T$. To state the proposed generalization, we introduce the following notations and results from \cite[Section 5]{KMR14}.

\begin{definition}[{\cite[Definition 5.8, Lemma 5.9]{KMR14}} and  {\cite[Definition 4.2]{park2022prime}}]
Let $T$ be an admissible $\mathrm{Gal}(\overline{K}/K)$- module. We define the following subsets of places $v$ of $K$.
\begin{itemize}
    \item $\Sigma$: a set of places which includes places at which $T$ is ramified.
    \item $\Sigma_T$: the set of places whose elements are exactly places at which $T$ is ramified.
    \item $\Sigma(\mathfrak{d})$: a set of places which includes places in $\Sigma$ and places dividing a square-free product of places $\mathfrak{d}$
    \item $\mathcal{P}_0$: the set of all places $v \not\in \Sigma$ and $\dim_{\mathbb{F}_\ell} T^{\Gal(\overline{K}_v/K_v)} = 0$.
    \item $\mathcal{P}_1$: the set of all places $v \not\in \Sigma$ and $\dim_{\mathbb{F}_\ell} T^{\Gal(\overline{K}_v/K_v)} = 1$.
    \item $\mathcal{P}_2$: the set of all places $v \not\in \Sigma$ and $\dim_{\mathbb{F}_\ell} T^{\Gal(\overline{K}_v/K_v)} = 2$.
\end{itemize}
For each $i$, we denote by $\delta(\mathcal{P}_i)$ the density of places $v$ of $T$ such that $v \in \mathcal{P}_i$.
\end{definition}

\begin{definition}
Given a place $v$ of $K$, we denote by $\mu_{T,v}$ the pairing over $H_\et^1(K_v,T)$ induced from the pairing $e_T: T \times T \to \mu_\ell$ from Condition \ref{condition:admissible}:
\begin{align*}
    \mu_{T,v} &: H_\et^1(K_v,T) \times H_\et^1(K_v,T) \to \mathbb{F}_\ell.
\end{align*}
\end{definition}

\begin{definition}[{\cite[Definition 5.4]{KMR14}}]
    For each place $v$ of $K$, we denote by $\mathcal{H}(\mu_{T,v})$ the set of maximal isotropic subspaces of $H^1_{\et}(K_v,T)$ with respect to the pairing $\mu_{T,v}$. Note that $H^1_\et(\Oh_{K_v},T) \in \mathcal{H}(\mu_{T,v})$. We denote by $\mathcal{H}_{ram}(\mu_{T,v})$ the set of subspaces in $\mathcal{H} \in \mathcal{H}(\mu_{T,v})$ such that $\mathcal{H} \cap H^1_\et(\Oh_{K_v},T) = 0$.
\end{definition}

We first characterize the number of elements in $\mathcal{H}(\mu_{T,v})$. We note that depending on whether $T$ is alternating or symmetric, the cardinality of $\mathcal{H}(\mu_{T,v})$ is different. This difference turns out to play an important role in detecting differences in how Selmer groups associated to the module $T$ are distributed.

\begin{lemma} \label{lemma:H_ram}
    Let $T$ be an admissible $\mathrm{Gal}(\overline{K}/K)$-module. Then the following statements hold for any place $v$ of $K$ that is not in $\Sigma_T$.
    \begin{itemize}
        \item For every $\mathcal{H} \in \mathcal{H}(\mu_{T,v})$ we have $\dim_{\mathbb{F}_\ell} \mathcal{H} = \dim_{\mathbb{F}_\ell} T^{\Gal(\overline{K}_v/K_v)}$.
        \item Suppose $T$ is alternating. If $v \in \mathcal{P}_i$ for $i = 1,2$, then we have $\# \mathcal{H}_{ram}(\mu_{T,v}) = \ell^{i(i-1)/2}$.
        \item Suppose $T$ is symmetric. If $v \in \mathcal{P}_i$ for $i = 1,2$, then we have $\# \mathcal{H}_{ram}(\mu_{T,v}) = \ell^{i(i+1)/2}$.
    \end{itemize}
\end{lemma}
\begin{proof}
    The statement for $T$ an alternating admissible module follows from an adaptation of \cite[Lemma 5.5]{KMR14} to the function field setting. However, we provide the details of the proof that is applicable for $T$ a symmetric admissible module as well.

    By \cite[Lemma B.2.8]{Rubin14} to global function field setting, we have
    \begin{equation*}
        H^1_{\et}(\mathcal{O}_{K_v}, T) \cong \frac{T}{(\mathrm{Frob}_v-1)T}.
    \end{equation*}
    We then proceed as in \cite[Lemma 3.7]{KMR13}. The short exact sequence
    \begin{equation*}
        0 \to T^{\mathrm{Gal}(\overline{K}_v/K_v)} \to T \to T \to \frac{T}{(\mathrm{Frob}_v - 1) T} \to 0,
    \end{equation*}
    where the middle map is the map $\mathrm{Frob}_v-1$, implies $\dim_{\mathbb{F}_\ell} T^{\mathrm{Gal}(\overline{K}_v/K_v)} = \dim_{\mathbb{F}_\ell} H^1_{\et}(\mathcal{O}_{K_v},T)$. By \cite[Theorem I.2.6]{Milne86}, the group $H^1_{\et}(\mathcal{O}_{K_v},T)$ is its own orthogonal complement with respect to the pairing $\mu_{T,v}$. Hence $H^1_{\et}(\mathcal{O}_{K_v},T) \in \mathcal{H}(\mu_{T,v})$. Since every maximal isotropic subspace of $H^1_{\et}(K_v,T)$ have same dimensions, we obtain the first statement of the lemma.

    The second and third statements follow from counting the number of maximal isotropic subspaces of a $2$ and $4$ dimensional $\mathbb{F}_\ell$ vector space with respect to the pairing $\mu_{T,v}$.
\end{proof}
The lemma above gives us a leverage to parametrize subspaces in $\mathcal{H}(\mu_{T,v})$ and $\mathcal{H}_{ram}(\mu_{T,v})$ using order $\ell$ Galois extensions over $K_v$. Given that the number of elements in $\mathcal{H}(\mu_{T,v})$ is different depending on whether $T$ is alternating or symmetric, we give separate definitions of such parametrization.
\begin{definition}[{\cite[Definition 5.10]{KMR14}}]
    Let $T$ be an alternating admissible $\mathrm{Gal}(\overline{K}/K)$-module. A twisting data $\alpha := (\alpha_v)_v$ is a collection of maps $\alpha_v$ for each place $v$ of $K$ defined as follows.
    \begin{itemize}
        \item If $v \in \Sigma$, then $\alpha_v$ is a set-theoretic map
        \begin{equation*}
            \alpha_v: \frac{\mathrm{Hom}(\Gal(\overline{K}_v/K_v), \mu_\ell)}{\mathrm{Aut}(\mu_\ell)} \to \mathcal{H}(\mu_{T,v}).
        \end{equation*}
        \item If $v \in \mathcal{P}_0$, then $\alpha_v$ is the zero map
        \begin{equation*}
            \alpha_v: \frac{\mathrm{Hom}(\Gal(\overline{K}_v/K_v), \mu_\ell)}{\mathrm{Aut}(\mu_\ell)} \to \{0\}.
        \end{equation*}
        \item If $v \in \mathcal{P}_1$, then $\alpha_v$ is the surjection
        \begin{equation*}
            \alpha_v: \frac{\mathrm{Hom}(\Gal(\overline{K}_v/K_v), \mu_\ell)}{\mathrm{Aut}(\mu_\ell)} \to \{H^1_{\et}(\Oh_{K_v}, T)\} \cup \mathcal{H}_{ram}(\mu_{T,v}),
        \end{equation*}
        which maps the unramified $\mathbb{Z}/\ell \mathbb{Z}$ extension to $H^1_{\et}(\Oh_{K_v},T)$, and maps all other ramified extensions to the unique element of $\mathcal{H}_{ram}(\mu_{T,v})$.
        \item If $v \in \mathcal{P}_2$, then $\alpha_v$ is a bijection
        \begin{equation*}
            \alpha_v: \frac{\mathrm{Hom}(\Gal(\overline{K}_v/K_v), \mu_\ell)}{\mathrm{Aut}(\mu_\ell)} \to \{H^1_{\et}(\Oh_{K_v}, T)\} \cup \mathcal{H}_{ram}(\mu_{T,v}).
        \end{equation*}
    \end{itemize}
\end{definition}

\begin{definition}
    Let $T$ be a symmetric admissible $\mathrm{Gal}(\overline{K}/K)$-module. A twisting data $\alpha := (\alpha_v)_v$ is a collection of maps $\alpha_v$ for each place $v$ of $K$ defined as follows.
    \begin{itemize}
        \item If $v \in \Sigma$, then $\alpha_v$ is a set-theoretic map
        \begin{equation*}
            \alpha_v: \frac{\mathrm{Hom}(\Gal(\overline{K}_v/K_v), \mu_\ell)}{\mathrm{Aut}(\mu_\ell)} \to \mathcal{H}(\mu_{T,v}).
        \end{equation*}
        \item If $v \in \mathcal{P}_0$, then $\alpha_v$ is the zero map
        \begin{equation*}
            \alpha_v: \frac{\mathrm{Hom}(\Gal(\overline{K}_v/K_v), \mu_\ell)}{\mathrm{Aut}(\mu_\ell)} \to \{0\}.
        \end{equation*}
        \item If $v \in \mathcal{P}_1$, then $\alpha_v$ is a bijection
        \begin{equation*}
            \alpha_v: \frac{\mathrm{Hom}(\Gal(\overline{K}_v/K_v), \mu_\ell)}{\mathrm{Aut}(\mu_\ell)} \to \{H^1_{\et}(\Oh_{K_v}, T)\} \cup \mathcal{H}_{ram}(\mu_{T,v}),
        \end{equation*}
        which maps the unramified $\mathbb{Z}/\ell \mathbb{Z}$ extension to $H^1_{\et}(\Oh_{K_v},T)$, and maps all other ramified extensions to distinct elements of $\mathcal{H}_{ram}(\mu_{T,v})$.
        \item If $v \in \mathcal{P}_2$, then $\alpha_v$ is an injection
        \begin{equation*}
            \alpha_v: \frac{\mathrm{Hom}(\Gal(\overline{K}_v/K_v), \mu_\ell)}{\mathrm{Aut}(\mu_\ell)} \to \{H^1_{\et}(\Oh_{K_v}, T)\} \cup \mathcal{H}_{ram}(\mu_{T,v}).
        \end{equation*}
    \end{itemize}
\end{definition}

Using the twisting data for $T$ and a choice of a global order-$\ell$ character, we can define the Selmer structure of $T$, the central object of our interest.

\begin{definition}[c.f. {\cite[Definitions 5.6, 5.12]{KMR14}}]
    Given an equivalence class of characters $\overline{\chi} \in \mathrm{Hom}(\Gal(\overline{K}/K), \mu_\ell) / \mathrm{Aut}(\mu_\ell)$ and a choice of a twisting data $\alpha := (\alpha_v)_v$, we denote by $\Sel(T,\alpha : \overline{\chi})$ a finite dimensional $\mathbb{F}_\ell$ vector space defined as
    \begin{equation}
        \Sel(T,\alpha : \overline{\chi}) := \mathrm{Ker} \left( H_\et^1(K,T) \to \prod_v \frac{H_\et^1(K_v,T)}{\alpha_v(\overline{\chi}_v)} \right),
    \end{equation}
    where $\overline{\chi}_v$ is the restriction of $\overline{\chi}$ to an element of $\mathrm{Hom}(\Gal(\overline{K}_v/K_v), \mu_\ell) / \mathrm{Aut}(\mu_\ell)$.
\end{definition}

\begin{definition}
    Let $f \in \mathbb{F}_q[t]$ be a polynomial of degree $n$. We define the following notations associated to $f$:
    \begin{itemize}
        \item We denote by $\overline{\chi}_f \in \mathrm{Hom}(\Gal(\overline{K}/K), \mu_\ell) / \mathrm{Aut}(\mu_\ell)$ an equivalence class of cyclic order $\ell$ characters whose fixed field is $K(\sqrt[\ell]{f})$. We use $\chi_f \in \mathrm{Hom}(\Gal(\overline{K}/K), \mu_\ell)$ to denote a lift of $\overline{\chi}_f$.
        \item We denote by $\overline{\chi}_{f,v} \in \mathrm{Hom}(\Gal(\overline{K}_v/K_v), \mu_\ell) / \mathrm{Aut}(\mu_\ell)$ the restriction of $\overline{\chi}_f$ to the local character over $K_v$. Likewise, $\chi_{f,v} \in \mathrm{Hom}(\Gal(\overline{K}_v/K_v), \mu_\ell)$ is the restriction of $\chi_f$ over $K_v$.
    \end{itemize}
\end{definition}

We state two main results to be proven in this section. We distinguish these results based on whether $T$ is alternating or symmetric. A key feature is the difference in the distribution of Selmer structures of $T$ between the alternating case and the symmetric case.

We first state the alternating case.
\begin{theorem}[Alternating Case]
    \label{thm:Galois_module-alternating}
    Fix a prime number $\ell$. Let $K = \mathbb{F}_q(t)$ be a global function field whose characteristic is coprime to $2$, $3$ and $q \equiv 1 \text{ mod } \ell$. Let $T$ be an alternating admissible $\Gal(\overline{K}/K)$-module. Fix a twisting data $\alpha := (\alpha_v)_v$. Then for any small enough $\delta > 0$, there exist sufficiently large $n$ and a fixed constant $B_{alt} > 0$ independent of $n$ such that
    \begin{align}
        \left| \frac{\# \{f \in F_n(\mathbb{F}_q) \; | \; \dim_{\mathbb{F}_\ell} \Sel(T, \alpha : \overline{\chi}_f) = r \}}{\#F_n(\mathbb{F}_q)} - \rho(r) \cdot \prod_{k=1}^\infty \frac{1}{1+\ell^{-k}} \prod_{k=1}^r \frac{\ell}{\ell^k-1} \right| < \frac{B_{alt}}{n^{\alpha(T) - \delta}},
    \end{align}
    where
    \begin{equation}
        \alpha(T) := \sup_{0 < \rho < 1} \left( \min \left( \rho \log \rho + 1 - \rho, \; -\rho \log \left( 1 - \delta(\mathcal{P}_0) \right), \; - \rho \log \gamma_T \right) \right),
    \end{equation}
    with $\gamma_T$ an explicit constant depending only on $T$ (see Theorem \ref{thm:Markov-alternating} for its construction), and
    \begin{equation}
        \rho(r) := \lim_{n \to \infty} \left(\frac{1}{2} + (-1)^r \cdot \left(\frac{1}{2} - \sum_{\substack{k \geq 0 \\ k \equiv r \text{ mod } 2}} \frac{\#\{f \in F_n(\mathbb{F}_q) \cap \mathcal{P}_0 \; | \; \dim_{\mathbb{F}_\ell} \Sel(T, \alpha : \overline{\chi}_f) = r\}}{\#F_n(\mathbb{F}_q) \cap \mathcal{P}_0} \right) \right).
    \end{equation}
    We note that $\rho(r) = \frac{1}{2}$ for all non-negative integer $r \geq 0$ if and only if $\delta(\mathcal{P}_1) \neq 0$.
\end{theorem}
Next, we state the symmetric case.
\begin{theorem}[Symmetric Case]
    \label{thm:Galois_module-symmetric}
    Fix a prime number $\ell$. Let $K = \mathbb{F}_q(t)$ be a global function field whose characteristic is coprime to $2$, $3$ and $q \equiv 1 \text{ mod } \ell$. Let $T$ be a symmetric admissible $\Gal(\overline{K}/K)$-module such that $\delta(\mathcal{P}_1) \neq 0$. Fix a twisting data $\alpha := (\alpha_v)_v$. Then for any small enough $\delta > 0$, there exist sufficiently large $n$ and a fixed constant $B_{sym} > 0$ independent of $n$ such that
    \begin{align}
        \left| \frac{\# \{f \in F_n(\mathbb{F}_q) \; | \; \dim_{\mathbb{F}_\ell} \Sel(T, \alpha : \overline{\chi}_f) = r \}}{\#F_n(\mathbb{F}_q)} - \prod_{k=1}^\infty \frac{1}{1 + \ell^{-k}} \cdot \prod_{i=1}^j \frac{1}{\ell^{i} - 1} \right| < \frac{B_{sym}}{n^{\alpha(T) - \delta}},
    \end{align}
    where
    \begin{align}
    \begin{split}
        \alpha(T) := \sup_{0 < \rho < 1} \biggl( \min \biggl( & \rho \log \rho + 1 - \rho, \; -\rho \cdot \log(1-\delta(\mathcal{P}_0)), \\
        & - \delta(\mathcal{P}_1) \cdot \rho \cdot \log \gamma_T, \; - \rho \cdot(1-\delta(\mathcal{P}_1)) \cdot \log(1-\delta(\mathcal{P}_1))  \biggr) \biggr),
    \end{split}
    \end{align}
    with $\gamma_T$ an explicit constant depending only on $T$ (see Theorem \ref{thm:Markov-symmetric} for its construction).
\end{theorem}

The proof of the theorem follows from an adaptation of the proof of the main result of \cite{KMR14} and \cite{park2022prime}. We borrow all the ingredients and relevant results from \cite{KMR14} and \cite{park2022prime}, filling in the gaps of the proof if needed. For the interest of potential readers who are aware of the two references, we will adhere to the notations utilized in both references, and introduce new notations when necessary. For the rest of the section, we assume that $T$ is an admissible $\Gal(\overline{K}/K)$-module.

\subsection{Effective Chebotarev density theorem}
We recall the two types of effective Chebotarev density theorems from Section 3.1 of \cite{park2022prime}.
\begin{theorem}[{\cite[Proposition 6.4.8]{FJ08}} and {\cite[Theorem 3.1]{park2022prime}}] \label{theorem:effective_chebotarev}
Let $L/K$ be a Galois extension of global function fields over $\F_q(t)$. Pick a conjugacy class $C \subset G = \Gal(L/K)$. Denote by $g_L$ and $g_K$ genera of $L$ and $K$. If the constant fields of $L$ and $K$ are both equal to $\F_q$, then
\begin{align*}
    & \left| \# \{ v \; \text{a place over} \; K \; | \; \Frob_v \in C, \; \dim_{\F_q} (\Oh_K/v) = n \} - \frac{|C|}{|G|} \frac{q^n}{n} \right| \\
    & < \frac{2|C|}{n|G|} \left[ (|G| + g_L) q^{\frac{n}{2}} + |G|(2g_K + 1)q^{\frac{n}{4}} + (|G| + g_L) \right].
\end{align*}
\end{theorem}

\begin{corollary}[{\cite[Corollary 3.2]{park2022prime}}]
\label{corollary:effective_chebotarev}
Assume all conditions and notations from Theorem \ref{theorem:effective_chebotarev}. Pick two non-empty subsets $S, S' \subset G = \Gal(L/K)$ stable under conjugation. Suppose that the size of the constant field $q$ satisfies
    \begin{equation*}
        q^{\frac{n}{2}} - q^{\frac{n}{4}} > 2(|G| + g_L + 2g_K).
    \end{equation*}
Then the following inequality holds.
\begin{align*}
    & \left| \frac{\{ v, \; \text{a place over} \; K \; | \; \Frob_v \in S, \; \dim_{\F_q} (\Oh_K/v) = n \}}{\{ v, \; \text{a place over} \; K \; | \; \Frob_v \in S', \; \dim_{\F_q} (\Oh_K/v) = n \}} - \frac{|S|}{|S'|} \right| \\
    & < 4\frac{|S|}{|S'|} (|G| + g_L + 2g_K) \left[ \frac{1}{q^\frac{n}{2} - q^\frac{n}{4} - 2(|G| + g_L + 2g_K) } \right].
\end{align*}
In particular, if $n \geq 2 \cdot (\log 8 + \log(|G| + g_L + 2g_K)) \cdot (\log q)^{-1}$, then
\begin{equation*}
    \left| \frac{\{ v, \; \text{a place over} \; K \; | \; \Frob_v \in S, \; \dim_{\F_q} (K/v) = n \}}{\{ v, \; \text{a place over} \; K \; | \; \Frob_v \in S', \; \dim_{\F_q} (K/v) = n \}} - \frac{|S|}{|S'|} \right| <  16 \frac{|S|}{|S'|} (|G| + g_L + 2g_K) q^{-\frac{n}{2}}.
\end{equation*}
\end{corollary}

\subsection{Irreducible factors}

Throughout this section, we define two constants $m_{n,q}$ and $\mathfrak{n}$ as
\begin{align}
    \begin{split}
        m_{n,q} &:= \log n + \log \log q, \\
        \mathfrak{n} &:= \frac{4m_{n,q}^2}{\log q}.
    \end{split}
\end{align}
We denote by $g_T$ the genus of the geometric Galois extension $K(T)/K$.

\begin{definition}[{\cite[Definition 4.11]{park2022prime}}]
We denote by $\hat{F}_{(n,N),(w,w')}(\F_q)$ the set of polynomials $f \in \mathbb{F}_q[t]$ of degree $n$ which satisfies the following conditions:
\begin{itemize}
    \item The number of distinct irreducible factors of $f$ is equal to $w$.
    \item The number of distinct irreducible factors of $f$ of degrees greater than $\mathfrak{n}$ is equal to $w'$.
    \item The degree of products of all irreducible factors of $f$ (including multiplicities) is equal to $N$.
    \item There exists at least one irreducible factor of $f$ of degree greater than $\mathfrak{n}$ that lies in $\mathcal{P}_0$.
\end{itemize}
\end{definition}

With appropriate choices of $w, w'$, and $N$, we can approximate the set $F_n(\mathbb{F}_q)$ with disjoint unions of $\hat{F}_{(n,N),(w,w')}(\F_q)$, with explicit error term computed as below.
\begin{proposition}[{\cite[Proposition 4.14]{park2022prime}}] \label{proposition:fan_approximation}
    Let $\rho \in (0,1)$ be a positive number. Suppose $n$ is a positive integer such that $m_{n,q} > \max \{e^{e^e}, \log6 + \log(\ell^3 + g_{T})\}$. Let $\epsilon = \frac{1}{\log \log m _{n,q}}$. Then
    \begin{align}
    \begin{split}
        & \; \; \; \; \# F_n(\F_q) - \sum_{w = \rho m_{n,q}}^{2 m_{n,q}} \sum_{w' = (1-\epsilon)w}^w \sum_{N = w'\mathfrak{n}}^n \# \hat{F}_{(n,N),(w,w')}(\F_q) \\
        &\leq 4 \cdot q^n \cdot \max \left(n^{-\rho \log \rho - 1 + \rho}, 3 m_{n,q}^2 \cdot \left( 1 - \delta(\mathcal{P}_0) \right)^{(1-\epsilon)\rho m_{n,q}}  \right).
    \end{split}
    \end{align}
\end{proposition}
The first component of the error term originates from large and moderate deviation principle for distinct irreducible factors of polynomials over $\mathbb{F}_q$ \cite{FWY20}, whereas the second component of the error term originates from computing the probability that all $w'$ irreducible factors of $f \in \hat{F}_{(n,N),(w,w')}(\F_q)$ do not lie in $\mathcal{P}_0$.

To extend the previous works to computing distribution of Selmer groups of symmetric admissible modules, we introduce a number of notations from \cite[Sections 4, 5]{park2022prime} which feasibly decomposes the set $\hat{F}_{(n,N),(w,w')}(\mathbb{F}_q)$.

\begin{definition}[{\cite[Definition 4.5, 4.7]{park2022prime}}]
    Fix two positive integers $n$ and $w$. Let $\lambda_n$ be a set of $3n^2$ many $4$-tuples of non-negative integers
    \begin{equation*}
        \lambda_n := \{(\lambda_{i,j,k},i,j,k)\}_{1 \leq i \leq n, 1 \leq j \leq n, 0 \leq k \leq 2}.
    \end{equation*}
    We say that $\lambda_n$ is a splitting partition with respect to $(n,w)$ if it satisfies
    \begin{equation*}
        \sum_{i,j,k} \lambda_{i,j,k} \cdot i \cdot j = n, \; \sum_{i,j,k} \lambda_{i,j,k} = w.
    \end{equation*}
    We say that a polynomial $f$ over $\mathbb{F}_q$ admits a splitting partition with respect to $\lambda_n$ if it satisfies the following conditions:
    \begin{itemize}
        \item We have $\deg f = n$ and $\omega(f) = w$.
        \item For each $1 \leq i \leq n$, there are $\lambda_{i,j,k}$ many irreducible factors $g$ of $f$ of degree $i$ and multiplicity $j$ such that $g \in \mathcal{P}_k$.
    \end{itemize}
\end{definition}

\begin{definition}[{\cite[Definitions 4.8, 4.9]{park2022prime}}]
    Given two integers $n, w$, we define the following subsets of splitting partitions.
    \begin{itemize}
        \item $\Lambda_{n,w} := \{\lambda_n : \lambda_n \text{ is a splitting partition}\}.$
        \item $\Lambda_{n,w}^{for} := \{\lambda_n \in \Lambda_{n,w}: \lambda_{i,j,k} = 0 \text{ whenever } i > \mathfrak{n}\}.$
        \item $\Lambda_{n,w}^{la} := \{\lambda_n \in \Lambda_{n,w} : \lambda_{i,j,k} = 0 \; \forall i \leq \mathfrak{n}, \text{ and } \exists \; i > \mathfrak{n}, j \not\equiv 0 \pmod \ell \text{ s.t. } \lambda_{i,j,0} \neq 0\}.$
    \end{itemize}
\end{definition}

\begin{definition}[{\cite[Definitions 4.10, 4.11]{park2022prime}}] \label{defn:Fhat_decomposition}
    Given a polynomial $f \in \mathbb{F}_q[t]$ and an irreducible polynomial $g \in \mathbb{F}_q[t]$, we denote by $v_g(f)$ the multiplicity of $g$ dividing $f$.

    Given a polynomial $f \in \mathbb{F}_q[t]$, we denote by $f_*$ and $f^*$ as
    \begin{equation}
        f_* := \prod_{i=1}^{\mathfrak{n}}\prod_{g \in \mathcal{P}(i)} g^{v_g(f)}, \hspace{15pt} f^* := \prod_{i=\mathfrak{n} + 1}^{\deg f} \prod_{g \in \mathcal{P}(i)} g^{v_g(f)}.
    \end{equation}
    By construction, we have $f = f_* f^*$.

    Given two partitions $\lambda \in \Lambda_{N,w'}^{la}$ and $\eta \in \Lambda_{n-N,w-w'}^{for}$, we define 
    \begin{equation}
        F_{(n,N),(w,w')}^{(\lambda,\eta)}(\mathbb{F}_q) := \{f \in F_{n}(\mathbb{F}_q) : f_* \text{ admits } \eta, \text{ and } f^* \text{ admits } \lambda\}.
    \end{equation}
    One can check that 
    \begin{equation*}
        \hat{F}_{(n,N),(w,w')}(\mathbb{F}_q) = \bigsqcup_{\lambda \in \Lambda_{N,w'}^{la}} \bigsqcup_{\eta \in \Lambda_{n-N,w-w'}^{for}} F_{(n,N),(w,w')}^{(\lambda,\eta)}(\mathbb{F}_q).
    \end{equation*}
\end{definition}

Using the notion of unordered configuration spaces, we describe an alternate way to further decompose the set $F_{(n,N),(w,w')}^{(\lambda,\eta)}(\mathbb{F}_q)$.

\begin{definition}[{\cite[p. 3301]{park2022prime}}]
    Given a set $X$ and a positive integer $n$, we denote by 
    \begin{equation*}
        \mathrm{PConf}_n(X) := \{(x_1, \cdots, x_n) \in X^{\oplus n} : x_i \neq x_j \forall 1 \leq i < j \leq n\}.
    \end{equation*}
    The set above admits a natural permutation action by the symmetric group in $n$ letters $S_n$. This permits us to define
    \begin{equation*}
        \mathrm{Conf}_n(X) := \mathrm{PConf}_n(X)/S_n.
    \end{equation*}
    With these notations, given $\lambda \in \Lambda_{N,w'}^{la}$ and $\eta \in \Lambda_{n-N,w-w'}^{for}$, we have
    \begin{equation*}
        F_{(n,N),(w,w')}^{(\lambda,\eta)}(\mathbb{F}_q) = \prod_{i,j,k} \mathrm{Conf}_{\lambda_{i,j,k}}(\mathcal{P}_k(i)) \times  \prod_{\hat{i},\hat{j},\hat{k}} \mathrm{Conf}_{\eta_{\hat{i},\hat{j},\hat{k}}}(\mathcal{P}_{\hat{k}}(\hat{i})).
    \end{equation*}
\end{definition}

Next, we recall the notations which will be useful in constructing the set $F_{(n,N),(w,w')}^{(\lambda,\eta)}(\mathbb{F}_q)$ by iteratively adding the components corresponding to unordered configuration space.

\begin{definition}[{\cite[Definitions 5.7, 5.9, 5.10]{park2022prime}}]
    Given $f \in F_n(\mathbb{F}_q)$, we define the following three square-free polynomials:
    \begin{equation*}
        \overline{f} := \prod_{g \mid f, g \in \mathcal{P}_1 \cup \mathcal{P}_2} g, \hspace{15pt} \overline{f}_* := \prod_{g \mid f_*, g \in \mathcal{P}_1 \cup \mathcal{P}_2} g, \hspace{15pt} \overline{f}^* := \prod_{g \mid f^*, g \in \mathcal{P}_1 \cup \mathcal{P}_2} g.
    \end{equation*}
    We denote by $\Sigma_f$ the set of places
    \begin{equation*}
        \Sigma_f := \Sigma_T \cup \{v \in \mathcal{P}: v \mid f_*\}.
    \end{equation*}
    Finally, we use the abbreviation $\Omega_1$ and $\Omega_{\overline{f}^*}$ to denote the following sets of finite Cartesian products of local characters:
    \begin{align*}
        \Omega_1 &:= \prod_{v \in \Sigma_f} \mathrm{Hom}(\mathrm{Gal}(\overline{K}/K), \mu_\ell), \\
        \Omega_{\overline{f}^*} &:= \prod_{v \in \Sigma_f} \mathrm{Hom}(\mathrm{Gal}(\overline{K}/K), \mu_\ell) \times \prod_{v \mid \overline{f}^*} \mathrm{Hom}_{ram}(\mathrm{Gal}(\overline{K}_v/K_v), \mu_\ell).
    \end{align*}
\end{definition}

From the construction of $F_{(n,N),(w,w')}^{(\lambda,\eta)}(\mathbb{F}_q)$, we required that any polynomial $f$ in this set must contain at least one polynomial factor of large enough degree in $\mathcal{P}_0$. The following definition gives a canonical choice of such irreducible factors of $f$.

\begin{definition}[{\cite[Definition 5.8, pp. 3301 -- 3302]{park2022prime}}]
    Given positive integers $n > N$ and $w > w'$, we choose $\lambda \in \Lambda_{N,w'}^{la}$ and $\eta \in \Lambda_{n-N,w-w'}^{for}$. Let $f \in F_{(n,N),(w,w')}^{(\lambda, \eta)}(\mathbb{F}_q)$.
    \begin{itemize}
        \item An auxiliary place of $f$ is an irreducible factor $g \in \mathcal{P}_0$ of $f$ of maximal degree. We denote by $d_a$ the maximal degree of such factor (or factors) of $f$.
        \item The auxiliary factor of $f$, denoted as $f_a$, is
        \begin{equation*}
            f_a := \prod_{g \mid f, g \in \mathcal{P}_0(d_a)} g^{v_g(f)}.
        \end{equation*}
        The quantity $d_{a^*}$ is given by $\deg(f_a)$. Note that $d_{a^*} = d_a \cdot \sum_{j=1}^{\ell - 1} \lambda_{d_a,j,0}$.
        \item Let $\phi_{d_a}$ be the natural projection map
    \begin{equation}
        \phi_{d_a}: F_{(n,N),(w,w')}^{(\lambda,\eta)}(\mathbb{F}_q) \to \left[ \prod_{\substack{i,j,k \\ (i,k) \neq (d_a,0)}} \mathrm{Conf}_{\lambda_{i,j,k}}(\mathcal{P}_k(i)) \right] \times \left[ \prod_{\substack{\hat{i},\hat{j},\hat{k}}} \mathrm{Conf}_{\eta_{\hat{i},\hat{j},\hat{k}}}(\mathcal{P}_{\hat{k}}(\hat{i})) \right],
    \end{equation}
    which forgets the auxiliary factor $f_a$ of $f$.
    \end{itemize}
\end{definition}

The motivation for defining the auxiliary factor of $f$ rests upon obtaining an equidistribution over Cartesian products of local characters ramified at places dividing $\overline{f}$ by restricting global characters. The following statement, which holds for any admissible module $T$, gives the desired equidistribution results.

\begin{lemma}[{\cite[Corollary 5.11]{park2022prime}}] \label{lemma:equidistribution}
    Fix $\lambda \in \Lambda_{N,w'}^{la}$ and $\eta \in \Lambda_{n-N,w-w'}^{for}$. Let $T$ be an admissible $\mathrm{Gal}(\overline{K}/K)$-module. Given $h \in F_{n-d_{a^*}}(\mathbb{F}_q)$, suppose the set
    \begin{equation*}
        \phi_{d_a}^{-1}(h) \subset F_{(n,N),(w,w')}^{(\lambda, \eta)}(\mathbb{F}_q)
    \end{equation*}
    is non-empty and $w \leq 2 m_{n,q}$. Denote by $h_1, h_2, \cdots, h_{w(h)}$ the irreducible factors of $h$. We let
    \begin{equation}
        \hat{w}(h) := \begin{cases}
            w(h) &\text{ if } K(\sqrt[\ell]{h_1}, \cdots, \sqrt[\ell]{h_{w(h)}}) \cap K(T) = K \text{ or } \ell \geq 5, \\
            w(h) - 1 &\text{ if } K(\sqrt[\ell]{h_1}, \cdots, \sqrt[\ell]{h_{w(h)}}) \cap K(T) \neq K \text{ and } \ell \leq 5.
        \end{cases} 
    \end{equation}
    Let $\chi := (\chi_v)_v \in \Omega_{\overline{f}^*}$ be any product of local characters, such that at every $v \mid h$ the component $\chi_v$ is ramified, and otherwise the component $\chi_v$ is unramified and there exists $f \in \phi_{d_a}^{-1}(h)$ such that $\mathrm{Ker}(\chi_{f,v}) = \mathrm{Ker}(\chi_v)$. Then there exists an explicit constant $C_{T,\ell} > 0$ depending on $T$ and $\ell$ such that
    \begin{equation}
        \left| \frac{\#\{f \in \phi^{-1}_{d_a}(h) : \mathrm{Ker}(\chi_{f,v}) = \mathrm{Ker}(\chi_v) \forall v \in \Sigma(\overline{h})\}}{\# \phi^{-1}_{d_a}(h)} - \frac{1}{\ell^{\hat{w}(h)}} \right| < \hat{C}_{T,\ell} \cdot (n \log q)^{-2 m_{n,q} + 2 \log \ell}.
    \end{equation}
\end{lemma}
\begin{proof}
    The proof follows from a line-by-line adaptation of \cite[Corollary 4.17, Corollary 5.11]{park2022prime}. We determine the explicit constant to be $\hat{C}_{T,\ell} := 12 g_{T} + 12 \ell^3 - \ell - 2$, where $g_T$ is the genus of the Galois extension $K(T)/K$.
\end{proof}

Lastly, we define a projection map which accesses each component of unordered configuration spaces of $F_{(n,N),(w,w')}^{(\lambda, \eta)}(\mathbb{F}_q)$.
\begin{definition}[{\cite[Definition 5.12]{park2022prime}}] \label{def:phi_d*j*k*}
    Fix $\lambda \in \Lambda_{N,w'}^{la}$ and $\eta \in \Lambda_{n-N,w-w'}^{for}$. Fix $1 \leq j^* \leq \ell - 1$ and $0 \leq k^* \leq 2$. Let $d^*$ be an integer such that $d^* \neq d_a$ and $\lambda_{d^*,j^*,k^*} \neq 0$.

    We denote by $\phi_{d^*,j^*,k^*}$ the projection map
    \begin{equation}
        \phi_{d^*,j^*,k^*}: F_{(n,N),(w,w')}^{(\lambda,\eta)}(\mathbb{F}_q) \to \left[ \prod_{\substack{i,j,k \\ (i,k) \neq (d_a,0) \\ (i,j,k) \neq (d^*,j^*,k^*)}} \mathrm{Conf}_{\lambda_{i,j,k}}(\mathcal{P}_k(i)) \right] \times \left[ \prod_{\hat{i},\hat{j},\hat{k}} \mathrm{Conf}_{\eta_{\hat{i},\hat{j},\hat{k}}}(\mathcal{P}_{\hat{k}}(\hat{i})) \right]
    \end{equation}
    which forgets the auxiliary factor $f_a$ and any irreducible factors of $f$ lying in the set $\mathrm{Conf}_{\lambda_{d^*,j^*,k^*}}(\mathcal{P}_{k^*}(d^*))$.
\end{definition}

\subsection{Markov chains}

A key insight from \cite{SD08}, \cite{KMR14}, and \cite{park2022prime} is that there is a Markov operator over the countable state space of non-negative integers whose stationary distribution determines the distribution of Selmer groups of twist families of elliptic curves over global fields. 

\subsubsection{Alternating case}
Let $T$ be an alternating admissible module. We recall the relevant facts about the Markov operator of our interest.

\begin{definition}[{\cite[Definition 6.1]{park2022prime}}] \label{defn:Markov}
    Let $M_\ell := [\ell_{r,s}]$ be the operator over the state space of non-negative integers $\mathbb{Z}_{\geq 0}$ defined as
    \begin{equation*}
        \ell_{r,s} := \begin{cases}
            1 - \ell^{-r} &\text{ if } s = r-1 \geq 0 \\
            \ell^{-r} &\text{ if } s = r+1 \\
            0 &\text{ else}.
        \end{cases}
    \end{equation*}
\end{definition}

\begin{definition}
    Given an alternating admissible $\Gal(\overline{K}/K)$-module $T$, we denote by $M_T$ the Markov operator over the state space of non-negative integers $\mathbb{Z}_{\geq 0}$ defined as
    \begin{equation}
        M_T := \delta(\mathcal{P}_0) \cdot I + \delta(\mathcal{P}_1) \cdot M_\ell + \delta(\mathcal{P}_2) \cdot M_\ell^2,
    \end{equation}
    where $I$ is the identity operator.
\end{definition}
We note that in the case where $\ell = 2$, a necessary condition for a module $T$ to be alternating and admissible is to have $\Gal(K(T)/K) \cong S_3$. In such a case, the Markov operator $M_T$ is given by
\begin{equation}
    M_T := \frac{1}{3} \cdot I + \frac{1}{2} \cdot M_2 + \frac{1}{6} \cdot M_2^2.
\end{equation}
We note that $M_T$ is an irreducible, aperiodic, and positive-recurrent Markov operator if and only if $\delta(\mathcal{P}_1) \neq 0$. Note that if $\delta(\mathcal{P}_1) = 0$, then $M_T$ has two non-communicating states: the set of odd integers and the set of even integers. Nevertheless, in each non-communicating states, the operator $M_T$ is irreducible, aperiodic, and positive-recurrent. To determine stationary distributions of $M_T$ based on the value of $\delta(\mathcal{P}_1)$, we recall the notion of the parity of a probability distribution.
\begin{definition}[{\cite[Definition 6.3]{park2022prime}}]
    Let $\mu:\mathbb{Z}_{\geq 0} \to [0,1]$ be a probability distribution. The parity of $\mu$, denoted as $\rho(\mu)$, is the sum of probabilities supported over odd integer states, i.e. $\rho(\mu) := \sum_{n \equiv 1 \text{ mod } 2} \mu(n)$.
\end{definition}
The following theorem elucidates the geometric rate of convergence of the Markov operator $M_T$ to a stationary distribution given the parity of the initial probability distribution $\mu$.
\begin{theorem} \label{thm:Markov-alternating}
    Suppose $T$ is an alternating admissible $\Gal(\overline{K}/K)$-module. Given a probability distribution $\mu: \mathbb{Z}_{\geq 0} \to [0,1]$, we denote by $\mathbb{E}[\ell^\mu] := \sum_{z \in \mathbb{Z}_{\geq 0}} \ell^z \cdot \mu(z)$.
    \begin{enumerate}
        \item Suppose $\delta(\mathcal{P}_1) \neq 0$. Then for every $n \in \mathbb{N}$ and any probability distribution $\mu$, there exist fixed constants $c_T > 0$ and $0 < \gamma_T < 1$ depending only on $T$ such that
        \begin{equation}
            \sup_{r \in \mathbb{Z}_{\geq 0}} \left| M_T^n \mu (r) - \prod_{j=0}^\infty \frac{1}{1+\ell^{-j}} \prod_{j=1}^r \frac{\ell}{\ell^j - 1} \right| < c_T \gamma_T^n (\mathbb{E}[\ell^\mu] + 1).
        \end{equation}
        \item Suppose $\delta(\mathcal{P}_1) = 0$. Then for every $n \in \mathbb{N}$ and any probability distribution $\mu$, there exist fixed constants $c_T > 0$ and $0 < \gamma_T < 1$ depending only on $T$ such that
        \begin{equation}
            \sup_{r \in \mathbb{Z}_{\geq 0}} \left| M_T^n \mu (z) - \left( \frac{1}{2} + (-1)^r \cdot \left(\frac{1}{2} - \rho(\mu) \right) \right)\prod_{j=1}^\infty \frac{1}{1+\ell^{-j}} \prod_{j=1}^r \frac{\ell}{\ell^j - 1} \right| < c_T \gamma_T^n (\mathbb{E}[\ell^\mu] + 1).
        \end{equation}
    \end{enumerate}
    In particular, the constant $\gamma_T$ is the geometric rate of convergence of the Markov operator $M_T$ to its stationary distribution.
\end{theorem}
\begin{proof}
    The theorem follows from \cite[Proposition 2.4]{KMR14}, \cite[Theorem 15.0.1]{MT93}, and \cite[Corollary 6.7, Lemma 6.8]{park2022prime}. For the case where $\delta(\mathcal{P}_1) = 0$, we take $\gamma_T$ to be the maximum of geometric rates of convergence over each non-communicating state of $M_T$.
\end{proof}

\subsubsection{Symmetric case}

In Lemma \ref{lemma:H_ram}, there is a difference in cardinality of $\mathcal{H}_{ram}(\mu_{T,v})$ depending on whether $T$ is alternating or symmetric. This difference is closely associated to difference in the auxiliary Markov operator that we will construct for a symmetric admissible module $T$. Given that the previously mentioned works \cite{SD08, KMR14, park2022prime} utilize the Lagrangian Markov operator where the Galois module $T$ is alternating, we will consider a different operator instead:

\begin{definition} \label{defn:Markov-symmetric}
    Let $M_\ell := [\ell_{r,s}]$ be the operator over the state space of non-negative integers $\mathbb{Z}_{\geq 0}$ defined as
    \begin{equation*}
        \ell_{r,s} := \begin{cases}
            1 - \ell^{-r} &\text{ if } s = r-1 \geq 0, \\
            \ell^{-r} \cdot (1 - \ell^{-1}) &\text{ if } s = r, \\
            \ell^{-r - 1} &\text{ if } s = r+1, \\
            0 &\text{ else}.
        \end{cases}
    \end{equation*}
\end{definition}

We note that $M_\ell$ is an irreducible, aperiodic, and positive-recurrent Markov operator. Hence, there exists a unique stationary distribution $\pi: \mathbb{Z}_{\geq 0} \to [0,1]$. As Theorem \ref{thm:Galois_module-symmetric} suggests, it turns out that the stationary distribution is given by
\begin{equation}
    \pi(j) := \prod_{k=1}^\infty \frac{1}{1 + \ell^{-k}} \cdot \prod_{i=1}^j \frac{1}{\ell^i - 1},
\end{equation}
where if $j = 0$, we regard the second product term to be $1$.

The following theorem elucidates the geometric rate of convergence of the Markov operator $M_\ell$ to the stationary distribution $\pi$.
\begin{theorem} \label{thm:Markov-symmetric}
    Suppose $T$ is a symmetric admissible $\Gal(\overline{K}/K)$-module. Given a probability distribution $\mu: \mathbb{Z}_{\geq 0} \to [0,1]$, we denote by $\mathbb{E}[\ell^\mu] := \sum_{z \in \mathbb{Z}_{\geq 0}} \ell^z \cdot \mu(z)$. Then for every $n \in \mathbb{N}$ and any probability distribution $\mu$, there exist fixed constants $c_T > 0$ and $0 < \gamma_T < 1$ such that
    \begin{equation}
        \sup_{j \in \mathbb{Z}_{\geq 0}} \biggl| M_\ell^n \mu (j) - \prod_{k=1}^\infty \frac{1}{1 + \ell^{-k}} \cdot \prod_{i=1}^j \frac{1}{\ell^i - 1} \biggr| < c_T \cdot \gamma_T^n \cdot \mathbb{E}[\ell^\mu].
    \end{equation}
\end{theorem}
\begin{proof}
    Given a non-negative integer $z \in \mathbb{Z}$, let $\mu_z:\mathbb{Z}_{\geq 0} \to [0,1]$ be a probability distribution defined as $\mu_z(z) = 1$, and $\mu_z(j) = 0$ if $j \neq z$.
    
    The theorem follows from \cite[Theorem 15.0.1]{MT93} and \cite[Theorem 1.1]{Bax05}. We will follow the treatment appearing in the latter reference, and check the following three conditions.
    \begin{itemize}
        \item \textit{(A1) Minorization condition}: We pick $C = \{0,1\}$, $\nu: \mathbb{Z}_{\geq 0} \to [0,1]$ a probability measure such that $\nu(\{0\}) = 1 - \ell^{-1}$ and $\nu(\{1\}) = \ell^{-1},$ and $\tilde{\beta} = 1 - \ell^{-1}$. Then for any subset $A \subset \mathbb{Z}_{\geq 0}$, we have
        \begin{equation*}
            \sum_{z \in A} (M_\ell \mu_0)(z) \geq \tilde{\beta} \cdot \nu(A), \hspace{15pt} \sum_{z \in A} (M_\ell \mu_1)(z) \geq \tilde{\beta} \cdot \nu(A).
        \end{equation*}
        \item \textit{(A2) Drift Condition}: We have
    \begin{align*}
        \mathbb{E}[\ell^{M_\ell \mu_z}] &= \ell^{z-1} \cdot (1 - \ell^{-z}) + \ell^z \cdot \ell^{-z} \cdot (1 - \ell^{-1}) + \ell^{z+1} \cdot \ell^{-z-1} \\
        &= \ell^{-1} \cdot \ell^z + 2 - 2 \cdot \ell^{-1} \\
        &= \ell^{-1} \cdot \mathbb{E}[\ell^{\mu_z}] + (2 - 2 \cdot \ell^{-1}).
    \end{align*}
    The condition then holds by setting $V(x) := \ell^x$, $K = 3 - 2 \cdot \ell^{-1}$, and $\gamma = \frac{3}{4}$.
        \item \textit{(A3) Strong aperiodicity condition}: We choose $\beta = \tilde{\beta}$ because $\nu(C) = 1$.
    \end{itemize}
    Hence, we have for all $x \in \mathbb{Z}_{\geq 0}$ and $n \geq 0$,
    \begin{equation*}
        \sup_{|g| \leq V} \left| (M_\ell^n g)(x) - \int g d \pi \right| < c_T \cdot V(x) \cdot \gamma_T^{-n},
    \end{equation*}
    where the supremum is taken over all measurable functions $g: \mathbb{Z}_{\geq 0} \to \mathbb{R}$ such that $|g(x)| \leq V(x)$ for all $v \in \mathbb{Z}_{\geq 0}$. Since $V(x) \geq 1$, we can always choose $g$ to be a probability distribution function over $\mathbb{Z}_{\geq 0}$ to conclude.
\end{proof}

\subsection{Local Selmer groups}

To relate the distribution of Selmer groups $\Sel(T, \alpha : \overline{\chi})$ with stationary behavior of the Markov operator $M_T$, we utilize the notion of local Selmer structures associated to the module $T$.
\begin{definition}[{\cite[Section 5]{KMR14}} and {\cite[Definition 5.1]{park2022prime}}]
Given a set $\Sigma$ of places of $K$, we denote by $\mathfrak{d}$ a square-free product of places coprime to elements in $\Sigma$.
\begin{itemize}
    \item $\Omega_T$: the set of finite Cartesian products of local characters
    \begin{equation*}
        (\chi_v)_v \in \Omega_T := \prod_{v \in \Sigma_T} \mathrm{Hom}(\Gal(\overline{K}_v/K_v), \mu_\ell)
    \end{equation*}
    \item $\Omega_\mathfrak{d}$: the set of finite Cartesian products of local characters
    \begin{equation*}
        (\chi_v)_v \in \Omega_\mathfrak{d} := \prod_{v \in \Sigma(\mathfrak{d})} \mathrm{Hom}(\Gal(\overline{K}_v/K_v), \mu_\ell)
    \end{equation*}
    such that the local character at places $v \mid \mathfrak{d}$ is ramified.
    \item $\Omega_{\chi, \mathfrak{v}}$: the set of finite Cartesian products of local characters $\chi' \in \Omega_{\mathfrak{d}\mathfrak{v}}$ given an element $\chi \in \Omega_\mathfrak{d}$ that satisfies the following two conditions:
    \begin{itemize}
        \item For any $v \in \Sigma(\mathfrak{d})$, we have $\chi'_v = \chi_v$.
        \item The component $\chi'_\mathfrak{v}$ is ramified at $\mathfrak{v}$.
    \end{itemize}
\end{itemize}
\end{definition}

\begin{definition}[{\cite[Definition 5.2]{park2022prime}}]
    Given an element $(\chi_v)_v \in \Omega_\mathfrak{d}$ and a choice of a twisting data $\alpha$, the local Selmer group of $T$ associated to $\alpha$ and $(\chi_v)_v$ is defined as
    \begin{equation}
        \Sel(T, \alpha : (\chi_v)_v) := \mathrm{Ker} \left( H_\et^1(K, T) \to \prod_v \frac{H^1_\et(K_v, T)}{\alpha_v(\chi_v)}\right).
    \end{equation}
    The dimension of $\Sel(T, \alpha : (\chi_v)_v)$, as an $\mathbb{F}_\ell$-vector space, is denoted as $\mathrm{rk}((\chi_v)_v)$. Lastly, given an element $\chi := (\chi_v)_v \in \Omega_\mathfrak{d}$ and a place $\mathfrak{v} \not\in \Sigma(\mathfrak{d})$, we denote by $t_{\chi}(\mathfrak{v})$ the dimension of the following subspace:
    \begin{equation}
        t_\chi(\mathfrak{v}) := \dim_{\mathbb{F}_\ell} \mathrm{Image} \left( \mathrm{loc}_\mathfrak{v}: \Sel(T, \alpha : (\chi_v)_v) \to H^1_\et(\Oh_{K_v}, T) \right).
    \end{equation}
\end{definition}
The variation of the dimensions of local Selmer groups can be detected from the quantity $t_\chi(\mathfrak{v})$. Once again, Lemma \ref{lemma:H_ram} becomes relevant in distinguishing this relation when $T$ is alternating or symmetric. We first recall the proposed relation when $T$ is alternating.
\begin{proposition}[{\cite[Proposition 7.2]{KMR14}} and {\cite[Proposition 5.3]{park2022prime}}]
    Let $T$ be an alternating admissible $\Gal(\overline{K}/K)$-module. Fix a set of places $\Sigma$, a square-free product of places $\mathfrak{d}$ coprime to any places in $\Sigma$, an element $\chi \in \Omega_\mathfrak{d}$, and a place $\mathfrak{v} \not\in \Sigma(\mathfrak{d})$. Then for any $\chi' \in \Omega_{\chi, \mathfrak{v}}$, we have
    \begin{equation*}
        \mathrm{rk}(\chi') - \mathrm{rk}(\chi) = \begin{cases}
            2 &\text{ if } \mathfrak{v} \in \mathcal{P}_2 \text{ and } t_\chi(\mathfrak{v}) = 0 \text{ for exactly } \ell - 1 \text{ out of } \ell(\ell - 1) \text{ many } \chi' \in \Omega_{\chi, \mathfrak{v}}, \\
            1 &\text{ if } \mathfrak{v} \in \mathcal{P}_1 \text{ and } t_\chi(\mathfrak{v}) = 0, \\
            -1 &\text{ if } \mathfrak{v} \in \mathcal{P}_1 \text{ and } t_\chi(\mathfrak{v}) = 1, \\
            -2 &\text{ if } \mathfrak{v} \in \mathcal{P}_2 \text{ and } t_\chi(\mathfrak{v}) = 2, \\
            2 &\text{ otherwise }.
        \end{cases}
    \end{equation*}
\end{proposition}

Next, we consider the case when $T$ is symmetric. Unlike the alternating case, we cannot fully determine the exact relation between variations in dimensions of local Selmer groups and $t_\chi(\mathfrak{v})$. Nevertheless, it is still possible to give bounds on such variations using $t_\chi(\mathfrak{v})$.
\begin{proposition} \label{prop:rank_diff}
    Let $T$ be a symmetric admissible $\Gal(\overline{K}/K)$-module. Fix a set of places $\Sigma$, a square-free product of places $\mathfrak{d}$ coprime to any places in $\Sigma$, an element $\chi \in \Omega_\mathfrak{d}$, and a place $\mathfrak{v} \not\in \Sigma(\mathfrak{d})$. 
    \begin{itemize}
    \item Suppose $\mathfrak{v} \in \mathcal{P}_0 \cup \mathcal{P}_1$. Then for any $\chi' \in \Omega_{\chi, \mathfrak{v}}$, we have
    \begin{equation*}
        \mathrm{rk}(\chi') - \mathrm{rk}(\chi) = \begin{cases}
            1 &\text{ if } \mathfrak{v} \in \mathcal{P}_1 \text{ and } t_\chi(\mathfrak{v}) = 0, \text{ for exactly } \ell - 1 \text{ out of } \ell(\ell-1) \text{ many } \chi' \in \Omega_{\chi,\mathfrak{v}}, \\
            -1 &\text{ if } \mathfrak{v} \in \mathcal{P}_1 \text{ and } t_\chi(\mathfrak{v}) = 1, \\
            0 &\text{ otherwise }.
        \end{cases}
    \end{equation*}
    \item Suppose $\mathfrak{v} \in \mathcal{P}_2$. Then for any $\chi' \in \Omega_{\chi,\mathfrak{v}}$, we have
    \begin{equation*}
        -t_{\chi}(\mathfrak{v}) \leq \mathrm{rk}(\chi') - \mathrm{rk}(\chi) \leq 2 - 2t_{\chi}(\mathfrak{v}).
    \end{equation*}
    \end{itemize}
\end{proposition}
\begin{proof}
    The proof follows from the ideas presented in {\cite[Proposition 7.2]{KMR14}} and {\cite[Proposition 5.3]{park2022prime}}. We define two types of Selmer structures
    \begin{align}
        \begin{split}
            \mathrm{Sel}(T, \alpha : (\chi_v)_v)^{(\mathfrak{v})} &:= \mathrm{Ker} \left( \oplus_v \mathrm{loc}_v: H^1(K,T) \to \bigoplus_{v \neq \mathfrak{v}} \frac{H^1_{\et}(K_v,T)}{\alpha_v(\chi_v)} \right), \\
            \mathrm{Sel}(T, \alpha : (\chi_v)_v)_{(\mathfrak{v})} &:= \mathrm{Ker} \left( \mathrm{loc}_{\mathfrak{v}}: \mathrm{Sel}(T, \alpha : (\chi_v)_v)^{(\mathfrak{v})} \to H^1_{\et}(K_\mathfrak{v}, T) \right).
        \end{split}
    \end{align}
    By definition, we have for any $\chi' \in \Omega_{\chi,\mathfrak{v}}$,
    \begin{equation*}
        \mathrm{Sel}(T, \alpha : (\chi'_v)_v)_{(\mathfrak{v})} \subset \mathrm{Sel}(T, \alpha : (\chi_v)_v) \subset \mathrm{Sel}(T, \alpha : (\chi'_v)_v)^{(\mathfrak{v})},
    \end{equation*}
    and we also have
    \begin{equation*}
        \mathrm{Sel}(T, \alpha : (\chi'_v)_v)_{(\mathfrak{v})} \subset \mathrm{Sel}(T, \alpha : (\chi'_v)_v) \subset \mathrm{Sel}(T, \alpha : (\chi'_v)_v)^{(\mathfrak{v})}.
    \end{equation*}
    By Poitou-Tate duality \cite[Theorem I.4.10]{Milne86}, the image 
    \begin{equation*}
        V_{\mathfrak{v}} := \mathrm{loc}_\mathfrak{v}\left( \mathrm{Sel}(T, \alpha : (\chi'_v)_v)^{(\mathfrak{v})} \right) \subset H^1(K_\mathfrak{v}, T)
    \end{equation*}
    is a maximal isotropic subspace with respect to the pairing $\mu_{T,\mathfrak{v}}$, so $V_{\mathfrak{v}} \in \mathcal{H}(\mu_{T,\mathfrak{v}})$. We then have two short exact sequences
    \begin{align*}
        0 \to \mathrm{Sel}(T, \alpha : (\chi'_v)_v)_{(\mathfrak{v})} \to &\mathrm{Sel}(T, \alpha: (\chi_v)_v) \to V_\mathfrak{v} \cap H^1_\et(\mathcal{O}_{K_\mathfrak{v}}, T) \to 0, \\
        0 \to \mathrm{Sel}(T, \alpha : (\chi'_v)_v)_{(\mathfrak{v})} \to &\mathrm{Sel}(T, \alpha: (\chi'_v)_v) \to V_\mathfrak{v} \cap \alpha_\mathfrak{v}(\chi'_{\mathfrak{v}}) \to 0.
    \end{align*}
    Notice that $t_\chi(\mathfrak{v}) = \dim_{\mathbb{F}_\ell}(V_\mathfrak{v} \cap H^1_\et(\mathcal{O}_{K_\mathfrak{v}}, T))$. Hence, by comparing the dimensions, we obtain
    \begin{equation}
        \dim_{\mathbb{F}_\ell} \mathrm{Sel}(T, \alpha: (\chi'_v)_v) - \dim_{\mathbb{F}_\ell} \mathrm{Sel}(T, \alpha: (\chi_v)_v) = \dim_{\mathbb{F}_\ell} (V_\mathfrak{v} \cap \alpha_\mathfrak{v}(\chi'_\mathfrak{v})) - t_\chi(\mathfrak{v}).
    \end{equation}
    Now we compute this difference based on which $0 \leq i \leq 2$ our place $\mathfrak{v}$ is an element of $\mathcal{P}_i$.
    \begin{itemize}
        \item If $i = 0$, then $H^1(K_\mathfrak{v},T) = 0$, hence $\dim_{\mathbb{F}_\ell} (V_\mathfrak{v} \cap \alpha_\mathfrak{v}(\chi'_\mathfrak{v})) = t_\chi(\mathfrak{v}) = 0$. The difference between the dimensions of two Selmer structures is hence always $0$.
        \item If $i = 1$, then we recall that the twisting data $\alpha_{\mathfrak{v}}$ induces a bijection
        \begin{equation*}
            \alpha_v: \frac{\mathrm{Hom}(\Gal(\overline{K}_v/K_v), \mu_\ell)}{\mathrm{Aut}(\mu_\ell)} \to \{H^1_{\et}(\Oh_{K_v}, T)\} \cup \mathcal{H}_{ram}(\mu_{T,v}).
        \end{equation*}
        If $t_\chi(\mathfrak{v}) = 1$, then $V_\mathfrak{v} = H^1_{\et}(\mathcal{O}_{K_v},T)$. Since $\alpha_\mathfrak{v}(\chi'_\mathfrak{v}) \cap H^1_{\et}(\mathcal{O}_{K_v},T) = 0$, it follows that the difference between the dimensions of two Selmer structures is equal to $-1$. 
        
        On the other hand, if $t_\chi(\mathfrak{v}) = 0$, then there exists $1$ out of $\ell$ many ramified maximal isotropic subspaces in $\mathcal{H}_{ram}(\mu_{T,\mathfrak{v}})$ that is equal to $V_{\mathfrak{v}}$. Hence, $\ell - 1$ many characters out of $\ell(\ell - 1)$ many characters in $\Omega_{\chi,\mathfrak{v}}$ satisfy $\alpha_\mathfrak{v}(\chi'_\mathfrak{v}) = V_\mathfrak{v}$. In such a case, the difference between the dimensions of two Selmer structures is equal to $1$. 
        
        For the rest of the characters, because all maximal isotropic subspaces in $\mathcal{H}_{ram}(\mu_{T,\mathfrak{v}})$ are of dimension $1$, we obtain $V_{\mathfrak{v}} \cap \alpha_\mathfrak{v}(\chi'_\mathfrak{v}) = 0$. Therefore, the difference between the dimensions of two Selmer structures is equal to $0$.
        \item If $i = 2$, we cannot use the reasoning provided above, because the twisting data $\alpha_\mathfrak{v}$ is no longer a bijection. Nevertheless, we still have the upper and lower bounds
        \begin{equation*}
            0 \leq \dim_{\mathbb{F}_\ell} V_{\mathfrak{v}} \cap \alpha_\mathfrak{v}(\chi'_\mathfrak{v}) \leq 2 - t_\chi(\mathfrak{v}).
        \end{equation*}
        This gives the desired upper and lower bounds on the differences of the dimensions between two Selmer structures.
    \end{itemize}
\end{proof}

Both previous propositions imply that the distribution of $t_\chi(\mathfrak{v})$ can be utilized to determine the variation of dimensions of local Selmer structures of $T$. The following proposition, which applies in common to alternating and symmetric admissible modules, can be obtained using Chebotarev density theorem.

\begin{proposition}[{\cite[Proposition 5.4]{park2022prime}}] \label{prop:local-Selmer}
    Let $T$ be an admissible $\Gal(\overline{K}/K)$-module. Fix a set of places $\Sigma$, a square-free product of places $\mathfrak{d}$ coprime to any places in $\Sigma$, and an element $\chi \in \Omega_\mathfrak{d}$.

    Let $d_{i,j}$ be given by the following table:
\begin{center}
    \begin{tabular}{|c||c|c|c|}
    \hline
    $d_{i,j}$ & $i = 0$ & $i = 1$ & $i = 2$ \\
    \hline
    \hline
    $j = -2$ & $\times$ & $\times$ & $1 - (\ell+1) \ell^{-\text{rk}(\chi)} + \ell^{1-2\text{rk}(\chi)}$ \\
    \hline
    $j = -1$ & $\times$ & $1-\ell^{-\text{rk}(\chi)}$ & $\times$ \\
    \hline
    $j = 0$ & $1$ & $\times$ & $(\ell+1)(\ell^{-\text{rk}(\chi)} -\ell^{-2\text{rk}(\chi)})$ \\
    \hline
    $j = 1$ & $\times$ & $ \ell^{-\text{rk}(\chi)}$ & $\times$ \\
    \hline 
    $j = 2$ & $\times$ & $\times$ & $\ell^{-2\text{rk}(\chi)}$ \\
    \hline
    \end{tabular}
\end{center}
Here, the term "$\times$" denotes the case where such a difference of ranks cannot occur. Let $D_{T,\ell,q} > 0$ be a constant defined as
\begin{equation}
    D_{T,\ell,q} := \ell^{\text{max}_{\chi \in \Omega_T} \left( \text{rk}(\chi) \right)}
\end{equation}
Then for every $$d > \frac{12 \log \ell + 2\log D_{T,\ell,q} + (6 \log \ell) \cdot \# \Sigma_T(\mathfrak{d})}{\log q},$$ we have
\begin{equation}
    \left| \frac{\# \{ \mathfrak{v} \in \mathcal{P}_i \; | \; \deg \mathfrak{v} = d, \mathfrak{v} \not\in \Sigma_T(\sigma) \text{ and } t_\chi(\mathfrak{v}) = j \}}{\# \{\mathfrak{v} \in \mathcal{P}_i \; | \; \deg \mathfrak{v} = d, \mathfrak{v} \not\in \Sigma_T(\sigma)\}} - d_{i,j} \right| < (16 \cdot D_{T,\ell,q} \cdot \ell^4) \cdot \ell^{3 \# \Sigma_T(\sigma)} \cdot q^{-\frac{d}{2}}.
\end{equation}
\end{proposition}
\begin{proof}
    The statement follows from adapting the proof presented in \cite[Proposition 5.4]{park2022prime}, because the arguments presented in the proof does not depend on the nature of the pairing $e_T: T \times T \to \mu_\ell$. We present a short summary of the proof. We denote by $\mathrm{Res}_T$ the following injection obtained from the inflation-restriction sequence:
    \begin{equation*}
        \mathrm{Res}_T: H^1_\et(K, T) \to \mathrm{Hom}(\Gal(\overline{K}/K(T)), T)^{\Gal(K(T)/K)}.
    \end{equation*}
    Given an element $\chi := (\chi_v)_v \in \Omega_\mathfrak{d}$, let $F_{\mathfrak{d},\chi}$ be the fixed field of the following subgroup of $\Gal(\overline{K}/K(T))$:
    \begin{equation*}
        \bigcap_{c \in \Sel(T, \alpha : (\chi_v)_v)} \mathrm{Ker} \left( \mathrm{Res}_T(c): \Gal(\overline{K}/K(T)) \to T \right).
    \end{equation*}
    As shown in \cite[Proposition 9.3]{KMR14}, the field $F_{\mathfrak{d}, \chi}$ is a Galois extension over $K$ satisfying the following properties.
    \begin{itemize}
        \item $\Gal(F_{\mathfrak{d},\chi}/K(T)) \cong T^{\mathrm{rk}(\chi)}$,
        \item $F_{\mathfrak{d},\chi}/K$ is unramified outside of places of $\Sigma_T(\mathfrak{d})$.
    \end{itemize}
    The fact that the constant field of $F_{\mathfrak{d}, \chi}$ is identical to that of $K$ follows from statement (7) of Condition \ref{condition:admissible}, from which one can deduce that any basis element $c \in \Sel(T, \alpha : (\chi_v)_v)$ maps the arithmetic Frobenius $\Frob_q \in \Gal(\overline{\mathbb{F}}_q/\mathbb{F}_q)$ to the identity element of $T$. This allows us to apply the effective Chebotarev density theorem over conjugacy-invariant subsets of $\Gal(F_{\mathfrak{d},\chi}/K)$, which relates the densities $d_{i,j}$ with according two conditions on $\mathfrak{v} \in \mathcal{P}_i$ and $t_\chi(\mathfrak{v})=j$, see \cite[Proposition 9.4]{KMR14} for more details.
\end{proof}

\begin{remark} \label{remark:motivation-Markov}
Using the Chebotarev density theorem provided above, we can motivate the construction of the Markov operator $M_\ell$ from the previous section.
\begin{itemize}
    \item \textbf{Alternating Case}: For places $\mathfrak{v} \in \mathcal{P}_1$, recall that for any $\chi' \in \Omega_{\chi, \mathfrak{v}}$, 
    \begin{equation*}
        \mathrm{rk}(\chi') - \mathrm{rk}(\chi) = \begin{cases}
            1 &\text{ if } t_\chi(\mathfrak{v}) = 0, \\
            -1 &\text{ if } t_\chi(\mathfrak{v}) = 1.
        \end{cases}
    \end{equation*}
Assuming that the global characters $\chi'$ ramified at $\Sigma(\mathfrak{d} \mathfrak{v})$ that restricts to $\chi \in \Omega_{\mathfrak{d}}$ equidistributes over $\Omega_{\chi, \mathfrak{v}}$ (which will be shown in the upcoming section), we can combine the probabilities appearing in the previous proposition to obtain an educated guess
    \begin{equation*}
        \mathbb{P}[\dim_{\mathbb{F}_\ell} \mathrm{Sel}(T, \alpha : \chi') - \dim_{\mathbb{F}_\ell} \mathrm{Sel}(T, \alpha : \chi) = j] =_? \begin{cases}
            1 - \ell^{-\mathrm{rk}(\chi)} &\text{ if } j = -1, \\
            \ell^{-\mathrm{rk}(\chi)} &\text{ if } j = 1, \\
            0 &\text{ otherwise }.
        \end{cases}
    \end{equation*}
When we substitute $r = \mathrm{rk}(\chi)$ and $s = \mathrm{rk}(\chi) - j$, then we obtain the Markov operator $M_\ell$ for the alternating case.

For places $\mathfrak{v} \in \mathcal{P}_2$, recall that for any $\chi' \in \Omega_{\chi, \mathfrak{v}}$, 
    \begin{equation*}
        \mathrm{rk}(\chi') - \mathrm{rk}(\chi) = \begin{cases}
            2 &\text{ if } t_\chi(\mathfrak{v}) = 0, \text{ for exactly } \ell - 1 \text{ out of } \ell(\ell-1) \text{ many } \chi' \in \Omega_{\chi,\mathfrak{v}}, \\
            -2 &\text{ if } t_\chi(\mathfrak{v}) = 2, \\
            0 &\text{ otherwise }.
        \end{cases}
    \end{equation*}
Assuming equidistribution results, we can combine the probabilities appearing in the previous proposition to obtain an educated guess
    \begin{equation*}
        \mathbb{P}[\dim_{\mathbb{F}_\ell} \mathrm{Sel}(T, \alpha : \chi') - \dim_{\mathbb{F}_\ell} \mathrm{Sel}(T, \alpha : \chi) = j] =_? \begin{cases}
            (1 - \ell^{-\mathrm{rk}(\chi)})(1-\ell^{-\mathrm{rk}(\chi)-1}) &\text{ if } j = -2, \\
            (\ell + 1)\ell^{-\mathrm{rk}(\chi)} - (\ell + 1/\ell) \ell^{-2 \mathrm{rk}(\chi)} &\text{ if } j = 0, \\
            \ell^{-2 \cdot \mathrm{rk}(\chi) - 1} &\text{ if } j = 2, \\
            0 &\text{ otherwise }.
        \end{cases}
    \end{equation*}
When we substitute $r = \mathrm{rk}(\chi)$ and $s = \mathrm{rk}(\chi) - j$, then we obtain the Markov operator $M_\ell^2$ for the alternating case.

    \item \textbf{Symmetric Case}: For places $\mathfrak{v} \in \mathcal{P}_1$, recall that for any $\chi' \in \Omega_{\chi, \mathfrak{v}}$, 
    \begin{equation*}
        \mathrm{rk}(\chi') - \mathrm{rk}(\chi) = \begin{cases}
            1 &\text{ if } t_\chi(\mathfrak{v}) = 0, \text{ for exactly } \ell - 1 \text{ out of } \ell(\ell-1) \text{ many } \chi' \in \Omega_{\chi,\mathfrak{v}}, \\
            0 &\text{ if } t_\chi(\mathfrak{v}) = 0, \text{ for other } (\ell - 1)^2 \text{ out of } \ell(\ell-1) \text{ many } \chi' \in \Omega_{\chi,\mathfrak{v}},\\
            -1 &\text{ if } t_\chi(\mathfrak{v}) = 1.
        \end{cases}
    \end{equation*}
Assuming that the global characters $\chi'$ ramified at $\Sigma(\mathfrak{d} \mathfrak{v})$ that restricts to $\chi \in \Omega_{\mathfrak{d}}$ equidistributes over $\Omega_{\chi, \mathfrak{v}}$ (which will be shown in the upcoming section), we can combine the probabilities appearing in the previous proposition to obtain an educated guess
    \begin{equation*}
        \mathbb{P}[\dim_{\mathbb{F}_\ell} \mathrm{Sel}(T, \alpha : \chi') - \dim_{\mathbb{F}_\ell} \mathrm{Sel}(T, \alpha : \chi) = j] =_? \begin{cases}
            1 - \ell^{-\mathrm{rk}(\chi)} &\text{ if } j = -1, \\
            \ell^{-\mathrm{rk}(\chi)} \cdot (1 - \ell^{-1}) &\text{ if } j = 0, \\
            \ell^{-\mathrm{rk}(\chi) - 1} &\text{ if } j = 1, \\
            0 &\text{ otherwise }.
        \end{cases}
    \end{equation*}
When we substitute $r = \mathrm{rk}(\chi)$ and $s = \mathrm{rk}(\chi) - j$, then we obtain the Markov operator $M_\ell$ for the symmetric case.

For places $\mathfrak{v} \in \mathcal{P}_2$, the fact that the twisting data $\alpha_\mathfrak{v}$ is not a bijection prevents us from directly relating the predicted variations in dimensions of Selmer groups with respect to increasing a place of ramification. Nevertheless, we may obtain the following heuristic: there exist parameters $s_1, s_2, s_3 \in [0,1]$ such that the probability $\mathbb{P}[\dim_{\mathbb{F}_\ell} \mathrm{Sel}(T, \alpha : \chi') - \dim_{\mathbb{F}_\ell} \mathrm{Sel}(T, \alpha : \chi) = j]$ is given by
\begin{equation*}
    =_? \begin{cases}
        1 - (\ell + 1)\ell^{-\mathrm{rk}(\chi)} + \ell^{-2\mathrm{rk}(\chi) + 1} &\text{ if } j = -2, \\
        (1-s_1) \cdot (\ell + 1) \cdot (\ell^{-\mathrm{rk}(\chi)} - \ell^{-2 \mathrm{rk}(\chi)}) &\text{ if } j = -1, \\
        s_1 \cdot (\ell + 1) \cdot (\ell^{-\mathrm{rk}(\chi)} - \ell^{-2 \mathrm{rk}(\chi)}) + s_2 \cdot \ell^{-2 \mathrm{rk}(\chi)} &\text{ if } j = 0, \\
        s_3 \cdot \ell^{-2 \mathrm{rk}(\chi)} &\text{ if } j = 1, \\
        (1-s_2-s_3) \cdot \ell^{-2 \mathrm{rk}(\chi)} &\text{ if } j = 2, \\
        0 &\text{ otherwise}.
    \end{cases}
\end{equation*}
We note that when $s_1 = \ell^{-1}$, $s_2 = 1 - \ell^{-1}$, and $s_3 = \frac{\ell^2-1}{\ell^3}$, then the above probabilities are transition probabilities of the Markov operator $M_\ell^2$ for the symmetric case. Even if we do not know the exact values of $s_1, s_2, s_3$, we can still hope to obtain that for sufficiently large $\mathrm{rk}(\chi)$, the first moment of differences of the dimensions between two Selmer structures is less than $\mathrm{rk}(\chi)$.
\end{itemize}
\end{remark}

We will rigorously prove these two heuristic observations in the next section using equidistribution of global characters to Cartesian product of local characters.

\subsection{Proof of Theorem \ref{thm:Galois_module-alternating}}

Proposition \ref{prop:local-Selmer} computes the variations of dimensions of local Selmer groups assuming that one can choose finite Cartesian products of local characters uniformly at random from $\Omega_\mathfrak{d}$. To utilize this proposition in computing the variation of dimensions of global Selmer groups, we require that the restriction of global characters $\chi \in \mathrm{Hom}(\Gal(\overline{K}/K), \mu_\ell)$ to finite Cartesian products of local characters $\Omega_\mathfrak{d}$ for some $\mathfrak{d}$ converges in distribution to a uniform distribution. This claim holds true for global characters induced from the set $\hat{F}_{(n,N),(w,w')}(\mathbb{F}_q)$ as shown in Lemma \ref{lemma:equidistribution}, because every polynomial in such a set has an irreducible factor in $\mathcal{P}_0$ of degree greater than $\mathfrak{n}$. Using this lemma, we obtain an adaptation of the analytic argument presented in \cite[Proposition 5.13, Proposition 6.11]{park2022prime}.
\begin{proposition} \label{prop:major_contribution-alternating}
    Suppose $w < 2 m_{n,q}$, and $w' = (1-\epsilon) w$ for some small enough $0 < \epsilon < 1$. Suppose that $n$ satisfies $m_{n,q} > \max \left(e^{e^e}, \; \sum_{v \in \Sigma_T} \deg v, \; 6 \log \ell + 2 \right).$

    Let $T$ be an alternating admissible $\mathrm{Gal}(\overline{K}/K)$-module with a twisting data $\alpha := (\alpha_v)_v$. Then there exist fixed constants $\tilde{B}_{T, \ell, q}$ depending only on $T, \ell, q$ and $\gamma_T$ depending only on $T$ such that
        \begin{align}
        \begin{split}
            & \; \; \; \; \left| \frac{\#\{ f \in \hat{F}_{(n,N),(w,w')}(\mathbb{F}_q) \; | \; \dim_{\mathbb{F}_\ell} \Sel(T, \alpha : \overline{\chi}_f) = r \}}{\# \hat{F}_{(n,N),(w,w')}(\mathbb{F}_q)} - \rho(r) \cdot \prod_{k=1}^\infty \frac{1}{1+\ell^{-k}} \prod_{k=1}^r \frac{\ell}{\ell^k-1} \right| \\
            &< \tilde{B}_{T, \ell, q} \cdot (n \log q)^{4 \epsilon \log \ell} \cdot \left( (n \log q)^{-m_{n,q}} + \gamma_T^{w' - 1} \right),
        \end{split}
        \end{align}
        where
        \begin{equation}
        \rho(r) := \lim_{n \to \infty} \left(\frac{1}{2} + (-1)^r \cdot \left(\frac{1}{2} - \sum_{\substack{k \geq 0 \\ k \equiv r \text{ mod } 2}} \frac{\#\{f \in F_n(\mathbb{F}_q) \cap \mathcal{P}_0 \; | \; \dim_{\mathbb{F}_\ell} \Sel(T, \alpha : \overline{\chi}_f) = r\}}{\#F_n(\mathbb{F}_q) \cap \mathcal{P}_0} \right) \right).
    \end{equation}
    We note that $\rho(r) = \frac{1}{2}$ for all non-negative integer $r \geq 0$ if and only if $\delta(\mathcal{P}_1) \neq 0$.
\end{proposition}
\begin{proof}
    We summarize the key ingredients needed to prove the desired statement. Further analytic details of the proof follows from a line-by-line adaptation of the proof of \cite[Proposition 5.13, Proposition 6.11]{park2022prime}, where we replace the $p$-torsion subgroup $E[p]$ with our chocie of an alternating admissible module $T$.
    
    The global Selmer group $\Sel(T, \alpha: \overline{\chi}_f)$ can be rewritten as a local Selmer group as shown in the following equation:
    \begin{equation}
        \Sel(T, \alpha : \overline{\chi}_f) := \Sel(T, \alpha : (({\chi}_{f,v})_{v \in \Sigma_T(\overline{f})}).
    \end{equation}
    
    Let $h$ be a squarefree polynomial of degree at most $n - \mathfrak{n}$ that satisfies the following conditions.
    \begin{itemize}
        \item $h$ has $w-1$ many distinct irreducible factors.
        \item $h$ has $w'-1$ many distinct irreducible factors of degree greater than $\mathfrak{n}$.
        \item All irreducible factors of $h$ are elements in the set $\Sigma_T \cup \mathcal{P}_1 \cup \mathcal{P}_2$.
    \end{itemize}
    Because every element in $\hat{F}_{(n,N), (w,w')}(\mathbb{F}_q)$ has an irreducible factor in $\mathcal{P}_0$ of degree greater than $\mathfrak{n}$, Lemma \ref{lemma:equidistribution} implies that the pushforward of the uniform distribution via the canonical projection
    \begin{equation}
    \begin{split}
        \left\{ f \in \hat{F}_{(n,N),(w,w')}(\mathbb{F}_q) \; : \; h \text{ divides } f, \; v \in \mathcal{P}_0 \text{ for every } v \mid \frac{f}{h} \right\} &\to \Omega_{\Sigma_T(h)}, \\
        f &\mapsto ({\chi}_{f,v})_{v \in \Sigma_T(h)},
    \end{split}
    \end{equation}
    converges to a uniform distribution as $n$ grows arbitrarily large with explicit rate of convergence. By Proposition \ref{prop:local-Selmer}, there exists a probability distribution $\mu^*: \mathbb{Z}_{\geq 0} \to [0,1]$ and an explicit constant ${B}_{T, \ell, q} > 0$ such that
    \begin{equation}
        \left| \frac{\#\{ f \in \hat{F}_{(n,N),(w,w')}(\mathbb{F}_q) \; | \; \dim_{\mathbb{F}_\ell} \Sel(T, \alpha : \overline{\chi}_f) = r \}}{\# \hat{F}_{(n,N),(w,w')}(\mathbb{F}_q)} - \left( M_T^{w'-1} \mu^* \right)(r) \right| < {B}_{T, \ell, q} \cdot (n \log q)^{-m_{n,q}}
    \end{equation}
    for all $r \geq 0$. The constant ${B}_{T,\ell,q}$ is given by
    \begin{equation}
        B_{T,\ell,q} := 10 \ell \cdot \hat{C}_{T,\ell} + 80 \cdot \ell^4 \cdot D_{T,\ell,q},
    \end{equation}
    where $\hat{C}_{T,\ell}$ is the explicit constant from Lemma \ref{lemma:equidistribution}, and $D_{T,\ell,q}$ is the explicit constant from Proposition \ref{prop:local-Selmer}.
    
    The statement of the proposition then follows from combining the above inequality with the geometric rate of convergence of the Markov operator $M_T$ computed in Theorem \ref{thm:Markov-alternating}. The explicit constant $\tilde{B}_{T,\ell,q}$ is given by (as deduced in \cite[p. 3316, before Equation (144)]{park2022prime})
    \begin{equation}
        \tilde{B}_{T,\ell,q} := B_{T,\ell,q} + c_T \cdot D_{T,\ell,q},
    \end{equation}
    where $c_T$ is the explicit constant appearing in Theorem \ref{thm:Markov-alternating}.
    
    To deal with the case when $\delta(\mathcal{P}_1) = 0$ (which can only happen when $\ell \neq 2$ if $T$ is an alternating admissible module), we proceed as in the proof of \cite[Theorem 8.2(i)]{KMR13}. The Markov operator $M_T$ satisfies the following relation for any probability distribution $\delta: \mathbb{Z}_{\geq 0} \to [0,1]$:
    \begin{equation}
        \rho(M_T \delta) = \rho(\delta).
    \end{equation}
    We denote by $\delta_T: \mathbb{Z}_{\geq 0} \to [0,1]$ the probability distribution
    \begin{equation}
        \delta_T(r) := \frac{\#\{\chi \in \Omega_T \; | \; \mathrm{rk}(\chi) = r\}}{\# \Omega_T}.
    \end{equation}
    We note that $\delta_T(r) \neq \frac{1}{2}$ for all non-negative $r \geq 0$, because the cardinality of $\Hom(\Gal(\overline{K}_v/K), \mu_\ell)$ is a power of $\ell$, in particular odd as long as $\ell \neq 2$. By Lemma \ref{lemma:equidistribution}, we obtain
    \begin{equation}
        \delta_T(r) = \lim_{n \to \infty} \frac{\#\{f \in F_n(\mathbb{F}_q) \cap \mathcal{P}_0 \; | \; \dim_{\mathbb{F}_\ell} \Sel(T, \alpha : \overline{\chi}_f) = r\}}{\#F_n(\mathbb{F}_q) \cap \mathcal{P}_0}.
    \end{equation}
    By induction on $w'$, we obtain
    \begin{equation}
        \rho(M_T^{w'-1} \mu^*) = \rho(\mu^*) = \rho(\delta_T).
    \end{equation}
    The construction of $\rho(r)$ follows from using the relation that $\rho(r) = \frac{1}{2} + (-1)^r \cdot \left( \frac{1}{2} - \rho(\delta_T) \right)$.
\end{proof}

Combining Proposition \ref{proposition:fan_approximation} and Proposition \ref{prop:major_contribution-alternating}, and optimizing the rate of convergence by choosing $\epsilon := (\log \log m_{n,q})^{-1}$ appearing in Proposition \ref{proposition:fan_approximation}, we prove Theorem \ref{thm:Galois_module-alternating}. Explicit computations on how such a choice of $\epsilon$ gives rise to the explicit rate of convergence can be obtained from a line-by-line adaptation of \cite[Proof of Theorem 1.2, pp.3316--3317]{park2022prime}. The implied constant $B_{alt}$ can be deduced from \cite[p. 3317, before Equation (151)]{park2022prime}: it is given by
\begin{equation}
    B_{alt} := 12 \cdot (B_{T,\ell,q} + c_T) \cdot D_{T,\ell,q}.
\end{equation}
As for the exponent $\alpha(T)$, the first two terms originate from Proposition \ref{proposition:fan_approximation}, whereas the last term originates from Proposition \ref{prop:major_contribution-alternating}.

\begin{remark}
We recall that the proof of Proposition \ref{prop:major_contribution-alternating} states that $\rho$ is a rational number whose denominator term is a power of $\ell$. When $\delta(\mathcal{P}_1) = 0$ and $\ell \geq 3$, the fact that $\rho(r) \neq \frac{1}{2}$ for all non-negative integers $r \geq 0$ implies that there is a disparity in the distribution of dimensions of Selmer structures $\mathrm{Sel}(T, \alpha : \overline{\chi}_f)$. The disparity in dimensions of Selmer groups was first observed in \cite[Theorem A]{KMR13}, where disparity in distribution of 2-Selmer groups of quadratic twists of elliptic curves may be observed if these curves are ordered in a non-canonical manner (for instance, with respect to a fan structure over number fields). Generalizations of the disparity phenomena for quadratic twist families of principally polarized abelian varieties can be also found in the works by Yu \cite{Yu16} and Morgan \cite{Mor19}.
\end{remark}

\subsection{Proof of Theorem \ref{thm:Galois_module-symmetric}}
Unlike the proof of Theorem \ref{thm:Galois_module-alternating}, we need further modifications to the strategy of the proof of \cite[Proposition 5.13, Proposition 6.11]{park2022prime}. Hence, instead of outlining the key ingredients of the proof as in Theorem \ref{thm:Galois_module-alternating}, we give respective (though non-trivial) modifications to previous relevant results.

The following two propositions are modifications to \cite[Proposition 5.13]{park2022prime}. We first consider variations in dimensions of local Selmer structures after changing local conditions at $\mathfrak{v} \in \mathcal{P}_0 \cup \mathcal{P}_1$.
\begin{proposition} \label{prop:P01_contribution}
    Let $T$ be a symmetric admissible $\mathrm{Gal}(\overline{K}/K)$-module. Assume notations and conditions as in Definition \ref{def:phi_d*j*k*}. Choose $k^*$ to be either $0$ or $1$, and suppose $w < 2 m_{n,q}$.

    Let $h$ be a polynomial of degree $D := n - d_{a^*} - d^* j^* \lambda_{d^*,j^*,k^*}$ such that $\phi_{d^*,j^*,k^*}^{-1}(h)$ is non-empty. Given $f \in \phi_{d^*,j^*,k^*}^{-1}(h)$, we let $\omega_f$ and $\omega_f'$ be characters defined as
    \begin{equation}
        \omega_f := (\chi_{f,v})_{v \in \Sigma_f(\overline{h}^*)} \in \Omega_{\overline{h}^*}, \hspace{15pt} \omega_f' := (\chi_{f,v})_{v \in \Sigma_f(\overline{f}^*)} \in \Omega_{\overline{f}^*}.
    \end{equation}
    We denote by $\delta_h: \mathbb{Z}_{\geq 0} \to [0,1]$ the probability distribution
    \begin{equation}
        \delta_h(J) := \frac{\#\{f \in \phi_{d^*,j^*,k^*}^{-1}(h) : \mathrm{rk}(\omega_f) = J\}}{\# \phi^{-1}_{d^*,j^*,k^*}(h)}.
    \end{equation}
    Let $\overline{k} := \lambda_{d^*,j^*,k^*} \cdot k^*$. Then for any $n$ with $m_{n,q} > \max\{\sum_{v \in \Sigma_T} \deg v, 3 \log \ell + 1\}$, there exists a fixed constant $B_{T,\ell,q}$ depending only on $T$, $\ell$, and $q$ such that 
    \begin{equation}
        \left| \frac{\# \{ f \in \phi^{-1}_{d^*,j^*,k^*}(h) : \mathrm{rk}(\omega_f') = J\}}{\# \phi_{d^*,j^*,k^*}^{-1}(h)} - (M_\ell^{\overline{k}} \delta_h)(J) \right| < \lambda_{d^*,j^*,k^*} \cdot B_{T,\ell,q} \cdot (n \log q)^{-2 m_{n,q} + 6 \log \ell + 1},
    \end{equation}
    where $M_\ell$ is the Markov operator constructed in Definition \ref{defn:Markov-symmetric}.
\end{proposition}
\begin{proof}
    The proposition follows from an adaptation of \cite[Proposition 5.13]{park2022prime}. We explain how this adaptation can be achieved for a symmetric admissible module $T$.

    Given two non-negative integers $J_0$ and $J_1$, we have
    \begin{align} \label{eqn:symmetric-1}
    \begin{split}
        & \hspace{15pt} \# \{f \in \phi_{d^*,j^*,k^*}^{-1}(h) : \mathrm{rk}(\omega_f') = J_1, \mathrm{rk}(\omega_f) = J_0\} \\
        &= \sum_{(g_a)_a \in \mathrm{Conf}_{\lambda_{d^*,j^*,k^*}}(\mathcal{P}_{k^*}(d^*))} \# \left\{ f \in \phi_{d^*,j^*,k^*}^{-1}(h) : \frac{f}{f_a h} = \prod_{a=1}^{\lambda_{d^*,j^*,k^*}} g_a^{j^*}, \; \mathrm{rk}(\omega_f') = J_1, \mathrm{rk}(\omega_f) = J_0 \right\}.
    \end{split}  
    \end{align}

    We use the abbreviation $g := \prod_{a=1}^{\lambda_{d^*,j^*,k^*}} g_a$. Denote by $\Omega_{\omega_f,g}$ the subset of local characters $\chi' \in \Omega_{\overline{f}^*}$ such that $\chi'_{g_a}$ is ramified for all $1 \leq a \leq \lambda_{d^*,j^*,k^*}$, and for any $v \in \Sigma_f(\overline{h}^*)$ we have $\omega'_v = (\omega_f)_v$. Furthermore, we use the notation $\widetilde{(g_a)}_a$ to denote elements in $\mathrm{PConf}_{\lambda_{d^*,j^*,k^*}}(\mathcal{P}_{k^*}(d^*))$. Then by using Lemma \ref{lemma:equidistribution} and following through \cite[Equations (91) -- (100)]{park2022prime}, we obtain
    \begin{align} \label{eqn:symmetric-0}
        \begin{split}
            & \hspace{15pt} \left| (\ref{eqn:symmetric-1}) - \frac{\sum_{\widetilde{(g_a)}_a} \# \{\omega' \in \Omega_{\omega_f,g} : \mathrm{rk}(\omega') - \mathrm{rk}(\omega_f) = J_1 - J_0\}}{\sum_{\widetilde{(g_a)}_a} \# \Omega_{\omega_f,g}} \right| \\
            &< \# \phi_{d^*,j^*,k^*}^{-1}(h) \cdot 2 \ell \cdot \hat{C}_{T,\ell} \cdot (n \log q)^{-2 m_{n,q} + 4 \log \ell},
        \end{split}
    \end{align}
    where $\hat{C}_{T,\ell}$ is the explicit constant obtained from Lemma \ref{lemma:equidistribution}, and the summations on both the numerator and denominator range over $\widetilde{(g_a)}_a \in \mathrm{PConf}_{\lambda_{d^*,j^*,k^*}}(\mathcal{P}_{k^*}(d^*))$.

    The next step is to use induction on $\lambda_{d^*,j^*,k^*}$ and iterate Proposition \ref{prop:rank_diff} and Proposition \ref{prop:local-Selmer} to obtain
    \begin{align} \label{eqn:symmetric-induction}
        \begin{split}
            & \hspace{15pt} \left| \frac{\sum_{\widetilde{(g_a)}_a} \# \{\omega' \in \Omega_{\omega_f,g} : \mathrm{rk}(\omega') - \mathrm{rk}(\omega_f) = J_1 - J_0\}}{\sum_{\widetilde{(g_a)}_a} \# \Omega_{\omega_f,g}} - (M_\ell^{\overline{k}} \delta_h)(J_1) \right| \\
            &< 48 \cdot \lambda_{d^*,j^*,k^*} \cdot \ell^4 \cdot D_{T,\ell,q} \cdot (n \log q)^{-2 m_{n,q} + 6 \log \ell + 1}.
        \end{split}
    \end{align}
    The base step (where $\lambda_{d^*,j^*,k^*} = 1$) concerns the case where $\mathrm{PConf}_{\lambda_{d^*,j^*,k^*}}(\mathcal{P}_{k^*}(d^*)) = \mathcal{P}_{k^*}(d)$ and $g = g_1$. We fix $\omega \in \Omega_{\overline{h}^*}$ such that $\mathrm{rk}(\omega) = J_0$. The strategy to prove the base case is identical to \cite[``Base Step'', pp. 3306--3307]{park2022prime}, but the only difference we take into account is the definition of the constants $c_{k^*,J_1 - J_0}$: these are obtained from a combination of Proposition \ref{prop:rank_diff} and Proposition \ref{prop:local-Selmer} (and their explicit values follow from a series of observations we made in Remark \ref{remark:motivation-Markov}):
    \begin{equation}
        c_{k^*,J_1-J_0} = \begin{cases}
            1 - \ell^{-J_0} &\text{ if } k^* = 1, J_1 - J_0 = -1, \\
            \ell^{-J_0} \cdot (1 - \ell^{-1}) &\text{ if } k^* = 1, J_1 - J_0 = 0, \\
            \ell^{-J_0 - 1} &\text{ if } k^* = 1, J_1 - J_0 = 1, \\
            1 &\text{ if } k^* = 0, \\
            0 &\text{ otherwise}.
        \end{cases}
    \end{equation}
    Then we obtain
    \begin{align}
    \begin{split}
        & \hspace{15pt} \left| \frac{\sum_{g_1 \in \mathcal{P}_{k^*}(d^*)} \# \{\omega' \in \Omega_{\omega, g_1} : \mathrm{rk}(\omega') - \mathrm{rk}(\omega) = J_1 - J_0\}}{\sum_{g_1 \in \mathcal{P}_{k^*}(d^*)} \# \Omega_{\omega,g_1}} - c_{k^*,J_1-J_0} \right| \\
        &< 16 \cdot \ell^4 \cdot D_{T,\ell,q} \cdot (n \log q)^{-2 m_{n,q} + 6 \log \ell + 1}
    \end{split}
    \end{align}

    Summing up $\delta_h(J_0)$ over $3$ possible values of $J_0$ such that $c_{k^*,J_1 - J_0} \neq 0, 1$ gives the base case.

    Now we consider the induction step. Suppose equation (\ref{eqn:symmetric-induction}) holds up to $\lambda_{d^*,j^*,k^*} = \overline{\lambda}$. Given $\widetilde{(g_a)}_a \in \mathrm{PConf}_{\overline{\lambda} + 1}(\mathcal{P}_{k^*}(d^*))$ (such that $1 \leq a \leq \overline{\lambda} + 1$), we define $\overline{g} := \prod_{a = 1}^{\overline{\lambda}} g_a = g / g_{\overline{\lambda} + 1}$. Once again, the strategy is identical to \cite[``Induction Step'', pp. 3307--3308]{park2022prime}, but the only difference we need to make is the definition of the constants $c_{k^*,J_1 - (J_0 + J_2)}$ for some integer $-\overline{\lambda} \leq J_2 \leq \overline{\lambda}$. Here, the integer $J_2$ is an auxiliary value such that $\mathrm{rk}(\overline{\omega}) = J_0 + J_2$ given a fixed choice of $\overline{\omega} \in \Omega_{\omega,\overline{g}}$. Defined analogously to the constants $c_{k^*, J_1 - J_0}$, the constants $c_{k^*,J_1 - (J_0 + J_2)}$ are defined as follows, also obtained from a combination of Proposition \ref{prop:rank_diff} and \ref{prop:local-Selmer}:
    \begin{equation}
        c_{k^*,J_1-(J_0 + J_2)} = \begin{cases}
            1 - \ell^{-(J_0 + J_2)} &\text{ if } k^* = 1, J_1 - (J_0 + J_2) = -1, \\
            \ell^{-(J_0 + J_2)} \cdot (1 - \ell^{-1}) &\text{ if } k^* = 1, J_1 - (J_0 + J_2) = 0, \\
            \ell^{-(J_0 + J_2) - 1} &\text{ if } k^* = 1, J_1 - (J_0 + J_2) = 1, \\
            1 &\text{ if } k^* = 0, \\
            0 &\text{ otherwise}.
        \end{cases}
    \end{equation}
    Then one obtains
    \begin{align}
    \begin{split}
        & \hspace{15pt} \left| \frac{\sum_{g_{\overline{\lambda} + 1}} \# \{\omega' \in \Omega_{\overline{\omega}, g_{\overline{\lambda} + 1}} : \mathrm{rk}(\omega') - \mathrm{rk}({\omega}) = J_1 - J_0, \mathrm{rk}(\overline{\omega}) = J_0 + J_2\}}{\sum_{g_{\overline{\lambda} + 1}} \# \Omega_{\overline{\omega},g_{\overline{\lambda} + 1}}} - c_{k^*,J_1-(J_0 + J_2)} \right| \\
        &< 16 \cdot \ell^4 \cdot D_{T,\ell,q} \cdot (n \log q)^{-2 m_{n,q} + 6 \log \ell + 1}
    \end{split}
    \end{align}
    where the summation $g_{\overline{\lambda} + 1}$ ranges over the set $\mathcal{P}_{k^*}(d^*) \setminus \{g_1, \cdots, g_{\overline{\lambda}}\}$. Now we consider the expression
    \begin{equation} \label{eqn:symmetric-2}
        \frac{1}{\# \mathrm{PConf}_{\overline{\lambda}}(\mathcal{P}_{k^*}(d^*))} \sum_{\widetilde{(g_a)}_a} \left( \frac{\sum_{J_2 = -\overline{\lambda}}^{\overline{\lambda}} \# \{\overline{\omega} \in \Omega_{\omega, \overline{g}} : \mathrm{rk}(\overline{\omega}) = J_0 + J_2\} \cdot c_{k^*,J_1 - (J_0 + J_2)} }{\# \Omega_{\omega, \overline{g}}} \right),
    \end{equation}
    where the summation $\widetilde{(g_a)}_a$ ranges over $\mathrm{PConf}_{\overline{\lambda}}(\mathcal{P}_{k^*}(d^*))$. Then we obtain
    \begin{align}
    \begin{split}
        & \hspace{15pt} \left| \frac{\sum_{\widetilde{(g_a)}_a} \# \{\omega' \in \Omega_{\omega_f,g} : \mathrm{rk}(\omega') - \mathrm{rk}(\omega_f) = J_1 - J_0\}}{\sum_{\widetilde{(g_a)}_a} \# \Omega_{\omega_f,g}} - (\ref{eqn:symmetric-2}) \right| \\
        &< 48 \cdot \ell^4 \cdot D_{T,\ell,q} \cdot (n \log q)^{-2 m_{n,q} + 6 \log \ell + 1}.
    \end{split}
    \end{align}
    By the induction hypothesis, we have
    \begin{equation}
        \left| (\ref{eqn:symmetric-2}) - (M_\ell^{k^* \cdot (\overline{\lambda} + 1)} \delta_h)(J_1) \right| < 48 \cdot \overline{\lambda} \cdot \ell^4 \cdot D_{T,\ell,q} \cdot (n \log q)^{-2 m_{n,q} + 6 \log p + 1}.
    \end{equation}
    Combining the two above equations gives equation (\ref{eqn:symmetric-induction}) for $\overline{\lambda} + 1$. We then combine equations (\ref{eqn:symmetric-0}) and (\ref{eqn:symmetric-induction}) to obtain the statement of the proposition, with the constant for the error term $B_{T,\ell,q}$ determined as
    \begin{equation}
        B_{T,\ell,q} := 3 \cdot (2 \ell \cdot \hat{C}_{T,\ell} + 16 \cdot \ell^4 \cdot D_{T,\ell,q}),
    \end{equation}
    where $\hat{C}_{T,\ell}$ is the constant obtained from Lemma \ref{lemma:equidistribution}.
\end{proof}

Next, we consider variations in dimensions of local Selmer structures after changing local conditions at $\mathfrak{v} \in \mathcal{P}_2$. This is the the most technical and involved part of the paper, which stems from the technicality we don't obtain a uniform description of explicit Markov operators as what we saw in Proposition \ref{prop:P01_contribution}. However, their exponential moments can be bounded above. The idea of using auxiliary Markov operators traces back to \cite{KP26}[Section 7], where authors used ``worst case'' Markov operators that govern variations in 2-Selmer groups of hyperelliptic curves $C$ with respect to changing local conditions at places $v$ with $\dim_{\mathbb{F}_2} \mathrm{Jac}(C)[2](K_v) > 2$. 
\begin{proposition} \label{prop:P2_contribution}
    Let $T$ be a symmetric admissible $\mathrm{Gal}(\overline{K}/K)$-module. Assume notations and conditions as in Definition \ref{def:phi_d*j*k*}. Fix $k^* = 2$, and suppose $w < 2 m_{n,q}$.

    Let $h$ be a polynomial of degree $D := n - d_{a^*} - d^* j^* \lambda_{d^*,j^*,2}$ such that $\phi_{d^*,j^*,2}^{-1}(h)$ is non-empty. Given $f \in \phi_{d^*,j^*,2}^{-1}(h)$, we let $\omega_f$ and $\omega_f'$ be characters defined as
    \begin{equation}
        \omega_f := (\chi_{f,v})_{v \in \Sigma_f(\overline{h}^*)} \in \Omega_{\overline{h}^*}, \hspace{15pt} \omega_f' := (\chi_{f,v})_{v \in \Sigma_f(\overline{f}^*)} \in \Omega_{\overline{f}^*}.
    \end{equation}
    We denote by $\delta_h: \mathbb{Z}_{\geq 0} \to [0,1]$ the probability distribution
    \begin{equation}
        \delta_h(J) := \frac{\#\{f \in \phi_{d^*,j^*,2}^{-1}(h) : \mathrm{rk}(\omega_f) = J\}}{\# \phi^{-1}_{d^*,j^*,2}(h)}.
    \end{equation}
    Let $\overline{k} := 2 \lambda_{d^*,j^*,2}$. Choose any $n$ with $m_{n,q} > \max\{\sum_{v \in \Sigma_T} \deg v, 3 \log \ell + 1\}$. 
    \begin{itemize}
        \item There exist auxiliary operators $\mathcal{W}_1, \mathcal{W}_2, \cdots, \mathcal{W}_{\lambda_{d^*,j^*,2}}$ depending on $\phi^{-1}_{d^*,j^*,2}(h)$ and a fixed constant $B_{T,\ell,q}$ depending only on $T$, $\ell$, and $q$ such that 
        \begin{align}
        \begin{split}
        & \hspace{15pt} \left| \frac{\# \{ f \in \phi^{-1}_{d^*,j^*,k^*}(h) : \mathrm{rk}(\omega_f') = J\}}{\# \phi_{d^*,j^*,k^*}^{-1}(h)} - (\mathcal{W}_1 \mathcal{W}_2 \cdots \mathcal{W}_{\lambda_{d^*,j^*,2}} \delta_h)(J) \right| \\
        &< \lambda_{d^*,j^*,2} \cdot B_{T,\ell,q} \cdot (n \log q)^{-2 m_{n,q} + 6 \log \ell + 1}.
        \end{split}
        \end{align}
        \item The auxiliary operators $\mathcal{W}_1, \mathcal{W}_2, \cdots, \mathcal{W}_{\lambda_{d^*,j^*,2}}$ satisfy the condition that $\mathbb{E}[\ell^{\mathcal{W}_i\delta_h}] \leq \mathbb{E}[\ell^{\mathcal{V} \delta_h}]$ for all $1 \leq i \leq \lambda_{d^*,j^*,2}$, where $\mathcal{V}$ is the Markov operator over $\mathbb{Z}_{\geq 0}$ defined by
    \begin{equation}
        \mathcal{V}(r,s) := \begin{cases}
            1 - (\ell + 1) \ell^{-r} + \ell^{-2r + 1} &\text{ if } s = r-2, \\
            (\ell + 1) \cdot (\ell^{-r} - \ell^{-2r}) &\text{ if } s = r, \\
            \ell^{-2r} &\text{ if } s = r + 2, \\
            0 &\text{ otherwise}.
        \end{cases}
    \end{equation}
    \end{itemize}
\end{proposition}
\begin{proof}

    \medskip

    \textbf{[A need for different approach]}

    \medskip

    We proceed as in the proof of Proposition \ref{prop:P01_contribution}, but there are many components in the proof that we need to change. Most crucially, we cannot use a line-by-line adaptation of the methods presented in the proof of the proposition to obtain direct adaptations of equation (\ref{eqn:symmetric-0}) and (\ref{eqn:symmetric-induction}). What is crucially used in obtaining them is the statement that for any place $v \not\in \Sigma(\overline{h}^*)$ and for any $\chi, \chi' \in \Omega_{\overline{h}^*}$ such that $\mathrm{rk}(\chi) = \mathrm{rk}(\chi')$, the governing Markov operators (or Chebotarev conditions) which determine the variations $\dim_{\mathbb{F}_\ell} \mathrm{Sel}(T, \alpha : \omega) - \dim_{\mathbb{F}_\ell} \mathrm{Sel}(T, \alpha : \chi)$ and $\dim_{\mathbb{F}_\ell} \mathrm{Sel}(T, \alpha : \omega') - \dim_{\mathbb{F}_\ell} \mathrm{Sel}(T, \alpha : \chi')$ as $\omega$ varies over $\Omega_{\chi,v}$ and $\omega'$ varies over $\Omega_{\chi',v}$ are identical. This is no longer guaranteed to be true for our case, because the twisting data $\alpha_v$ at places $v \in \mathcal{P}_2$ is no longer a bijection.

    Nevertheless, we can obtain an analogous statement to these equations by fixing a choice of a character $\omega \in \Omega_{\overline{h}^*}$ (which is a finite set given any choice of $n$), obtain auxiliary operators for each such characters, and take averages of the said operators to construct the desired operators $\mathcal{W}_1, \mathcal{W}_2, \cdots, \mathcal{W}_{\lambda_{d^*,j^*,2}}$.

    \medskip 

    \textbf{[Decomposition]}

    \medskip

    Throughout the proof of this proposition, we will abbreviate $\lambda_{d^*,j^*,2} = \overline{\lambda}$.
    As before, consider the following expression given two non-negative integers $J_0, J_{\overline{\lambda}}$:
    \begin{equation} \label{eqn:symmetric-1-2}
        \#\{f \in \phi_{d^*,j^*,2}^{-1}(h): \mathrm{rk}(\omega_f') = J_{\overline{\lambda}}, \mathrm{rk}(\omega_f) = J_0 \}.
    \end{equation}
    We have a set-theoretic bijection
    \begin{equation}
        \phi_{d^*,j^*,2}^{-1}(h) \cong \prod_{j=1}^{\ell-1}\mathrm{Conf}_{\lambda_{d_a,j,0}} \times \mathrm{Conf}_{\overline{\lambda}}(\mathcal{P}_{2}(d^*)) \cong \phi_{d_a}^{-1}(h) \times \frac{\mathrm{PConf}_{\overline{\lambda}}(\mathcal{P}_{2}(d^*))}{S_{\overline{\lambda}}}.
    \end{equation}
    Hence, there exists a natural $S_{\overline{\lambda}}$-cover
    \begin{equation}
        \pi_{d^*,j^*,2}: \phi_{d_a}^{-1}(h) \times \mathrm{PConf}_{\overline{\lambda}}(\mathcal{P}_{2}(d^*)) \to \phi_{d^*,j^*,2}^{-1}(h)
    \end{equation}
    corresponding to enumerating (with order) all irreducible factors of $f$ of degree $d^*$ and multiplicity $j^*$ lying in $\mathcal{P}_2$. By construction, all the fibers at a fixed point $f \in \phi_{d^*,j^*,2}^{-1}(h)$ have identical local Selmer structures $\mathrm{Sel}(T, \alpha : \omega_f')$. Hence, we have
    \begin{align} \label{eqn:symmetric-1-2-Pconf}
    \begin{split}
        & \hspace{15pt} \# S_{\overline{\lambda}} \cdot (\ref{eqn:symmetric-1-2}) \\
        &= \#\left\{(f_a, \widetilde{(g_a)}_a) \in \prod_{j=1}^{\ell-1} \mathrm{Conf}_{\lambda_{d_a,j,0}} \times \mathrm{PConf}_{\overline{\lambda}}(\mathcal{P}_{2}(d^*)): \mathrm{rk}(\omega_f') = J_{\overline{\lambda}}, \; \mathrm{rk}(\omega_f) = J_0 \right\},
    \end{split}
    \end{align}
    where we identify $f = h \cdot f_a \cdot \prod_{a=1}^{\overline{\lambda}} g_a^{j^*}$. To abbreviate the products of two configuration spaces, we will use the notation
    \begin{equation*}
        \pi_{d^*,j^*,2}^{-1}(\phi^{-1}_{d^*,j^*,2}(h)) := \phi_{d_a}^{-1}(h) \times \mathrm{PConf}_{\overline{\lambda}}(\mathcal{P}_{2}(d^*)).
    \end{equation*}
    
    Given a fixed choice of $\omega \in \Omega_{\overline{h}^*}$, consider the following expression:
    \begin{equation} \label{eqn:symmetric-1-2-local}
        \#\left\{(f_a, \widetilde{(g_a)}_a) \in \pi_{d^*,j^*,2}^{-1}(\phi^{-1}_{d^*,j^*,2}(h)): \omega_f = \omega, \mathrm{rk}(\omega_f') = J_{\overline{\lambda}}, \; \mathrm{rk}(\omega) = J_0 \right\}.
    \end{equation}
    By definition, $(\ref{eqn:symmetric-1-2-Pconf}) = \sum_{\omega \in \Omega_{\overline{h}^*}} (\ref{eqn:symmetric-1-2-local})$. We will focus on understanding equation (\ref{eqn:symmetric-1-2-Pconf}) using auxiliary Markov operators. Our goal is to prove via induction the following upper bound of the following expression for any fixed $J_{\overline{\lambda}} \in \mathbb{Z}_{\geq 0}$:
    \begin{equation} \label{eqn:symmetric-2-induction}
        \left| \frac{ \sum_{J_0 = -\infty}^{\infty} (\ref{eqn:symmetric-1-2-Pconf})}{\# \pi_{d^*,j^*,2}^{-1}(\phi_{d^*,j^*,2}^{-1}(h))} - (\mathcal{W}_{\overline{\lambda}} \cdots \mathcal{W}_1 \delta_{h})(J_{\overline{\lambda}}) \right| < \overline{\lambda} \cdot B_{T,\ell,q} \cdot (n \log q)^{-2m_{n,q} + 6 \log \ell + 1},
    \end{equation}
    where $B_{T,\ell,q} := (16 \cdot \ell^4 \cdot D_{T,\ell,q} + 18 \cdot \hat{C}_{T,\ell})$ and all operators $\mathcal{W}_i$ for $1 \leq i \leq \overline{\lambda}$ satisfy $\mathbb{E}[\ell^{\mathcal{W}_i \delta_h}] \leq \mathbb{E}[\ell^{\mathcal{V} \delta_h}]$. This statement in particular proves the statement of the proposition.

    \medskip 

    \textbf{[Relating global and local statistics]}

    \medskip 

    To obtain the inequality (\ref{eqn:symmetric-2-induction}), we first aim to obtain an analogue of equation (\ref{eqn:symmetric-0}). 
    
    We fix $\omega \in \Omega_{\overline{h}^*}$. Given $\widetilde{(g_a)}_a \in \mathrm{PConf}_{\overline{\lambda}}(\mathcal{P}_2(d^*))$, let $g := \prod_{a=1}^{\overline{\lambda}} g_a$. Given $1 \leq b \leq \overline{\lambda}$, we denote by $g_{[b]} := \prod_{a=1}^b g_a$. We denote by $\widetilde{(g_{[b]})_a} := (g_1, \cdots, g_{[b]}) \in \mathrm{PConf}_{b}(\mathcal{P}_2(d^*))$. Given $f \in \phi_{d^*,j^*,2}^{-1}(h)$, we let $\omega_f^{[b]}$ be a character defined as
    \begin{equation*}
        \omega_f^{[b]} := (\chi_{f,v})_{v \in \Sigma_f(\overline{h}^* \cdot g_{[b]})} \in \Omega_{\overline{h}^* \cdot g_{[b]}}.
    \end{equation*}
    We use the variable $J_{b}$ to denote $\mathrm{rk}(\omega_f^{[b]})$. Furthermore, given $\omega \in \Omega_{\overline{h}^*}$, we denote by $\Omega_{\omega,g_{[b]}}$ the subset of local characters $\omega^{[b]} \in \Omega_{\overline{h}^* \cdot g_{[b]}}$ such that $\omega^{[b]}_{g_a}$ is ramified for all $1 \leq a \leq b$, and for any $v \in \Sigma_f(\overline{h}^*)$ we have $\omega^{[b]}_v = \omega_v$.
    
    Our goal is to obtain relations between these two expressions (such that the first expression is nonzero):
    \begin{align*}
        & \frac{\#\left\{(f_a, \widetilde{(g_a)}_a) \in \pi_{d^*,j^*,2}^{-1}(\phi^{-1}_{d^*,j^*,2}(h)): \omega_f = \omega, \mathrm{rk}(\omega_f^{[b]}) = J_{b}, \; \mathrm{rk}(\omega) = J_0 \right\}}{\# \pi_{d^*,j^*,2}^{-1}(\phi_{d^*,j^*,2}^{-1}(h))}, \\
        &\frac{\sum_{\widetilde{(g_{[b]})}_a \in \mathrm{PConf}_{b}(\mathcal{P}_2(d^*))} \# \{ \omega^{[b]} \in \Omega_{\omega, g_{[b]}} : \mathrm{rk}(\omega^{[b]}) - \mathrm{rk}(\omega) = J_{b} - J_0\}}{\sum_{\widetilde{(g_{[b]})}_a \in \mathrm{PConf}_{b}(\mathcal{P}_2(d^*))} \# \Omega_{\omega, g_{[b]}}}.
    \end{align*}

    We start with the first expression. For each fixed $\widetilde{(g^*_a)}_a \in \mathrm{PConf}_{\overline{\lambda}}(\mathcal{P}_2(d^*))$, we have 
    \begin{align*}
        & \hspace{15pt} \frac{\#\{(f_a,\widetilde{(g^*_a)}_a) \in \pi_{d^*,j^*,2}^{-1}(\phi^{-1}_{d^*,j^*,2}(h)): \omega_f = \omega, \mathrm{rk}(\omega_f^{[b]}) - \mathrm{rk}(\omega) = J_{b} - J_0\}}{\#\{(f_a,\widetilde{(g^*_a)}_a) \in \pi_{d^*,j^*,2}^{-1}(\phi^{-1}_{d^*,j^*,2}(h)): \omega_f = \omega\}} \\
        &= \frac{\#\{(f_a,\widetilde{(g^*_a)}_a) \in \pi_{d^*,j^*,2}^{-1}(\phi^{-1}_{d^*,j^*,2}(h)): \omega_f = \omega, \mathrm{rk}(\omega_f^{[b]}) - \mathrm{rk}(\omega) = J_{b} - J_0\} / \# \phi^{-1}_{d_a}(h)}{\#\{(f_a,\widetilde{(g^*_a)}_a) \in \pi_{d^*,j^*,2}^{-1}(\phi^{-1}_{d^*,j^*,2}(h)): \omega_f = \omega\} / \# \phi^{-1}_{d_a}(h)}.
    \end{align*}
    Note that we are dividing both the numerator and denominator by $\# \phi^{-1}_{d_a}(h)$ (not by $\# \pi^{-1}_{d^*,j^*,2}(\phi^{-1}_{d^*,j^*,2}(h))$), because we have fixed a choice of $\widetilde{(g^*_a)}_a \in \mathrm{PConf}_{\overline{\lambda}}(\mathcal{P}_2(d^*))$.
    
    By applying Lemma \ref{lemma:equidistribution} to the numerator and the denominator of the quotient above (and the fact that $w = \omega(f) < 2m_{n,q}$), we obtain for such a given $\widetilde{(g^*_a)}_a \in \mathrm{PConf}_{\overline{\lambda}}(\mathcal{P}_2(d^*))$ (where we use the notation $g_{[b]}^* := \prod_{a=1}^b g_a^*$),
    \begin{align*}
    \begin{split}
        & \hspace{15pt} \biggl| \frac{\#\{(f_a,\widetilde{(g^*_a)}_a) \in \pi_{d^*,j^*,2}^{-1}(\phi^{-1}_{d^*,j^*,2}(h)): \omega_f = \omega, \mathrm{rk}(\omega_f^{[b]}) - \mathrm{rk}(\omega) = J_{b} - J_0\}}{\#\{(f_a,\widetilde{(g^*_a)}_a) \in \pi_{d^*,j^*,2}^{-1}(\phi^{-1}_{d^*,j^*,2}(h)): \omega_f = \omega\}} \\
        & \hspace{150pt} - \frac{\# \{\omega^{[b]} \in \Omega_{\omega, g_{[b]}^*}: \mathrm{rk}(\omega') - \mathrm{rk}(\omega) = J_{b} - J_0\}}{\# \Omega_{\omega, g_{[b]}^*}} \biggr| \\
        &< 2 \cdot \hat{C}_{T,\ell} \cdot (n \log q)^{-2 m_{n,q} + 4 \log \ell},
    \end{split}
    \end{align*}
    where the exponent for $(n \log q)$ grows by $2 \log \ell$ because the density of elements $f \in \phi_{d_a}^{-1}(h)$ such that $\omega_f = \omega$ is equal to $\ell^{-\hat{\omega}(h)} \geq \ell^{-2 m_{n,q}}$.

    We multiply both sides of the expression by $\#\{(f_a,\widetilde{(g^*_a)}_a) \in \pi_{d^*,j^*,2}^{-1}(\phi^{-1}_{d^*,j^*,2}(h)): \omega_f = \omega\}$ to obtain
    \begin{align*}
    \begin{split}
        & \hspace{15pt} \biggl| \#\{(f_a,\widetilde{(g^*_a)}_a) \in \pi_{d^*,j^*,2}^{-1}(\phi^{-1}_{d^*,j^*,2}(h)): \omega_f = \omega, \mathrm{rk}(\omega_f^{[b]}) - \mathrm{rk}(\omega) = J_{b} - J_0\} \\
        & \hspace{15pt} - \frac{\# \{\omega^{[b]} \in \Omega_{\omega, g_{[b]}^*}: \mathrm{rk}(\omega^{[b]}) - \mathrm{rk}(\omega) = J_{b} - J_0\}}{\# \Omega_{\omega, g_{[b]}^*}} \cdot {\#\{(f_a,\widetilde{(g^*_a)}_a) \in \pi_{d^*,j^*,2}^{-1}(\phi^{-1}_{d^*,j^*,2}(h)): \omega_f = \omega\}} \biggr| \\
        &< 2 \cdot \hat{C}_{T,\ell} \cdot (n \log q)^{-2 m_{n,q} + 4 \log \ell} \cdot {\#\{(f_a,\widetilde{(g^*_a)}_a) \in \pi_{d^*,j^*,2}^{-1}(\phi^{-1}_{d^*,j^*,2}(h)): \omega_f = \omega\}}.
    \end{split}
    \end{align*}
    We note that
    \begin{align*}
        & \hspace{15pt} \#\{(f_a,\widetilde{(g^*_a)}_a) \in \pi_{d^*,j^*,2}^{-1}(\phi^{-1}_{d^*,j^*,2}(h)): \omega_f = \omega\} \\
        &= \frac{{\#\{(f_a,\widetilde{(g^*_a)}_a) \in \pi_{d^*,j^*,2}^{-1}(\phi^{-1}_{d^*,j^*,2}(h)): \omega_f = \omega\}}}{\# \phi_{d_a}^{-1}(h)} \cdot \# \phi_{d_a}^{-1}(h).
    \end{align*}
    Hence we can use Lemma \ref{lemma:equidistribution} (which is where the term $\hat{w}(h)$ is defined) to the first term of the right hand side above to obtain
    \begin{align*}
    \begin{split}
        & \hspace{15pt} \biggl| \#\{(f_a,\widetilde{(g^*_a)}_a) \in \pi_{d^*,j^*,2}^{-1}(\phi^{-1}_{d^*,j^*,2}(h)): \omega_f = \omega, \mathrm{rk}(\omega_f^{[b]}) - \mathrm{rk}(\omega) = J_{b} - J_0\} \\
        & \hspace{15pt} - \frac{\# \{\omega^{[b]} \in \Omega_{\omega, g_{[b]}^*}: \mathrm{rk}(\omega^{[b]}) - \mathrm{rk}(\omega) = J_{b} - J_0\}}{\# \Omega_{\omega, g_{[b]}^*}} \cdot \frac{\# \phi_{d_a}^{-1}(h)}{(\ell(\ell-1))^{\hat{\omega}(h)}} \biggr| \\
        &< 3 \cdot \hat{C}_{T,\ell} \cdot (n \log q)^{-2 m_{n,q} + 4 \log \ell} \cdot \# \phi_{d_a}^{-1}(h).
    \end{split}
    \end{align*}
    We then divide both sides by $\# \phi_{d_a}^{-1}(h)$ to obtain
    \begin{align*}
    \begin{split}
        & \hspace{15pt} \biggl| \frac{\#\{(f_a,\widetilde{(g^*_a)}_a) \in \pi_{d^*,j^*,2}^{-1}(\phi^{-1}_{d^*,j^*,2}(h)): \omega_f = \omega, \mathrm{rk}(\omega_f^{[b]}) - \mathrm{rk}(\omega) = J_{b} - J_0\}}{\# \phi_{d_a}^{-1}(h)} \\
        & \hspace{15pt} - \frac{\# \{\omega' \in \Omega_{\omega, g_{[b]}^*}: \mathrm{rk}(\omega^{[b]}) - \mathrm{rk}(\omega) = J_{b} - J_0\}}{\# \Omega_{\omega, g_{[b]}^*}} \cdot \frac{1}{(\ell(\ell-1))^{\hat{\omega}(h)}} \biggr| \\
        &< 3 \cdot \hat{C}_{T,\ell} \cdot (n \log q)^{-2 m_{n,q} + 4 \log \ell}.
    \end{split}
    \end{align*}
    We then take an average of the expressions above over $\widetilde{(g_a)}_a \in \mathrm{PConf}_{\overline{\lambda}}(\mathcal{P}_2(d^*))$ , use triangle inequality, and recall that $\# \pi_{d^*,j^*,2}^{-1}(\phi_{d^*,j^*,2}^{-1}(h)) = \# \phi_{d_a}^{-1}(h) \cdot \# \mathrm{PConf}_{\overline{\lambda}}(\mathcal{P}_2(d^*))$ to obtain
    \begin{align*}
    \begin{split}
        & \hspace{15pt} \biggl| \frac{\#\left\{(f_a, \widetilde{(g_a)}_a) \in \pi_{d^*,j^*,2}^{-1}(\phi^{-1}_{d^*,j^*,2}(h)): \omega_f = \omega, \mathrm{rk}(\omega_f^{[b]}) - \mathrm{rk}(\omega) = J_{b} - J_0 \right\}}{\# \pi_{d^*,j^*,2}^{-1} (\phi_{d^*,j^*,2}^{-1}(h))} \\
        & \hspace{30pt} - \sum_{\widetilde{(g_a)}_a \in \mathrm{PConf}_{\overline{\lambda}}(\mathcal{P}_2(d^*))} \frac{\# \{\omega^{[b]} \in \Omega_{\omega, g_{[b]}}: \mathrm{rk}(\omega^{[b]}) - \mathrm{rk}(\omega) = J_{b} - J_0\}}{\# \mathrm{PConf}_{\overline{\lambda}}(\mathcal{P}_2(d^*)) \cdot \# \Omega_{\omega, g_{[b]}}} \cdot \frac{1}{(\ell(\ell-1))^{\hat{\omega}(h)}} \biggr| \\
        &< 3 \cdot \hat{C}_{T,\ell} \cdot (n \log q)^{-2 m_{n,q} + 4 \log \ell}.
    \end{split}
    \end{align*}
    Because the set of elements in $\Omega_{\omega, g_{[b]}}$ only depend on the first $b$ many elements of $\widetilde{(g_a)}_a \in \mathrm{PConf}_{\overline{\lambda}}(\mathcal{P}_2(d^*))$, we have
    \begin{align*}
        & \hspace{15pt} \sum_{\widetilde{(g_a)}_a \in \mathrm{PConf}_{\overline{\lambda}}(\mathcal{P}_2(d^*))} \# \{\omega^{[b]} \in \Omega_{\omega, g_{[b]}}: \mathrm{rk}(\omega^{[b]}) - \mathrm{rk}(\omega) = J_{b} - J_0\} \\
        &= \frac{\# \mathrm{PConf}_{\overline{\lambda}}(\mathcal{P}_2(d^*))}{\# \mathrm{PConf}_{b}(\mathcal{P}_2(d^*))} \cdot \sum_{\widetilde{(g_{[b]})_a} \in \mathrm{PConf}_b(\mathcal{P}_2(d^*))} \# \{\omega^{[b]} \in \Omega_{\omega, g_{[b]}}: \mathrm{rk}(\omega^{[b]}) - \mathrm{rk}(\omega) = J_{b} - J_0\}.
    \end{align*}
    Since for every $g_{[b]}$ we have $\# \Omega_{\omega, g_{[b]}} = (\ell(\ell-1))^{\overline{\lambda} + b}$, the following equation holds:
    \begin{align*}
        \begin{split}
            & \hspace{15pt} \sum_{\widetilde{(g_a)}_a \in \mathrm{PConf}_{\overline{\lambda}}(\mathcal{P}_2(d^*))} \frac{\# \{\omega^{[b]} \in \Omega_{\omega, g_{[b]}}: \mathrm{rk}(\omega^{[b]}) - \mathrm{rk}(\omega) = J_{b} - J_0\}}{\# \mathrm{PConf}_{\overline{\lambda}}(\mathcal{P}_2(d^*)) \cdot \# \Omega_{\omega, g_{[b]}}} \\
            &= \frac{\sum_{\widetilde{(g_{[b]})}_a \in \mathrm{PConf}_{b}(\mathcal{P}_2(d^*))} \# \{\omega^{[b]} \in \Omega_{\omega, g^{[b]}}: \mathrm{rk}(\omega^{[b]}) - \mathrm{rk}(\omega) = J_{b} - J_0\}}{\sum_{\widetilde{(g_{[b]})}_a \in \mathrm{PConf}_{b}(\mathcal{P}_2(d^*))} \# \Omega_{\omega, g_{[b]}}}.
        \end{split}
    \end{align*}
    Hence, we obtain the desired relation between the two quantities:
    \begin{align} \label{eqn:symmetric-2-global-local}
    \begin{split}
        & \hspace{15pt} \biggl| \frac{\#\left\{(f_a, \widetilde{(g_a)}_a) \in \pi_{d^*,j^*,2}^{-1}(\phi^{-1}_{d^*,j^*,2}(h)): \omega_f = \omega, \mathrm{rk}(\omega_f^{[b]}) = J_{b}, \; \mathrm{rk}(\omega) = J_0 \right\}}{\# \pi_{d^*,j^*,2}^{-1} (\phi_{d^*,j^*,2}^{-1}(h))} \\
        & \hspace{30pt} - \frac{\sum_{\widetilde{(g_{[b]})}_a \in \mathrm{PConf}_{b}(\mathcal{P}_2(d^*))} \# \{\omega^{[b]} \in \Omega_{\omega, g_{[b]}}: \mathrm{rk}(\omega^{[b]}) - \mathrm{rk}(\omega) = J_{b} - J_0\}}{(\ell(\ell-1))^{\hat{\omega}(h)} \cdot \sum_{\widetilde{(g_{[b]})}_a \in \mathrm{PConf}_{b}(\mathcal{P}_2(d^*))} \# \Omega_{\omega, g_{[b]}}} \biggr| \\
        &< 3 \cdot \hat{C}_{T,\ell} \cdot (n \log q)^{-2 m_{n,q} + 4 \log \ell}.
    \end{split}
    \end{align}
    We will take summations of the inequality above over the set 
    \begin{equation*}
        \widetilde{\Omega_{\overline{h}^*}} := \{\omega \in \Omega_{\overline{h}^*} : \exists f \in \pi_{d^*,j^*,2}^{-1}(\phi_{d^*,j^*,2}^{-1}(h)) \; \textrm{s.t. } \omega_f = \omega \}.
    \end{equation*}
    Using Lemma \ref{lemma:equidistribution}, we obtain $\# \widetilde{\Omega_{\overline{h}^*}} = (\ell(\ell-1))^{\hat{\omega}(h)}$. Hence we obtain the following inequality:
    \begin{align} \label{eqn:symmetric-2-global}
    \begin{split}
        & \hspace{15pt} \biggl| \frac{\sum_{J_0 = -\infty}^{\infty} \#\left\{(f_a, \widetilde{(g_a)}_a) \in \pi_{d^*,j^*,2}^{-1}(\phi^{-1}_{d^*,j^*,2}(h)): \mathrm{rk}(\omega_f^{[b]}) = J_{b}, \; \mathrm{rk}(\omega) = J_0 \right\}}{\# \pi_{d^*,j^*,2}^{-1} (\phi_{d^*,j^*,2}^{-1}(h))} \\
        & \hspace{30pt} - \sum_{\omega \in \widetilde{\Omega_{\overline{h}^*}}} \frac{\sum_{\widetilde{(g_{[b]})}_a \in \mathrm{PConf}_{b}(\mathcal{P}_2(d^*))} \# \{\omega^{[b]} \in \Omega_{\omega, g_{[b]}}: \mathrm{rk}(\omega^{[b]}) - \mathrm{rk}(\omega) = J_{b} - J_0\}}{\# \widetilde{\Omega_{\overline{h}^*}} \cdot \sum_{\widetilde{(g_{[b]})}_a \in \mathrm{PConf}_{b}(\mathcal{P}_2(d^*))} \# \Omega_{\omega, g_{[b]}}} \biggr| \\
        &< 3 \cdot \hat{C}_{T,\ell} \cdot (n \log q)^{-2 m_{n,q} + 4 \log \ell},
    \end{split}
    \end{align}
    which resembles inequality (\ref{eqn:symmetric-2-induction}) when $b = \overline{\lambda}$. It hence remains to show how we can rewrite the second term of the inequality above by suitably defining auxiliary operators $\mathcal{W}_1, \cdots, \mathcal{W}_{b}$.
    
    \medskip

    \textbf{[Base case]}

    \medskip

    We start with the base step. Let $g_1 \in \mathcal{P}_2(d^*)$ be the first coordinate element of $\widetilde{(g_a)}_a \in \mathrm{PConf}_{\overline{\lambda}}(\mathcal{P}_2(d^*))$ (Note $g_{[1]} = g_1$). Given $f \in \phi_{d^*,j^*,2}^{-1}(h)$, we let $\omega_f^{[1]}$ be a character defined as
    \begin{equation*}
        \omega_f^{[1]} := (\chi_{f,v})_{v \in \Sigma_f(\overline{h}^* \cdot g_1)} \in \Omega_{\overline{h}^* \cdot g_1}.
    \end{equation*}
    Our goal is to construct an auxiliary operator $\mathcal{W}_1$ over $\mathbb{Z}_{\geq 0}$ such that
    \begin{align} \label{eqn:symmetric-ind-base}
    \begin{split}
        & \hspace{15pt} \left| \frac{\sum_{J_0 = -\infty}^\infty \# \{f \in \pi_{d^*,j^*,2}^{-1}(\phi_{d^*,j^*,2}^{-1}(h)) : \mathrm{rk}(\omega_f^{[1]}) = J_1, \mathrm{rk}(\omega_f) = J_0\} }{\# \pi_{d^*,j^*,2}^{-1}(\phi_{d^*,j^*,2}^{-1}(h))} - (\mathcal{W}_{1} \delta_h)(J_{1}) \right| \\
        &< B_{T,\ell,q} \cdot (n \log q)^{-2m_{n,q} + 6 \log \ell + 1}.
    \end{split}
    \end{align}
    We follow through the steps we used to obtain the base case for equation (\ref{eqn:symmetric-induction}), but take a different definition of the constants $c_{2, J_1 - J_0}$ (where we recall that $k^* = 2$). As before, we combine Propositions \ref{prop:rank_diff} and \ref{prop:local-Selmer}. In addition, we make the following additional observations regarding the proof of Proposition \ref{prop:rank_diff}.
    \begin{itemize}
        \item If $t_{\omega}(g_1) = 1$, then $\dim_{\mathbb{F}_\ell} V_{g_1} \cap \alpha_{g_1}(\overline{\omega}_f^{[1]}) = 0$ or $1$. Because the twisting data $\alpha_{g_1}$ is injective, we obtain either one of the following two scenarios for each $g_1 \in \mathcal{P}_{2}(d^*)$:
        \begin{itemize}
            \item All ramified maximal isotropic subspaces in $\mathrm{Image}(\alpha_{g_1})$ satisfy $\dim_{\mathbb{F}_\ell} V_{g_1} \cap \alpha_{g_1}(\overline{\omega}_f^{[1]}) = 0$. In this case, for all characters $\omega' \in \Omega_{\omega,g_1}$, we have $\mathrm{rk}(\omega') - \mathrm{rk}(\omega_f) = -1$.
            \item There exists exactly one out of $\ell$ many ramified maximal isotropic subspaces in $\mathrm{Image}(\alpha_{g_1})$ that satisfies $\dim_{\mathbb{F}_\ell} V_{g_1} \cap \alpha_{g_1}(\overline{\omega}_f^{[1]}) = 1$. In this case, for $\ell - 1$ many ramified characters $\omega'$ in $\Omega_{\omega,g_1}$, $\mathrm{rk}(\omega') - \mathrm{rk}(\omega_f) = 0$, whereas for other $(\ell - 1)^2$ many ramified characters in $\Omega_{\omega,g_1}$, we have $\mathrm{rk}(\omega') - \mathrm{rk}(\omega_f) = -1$.
        \end{itemize}
        \item If $t_{\omega}(g_1) = 2$, then $\dim_{\mathbb{F}_\ell} V_{g_1} \cap \alpha_{g_1}(\overline{\omega}_f^{[1]}) = 0$, $1$, or $2$. We obtain either one of the following three scenarios for each $g_1 \in \mathcal{P}_{2}(d^*)$.
        \begin{itemize}
            \item All ramified maximal isotropic subspaces in $\mathrm{Image}(\alpha_{g_1})$ satisfy $\dim_{\mathbb{F}_\ell} V_{g_1} \cap \alpha_{g_1}(\overline{\omega}_f^{[1]}) = 0$. In this case, for all characters $\omega' \in \Omega_{\omega,g_1}$, we have $\mathrm{rk}(\omega') - \mathrm{rk}(\omega_f) = 0$.
            \item There exist some ramified maximal isotropic subspaces in $\mathrm{Image}(\alpha_{g_1})$ that satisfy $\dim_{\mathbb{F}_\ell} V_{g_1} \cap \alpha_{g_1}(\overline{\omega}_f') = 1$, whereas for others the intersection has dimension $0$. In this case, we have $\mathrm{rk}(\omega') - \mathrm{rk}(\omega_f) = 0$ or $1$.
            \item There exists exactly one out of $\ell$ many ramified maximal isotropic subspaces in $\mathrm{Image}(\alpha_{g_1})$ that satisfies $\dim_{\mathbb{F}_\ell} V_{g_1} \cap \alpha_{g_1}(\overline{\omega}_f') = 2$. Here, we use the fact that $\alpha_v$ is an injection. In this case, for $\ell - 1$ many ramified characters $\omega'$ in $\Omega_{\omega,g_1}$, $\mathrm{rk}(\omega') - \mathrm{rk}(\omega_f) = 2$, whereas for other $(\ell - 1)^2$ many ramified characters in $\Omega_{\omega,g_1}$, we have $\mathrm{rk}(\omega') - \mathrm{rk}(\omega_f) = 0$ or $1$.
        \end{itemize}
    \end{itemize}
    We define the following three parameters, which depends on our fixed choice of $\omega$. The notations we use for these parameters aligns with the parameters which appear in Remark \ref{remark:motivation-Markov}.
    \begin{align}
        \begin{split}
            s_{1,\omega} &:= \frac{\sum_{\substack{g_1 \in \mathcal{P}_2(d^*) \\ t_{\omega}(g_1) = 1}}\#\{\omega' \in \Omega_{\omega,g_1} : \mathrm{rk}(\omega') - \mathrm{rk}(\omega) = 0  \}}{\sum_{\substack{g_1 \in \mathcal{P}_2(d^*) \\ t_{\omega}(g_1) = 1}} \# \Omega_{\omega,g_1}}, \\
            s_{2,\omega} &:= \frac{\sum_{\substack{g_1 \in \mathcal{P}_2(d^*) \\ t_{\omega}(g_1) = 2}}\#\{\omega' \in \Omega_{\omega,g_1} : \mathrm{rk}(\omega') - \mathrm{rk}(\omega) = 0  \}}{\sum_{\substack{g_1 \in \mathcal{P}_2(d^*) \\ t_{\omega}(g_1) = 2}} \# \Omega_{\omega,g_1}}, \\
            s_{3,\omega} &:= \frac{\sum_{\substack{g_1 \in \mathcal{P}_2(d^*) \\ t_{\omega}(g_1) = 2}}\#\{\omega' \in \Omega_{\omega,g_1} : \mathrm{rk}(\omega') - \mathrm{rk}(\omega) = 1  \}}{\sum_{\substack{g_1 \in \mathcal{P}_2(d^*) \\ t_{\omega}(g_1) = 2}} \# \Omega_{\omega,g_1}}.
        \end{split}
    \end{align}
    Then we can define the constants $c_{2, J_1 \; - \;J_0, \omega}$ as follows:
    \begin{equation}
        c_{2,J_1 \; - \; J_0, \omega} := \begin{cases}
            1 - (\ell + 1)\ell^{-J_0} + \ell^{-2J_0 + 1} &\text{ if } J_1 - J_0 = -2, \\
            (1-s_{1,\omega}) \cdot (\ell + 1) \cdot (\ell^{-J_0} - \ell^{-2J_0}) &\text{ if } J_1 - J_0 = -1, \\
            s_{1,\omega} \cdot (\ell + 1) \cdot (\ell^{-J_0} - \ell^{-2J_0}) + s_{2,\omega} \cdot \ell^{-2J_0} &\text{ if } J_1 - J_0 = 0, \\
            s_{3,\omega} \cdot \ell^{-2J_0} &\text{ if } J_1 - J_0 = 1, \\
            (1 - s_{2,\omega} - s_{3,\omega}) \cdot \ell^{-2 J_0} &\text{ if } J_1 - J_0 = 2, \\
            0 &\text{ otherwise}.
        \end{cases}
    \end{equation}
    Then we obtain the following inequality by using Propositions \ref{prop:rank_diff} and \ref{prop:local-Selmer}:
    \begin{align} \label{eqn:symmetric-basecase-1}
    \begin{split}
        & \hspace{15pt} \left| \frac{\sum_{g_1 \in \mathcal{P}_2(d^*)} \# \{\omega' \in \Omega_{\omega, g_1}: \mathrm{rk}(\omega') - \mathrm{rk}(\omega) = J_1 - J_0\}}{\sum_{g_1 \in \mathcal{P}_2(d^*)} \# \Omega_{\omega,g_1}} - c_{2, J_1 \; - \; J_0, \omega} \right| \\
        &< 16 \cdot \ell^4 \cdot D_{T,\ell,q} \cdot (n \log q)^{-2 m_{n,q} + 6 \log \ell + 1}.
    \end{split}
    \end{align}
    To obtain the base case for inequality (\ref{eqn:symmetric-2-induction}), we take the average of inequality (\ref{eqn:symmetric-basecase-1}) over the set of characters $\omega \in \widetilde{\Omega_{\overline{h}^*}}$. This yields for each $J_1$:
    \begin{align*}
    \begin{split}
        & \hspace{15pt} \left| \sum_{\omega \in \widetilde{\Omega_{\overline{h}^*}}}\frac{\sum_{g_1 \in \mathcal{P}_2(d^*)} \# \{\omega' \in \Omega_{\omega, g_1}: \mathrm{rk}(\omega') - \mathrm{rk}(\omega) = J_1 - J_0\}}{\# \widetilde{\Omega_{\overline{h}^*}} \cdot \sum_{g_1 \in \mathcal{P}_2(d^*)} \# \Omega_{\omega,g_1}} - \sum_{\omega \in \widetilde{\Omega_{\overline{h}^*}}} \frac{c_{2,J_1 \; - \; J_0, \omega}}{\# \widetilde{\Omega_{\overline{h}^*}}} \right| \\
        &< 16 \cdot \ell^4 \cdot D_{T,\ell,q} \cdot (n \log q)^{-2 m_{n,q} + 6 \log \ell + 1}.
    \end{split}
    \end{align*}
    
    Given $\omega \in \Omega_{\overline{h}^*}$, let $\tilde{\delta}_\omega: \mathbb{Z}_{\geq 0} \to [0,1]$ be a function defined as
    \begin{equation}
        \tilde{\delta}_\omega(J) = \begin{cases}
            1 / \# \widetilde{\Omega_{\overline{h}^*}} &\text{ if } J = \mathrm{rk}(\omega) \\
            0 &\text{ otherwise}.
        \end{cases}
    \end{equation}
    Then the function $\sum_{\omega \in \widetilde{\Omega_{\overline{h}^*}}} \tilde{\delta}_\omega$ is a probability distribution over $\mathbb{Z}_{\geq 0}$. Furthermore,
    Lemma \ref{lemma:equidistribution} ensures that
    \begin{equation} \label{eqn:symmetric-2-base-2}
        \sup_{J \in \mathbb{Z}_{\geq 0}} \left| \delta_h - \sum_{\omega \in \widetilde{\Omega_{\overline{h}^*}}} \tilde{\delta}_\omega \right| < \hat{C}_{T,\ell} \cdot (n \log q)^{-2m_{n,q} + 4 \log \ell},
    \end{equation}
    For each $\omega \in \widetilde{\Omega_{\overline{h}^*}}$, we denote by $\mathcal{W}_{\omega}$ an auxiliary Markov operator over $\mathbb{Z}_{\geq 0}$ whose transition probabilities from $J_0 = \mathrm{rk}(\omega)$ to $J_1$ are given by constants $c_{2, J_1 \; - \; J_0, \omega}$. Then we have
    \begin{equation}
        \sum_{\omega \in \widetilde{\Omega_{\overline{h}^*}}} \frac{c_{2,J_1 \; - \; J_0, \omega}}{\# \widetilde{\Omega_{\overline{h}^*}}} = \sum_{\omega \in \widetilde{\Omega_{\overline{h}^*}}} (\mathcal{W}_\omega \delta_\omega)(J_1).
    \end{equation}
    Because $\widetilde{\Omega_{\overline{h}^*}}$ is a finite set, and the probability distribution $\sum_{\omega \in \widetilde{\Omega_{\overline{h}^*}}} \tilde{\delta}_\omega$ is nonzero only over a finite set of integers $J$, we can always construct an auxiliary operator $\mathcal{W}_{1}$ over $\mathbb{Z}_{\geq 0}$ which satisfies the equation
    \begin{equation}
        \mathcal{W}_{1} \sum_{\omega \in \widetilde{\Omega_{\overline{h}^*}}} \tilde{\delta}_\omega = \sum_{\omega \in \widetilde{\Omega_{\overline{h}^*}}} \mathcal{W}_\omega \delta_\omega.
    \end{equation}
    Applying the operator $\mathcal{W}_1$ to (\ref{eqn:symmetric-2-base-2}) gives
    \begin{equation} \label{eqn:W_omega-basecase-defn}
        \sup_{J \in \mathbb{Z}_{\geq 0}} \left|(\mathcal{W}_{1} \delta_h)(J) - \left(\mathcal{W}_1 \sum_{\omega \in \widetilde{\Omega_{\overline{h}^*}}} \tilde{\delta}_\omega \right)(J)\right| < 5 \cdot \hat{C}_{T,\ell} \cdot (n \log q)^{-2m_{n,q} + 4 \log \ell}.
    \end{equation}
    The constant $5$ is multiplied because $\mathcal{W}_1$ allows only at $5$ states $J-2, \cdots, J+2$ to move to state $J$. Using this auxiliary operator, we obtain the following inequality for all $J_1 \in \mathbb{Z}_{\geq 0}$:
    \begin{align*}
    \begin{split}
        & \hspace{15pt} \left| \sum_{\omega \in \widetilde{\Omega_{\overline{h}^*}}}\frac{\sum_{g_1 \in \mathcal{P}_2(d^*)} \# \{\omega' \in \Omega_{\omega, g_1}: \mathrm{rk}(\omega') - \mathrm{rk}(\omega) = J_1 - J_0\}}{\# \widetilde{\Omega_{\overline{h}^*}} \cdot \sum_{g_1 \in \mathcal{P}_2(d^*)} \# \Omega_{\omega,g_1}} - (\mathcal{W}_{1} \delta_h)(J_1) \right| \\
        &< (16 \cdot \ell^4 \cdot D_{T,\ell,q} + 5 \cdot \hat{C}_{T,\ell}) \cdot (n \log q)^{-2 m_{n,q} + 6 \log \ell + 1}.
    \end{split}
    \end{align*}
    We combine this inequality with (\ref{eqn:symmetric-2-global}) to prove the base case for (\ref{eqn:symmetric-2-induction}).

    It is straightforward to check that $\mathbb{E}[\ell^{\mathcal{W}_{1} \delta_h}]$ achieves its maximum value if for all $\omega \in \widetilde{\Omega_{\overline{h}^*}}$ the constants $c_{2, J_{\overline{\lambda} \; - \; J_0, \omega}}$ satisfy the conditions that $s_{1,\omega} = 1$ and $s_{2, \omega} = s_{3, \omega} = 0$. But using such constraints on the variables $s_{1, \omega}, s_{2, \omega}, s_{3, \omega}$ gives transition probabilities for the Markov operator $\mathcal{V}$ over $\mathbb{Z}_{\geq 0}$. We hence obtain the desired inequality $\mathbb{E}[\ell^{\mathcal{W}_{1} \delta_h}] \leq \mathbb{E}[\ell^{\mathcal{V} \delta_h}]$.

    \medskip 

    \textbf{[Induction case]}

    \medskip

    We proceed to proving the case where $b > 1$. We make the following induction hypothesis: suppose inequality (\ref{eqn:symmetric-2-induction}) holds for $b$. In particular, we require the following inequality to hold for every $J_b \in \mathbb{Z}_{\geq 0}$:
    \begin{align} \label{eqn:symmetric-ind-hypothesis}
    \begin{split}
        & \hspace{15pt} \left| \frac{\sum_{J_0 = -\infty}^\infty \# \{f \in \pi_{d^*,j^*,2}^{-1}(\phi_{d^*,j^*,2}^{-1}(h)) : \mathrm{rk}(\omega_f^{[b]}) = J_{b}, \mathrm{rk}(\omega_f) = J_0\} }{\# \pi_{d^*,j^*,2}^{-1}(\phi_{d^*,j^*,2}^{-1}(h))} - (\mathcal{W}_{b} \cdots \mathcal{W}_{1} \delta_h)(J_b) \right| \\
        &< b \cdot B_{T,\ell,q} \cdot (n \log q)^{-2m_{n,q} + 6 \log \ell + 1}.
    \end{split}
    \end{align}

    To proceed, we fix a choice of $\widetilde{(g_{[b]})}_a \in \mathrm{PConf}_{b}(\mathcal{P}_2(d^*))$. There is a canonical projection map 
    \begin{equation*}
        \mathrm{PConf}_{b+1}(\mathcal{P}_2(d^*)) \to \mathrm{PConf}_{b}(\mathcal{P}_2(d^*))
    \end{equation*}
    which forgets $g_{b+1}$, and the fiber of $\widetilde{(g_{[b]})}_a$ is isomorphic to the set $\mathcal{P}_2(d^*) \setminus \{g_1, g_2, \cdots, g_{b}\}$. With this choice made, we can fix a choice of $\omega^{[b]} \in \Omega_{\omega,g_{[b]}}$.

    Proceeding as the base case, we define the following three parameters.
    \begin{align}
        \begin{split}
            s_{1,\omega^{[b]}} &:= \frac{\sum_{\substack{g_{b+1} \in \mathcal{P}_2(d^*) \setminus \{g_1, \cdots, g_{b}\} \\ t_{\omega^{[b]}}(g_{b+1}) = 1}}\#\{\omega^{[b+1]} \in \Omega_{\omega^{[b]},g_{b+1}} : \mathrm{rk}(\omega^{[b+1]}) - \mathrm{rk}(\omega^{[b]}) = 0  \}}{\sum_{\substack{g_{b+1} \in \mathcal{P}_2(d^*) \setminus \{g_1, \cdots, g_{b}\} \\ t_{\omega^{[b]}}(g_{b+1}) = 1}} \# \Omega_{\omega^{[b]},g_{b+1}}}, \\
            s_{2,\omega^{[b]}} &:= \frac{\sum_{\substack{g_{b+1} \in \mathcal{P}_2(d^*) \setminus \{g_1, \cdots, g_{b}\} \\ t_{\omega^{[b]}}(g_{b+1}) = 2}}\#\{\omega^{[b+1]} \in \Omega_{\omega^{[b]},g_{b+1}} : \mathrm{rk}(\omega^{[b+1]}) - \mathrm{rk}(\omega^{[b]}) = 0  \}}{\sum_{\substack{g_{b+1} \in \mathcal{P}_2(d^*) \setminus \{g_1, \cdots, g_{b}\} \\ t_{\omega^{[b]}}(g_{b+1}) = 2}} \# \Omega_{\omega^{[b]},g_{b+1}}}, \\
            s_{3,\omega^{[b]}} &:= \frac{\sum_{\substack{g_{b+1} \in \mathcal{P}_2(d^*) \setminus \{g_1, \cdots, g_{b}\} \\ t_{\omega^{[b]}}(g_{b+1}) = 2}}\#\{\omega^{[b+1]} \in \Omega_{\omega^{[b]},g_{b+1}} : \mathrm{rk}(\omega^{[b+1]}) - \mathrm{rk}(\omega^{[b]}) = 1  \}}{\sum_{\substack{g_{b+1} \in \mathcal{P}_2(d^*) \setminus \{g_1, \cdots, g_{b}\} \\ t_{\omega^{[b]}}(g_{b+1}) = 2}} \# \Omega_{\omega^{[b]},g_{b+1}}}.
        \end{split}
    \end{align}
    Then we can define the constants $c_{2, J_{b+1} - J_b, \omega^{[b]}}$ as follows:
    \begin{equation}
        c_{2,J_{b+1} - J_{b}, \omega^{[b]}} := \begin{cases}
            1 - (\ell + 1)\ell^{-J_{b}} + \ell^{-2J_{b} + 1} &\text{ if } J_{b+1} - J_{b} = -2, \\
            (1-s_{1,\omega^{[b]}}) \cdot (\ell + 1) \cdot (\ell^{-J_{b}} - \ell^{-2J_{b}}) &\text{ if } J_{b+1} - J_{b} = -1, \\
            s_{1,\omega^{[b]}} \cdot (\ell + 1) \cdot (\ell^{-J_{b}} - \ell^{-2J_{b}}) + s_{2,\omega^{[b]}} \cdot \ell^{-2J_{b}} &\text{ if } J_{b+1} - J_{b} = 0, \\
            s_{3,\omega^{[b]}} \cdot \ell^{-2J_{b}} &\text{ if } J_{b+1} - J_{b} = 1, \\
            (1 - s_{2,\omega^{[b]}} - s_{3,\omega^{[b]}}) \cdot \ell^{-2 J_{b}} &\text{ if } J_{b+1} - J_{b} = 2, \\
            0 &\text{ otherwise}.
        \end{cases}
    \end{equation}
    Using these constants, the following inequality holds for a fixed $\omega^{[b]}$ and $g_{[b]}$ as a consequence of Propositions \ref{prop:rank_diff} and \ref{prop:local-Selmer}, and the fact that $\omega(g_{[b+1]}) < 2 m_{n,q}$:
    \begin{align} \label{eqn:symmetric-2-ind-1}
        \begin{split}
            & \hspace{15pt} \biggl| \frac{\sum_{g_{b+1} \in \mathcal{P}_2(d^*) \setminus \{g_1, \cdots, g_{b}\}} \#\{ \omega^{[b+1]} \in \Omega_{\omega^{[b]},g_{b+1}} : \mathrm{rk}(\omega^{[b+1]}) - \mathrm{rk}(\omega^{[b]}) = J_{b+1} - J_{b} \} }{\sum_{g_{b+1} \in \mathcal{P}_2(d^*) \setminus \{g_1, \cdots, g_{b}\}} \# \Omega_{\omega^{[b]}, g_{b+1}}} \\
            & \hspace{30pt} - c_{2,J_{b+1} - J_b, \omega^{[b]}} \biggr| < 16 \cdot \ell^4 \cdot D_{T,\ell,q} \cdot (n \log q)^{-2 m_{n,q} + 6 \log \ell + 1}.
        \end{split}
    \end{align}
    Using Lemma \ref{lemma:equidistribution}, we have $\# \Omega_{\omega, g_{[b]}} = (\ell(\ell-1))^{b}$. Hence by taking an average of inequality (\ref{eqn:symmetric-2-ind-1}) over the set $\omega^{[b]} \in \Omega_{\omega, g_{[b]}}$, we obtain for a fixed $g_{[b]}$:
    \begin{align*}
        \begin{split}
            & \hspace{15pt} \biggl| \sum_{\omega^{[b]} \in \Omega_{\omega,g_{[b]}}}\frac{\sum_{g_{b+1} \in \mathcal{P}_2(d^*) \setminus \{g_1, \cdots, g_{b}\}} \#\{ \omega^{[b+1]} \in \Omega_{\omega^{[b]},g_{b+1}} : \mathrm{rk}(\omega^{[b+1]}) - \mathrm{rk}(\omega^{[b]}) = J_{b+1} - J_{b} \} }{ \sum_{g_{b+1} \in \mathcal{P}_2(d^*) \setminus \{g_1, \cdots, g_{b}\}} \# \Omega_{\omega, g_{[b+1]}}} \\
            & \hspace{30pt} - \sum_{\omega^{[b]} \in \Omega_{\omega,g_{[b]}}} \frac{c_{2,J_{b+1} - J_{b}, \omega^{[b]}}}{\# \Omega_{\omega,g_{[b]}}} \biggr| < 16 \cdot \ell^4 \cdot D_{T,\ell,q} \cdot (n \log q)^{-2 m_{n,q} + 6 \log \ell + 1}.
        \end{split}
    \end{align*}
    Notice that for all $g_{b}$ and $\omega^{[b]} \in \Omega_{\omega, g_{[b]}}$, we have $\# \Omega_{\omega, g_{[b]}} \cdot \# \Omega_{\omega^{[b]}, g_{b+1}} = \# \Omega_{\omega, g_{[b+1]}}$. This allows us to rewrite the denominator term of the expression above as displayed above.

    Next, we take an average of the inequality above over the set $\widetilde{(g_{[b]})}_a \in \mathrm{PConf}_{b}(\mathcal{P}_2(d^*))$. This leads us to obtain for any fixed $\omega \in \Omega_{\overline{h}^*}$ (where we abbreviate the summand for $g_{b+1}$):
    \begin{align*}
        \begin{split}
            & \hspace{15pt} \biggl| \sum_{\widetilde{(g_{[b]})}_a \in \mathrm{PConf}_{b}(\mathcal{P}_2(d^*))}\sum_{\omega^{[b]} \in \Omega_{\omega,g_{[b]}}}\frac{\sum_{g_{b+1}} \#\{ \omega^{[b+1]} \in \Omega_{\omega^{[b]},g_{b+1}} : \mathrm{rk}(\omega^{[b+1]}) - \mathrm{rk}(\omega^{[b]}) = J_{b+1} - J_{b} \} }{\sum_{\widetilde{(g_{[b]})}_a \in \mathrm{PConf}_{b}(\mathcal{P}_2(d^*))} \sum_{g_{b+1}} \# \Omega_{\omega, g_{[b+1]}}} \\
            & \hspace{30pt} - \sum_{\widetilde{(g_{[b]})}_a \in \mathrm{PConf}_{b}(\mathcal{P}_2(d^*))} \sum_{\omega^{[b]} \in \Omega_{\omega,g_{[b]}}} \frac{c_{2,J_{b+1} - J_{b}, \omega^{[b]}}}{ \sum_{\widetilde{(g_{[b]})}_a \in \mathrm{PConf}_{b}(\mathcal{P}_2(d^*))}\# \Omega_{\omega,g_{[b]}}} \biggr| \\
            &< 16 \cdot \ell^4 \cdot D_{T,\ell,q} \cdot (n \log q)^{-2 m_{n,q} + 6 \log \ell + 1}.
        \end{split}
    \end{align*}
    By definition, we can rewrite the first term of the above inequality as follows for any fixed $\omega \in \Omega_{\overline{h}^*}$:
    \begin{align*}
        \begin{split}
            & \hspace{15pt} \biggl| \frac{\sum_{\widetilde{(g_{[b+1]})}_a \in \mathrm{PConf}_{b+1}(\mathcal{P}_2(d^*))} \# \{\omega^{[b+1]} \in \Omega_{\omega, g_{[b+1]}}: \mathrm{rk}(\omega^{[b+1]}) - \mathrm{rk}(\omega) = J_{b+1} - J_0\}}{ \sum_{\widetilde{(g_{[b+1]})}_a \in \mathrm{PConf}_{b+1}(\mathcal{P}_2(d^*))} \# \Omega_{\omega, g_{[b+1]}}} \\
            & \hspace{30pt} - \sum_{\widetilde{(g_{[b]})}_a \in \mathrm{PConf}_{b}(\mathcal{P}_2(d^*))} \sum_{\omega^{[b]} \in \Omega_{\omega,g_{[b]}}} \frac{c_{2,J_{b+1} - J_{b}, \omega^{[b]}}}{ \sum_{\widetilde{(g_{[b]})_a} \in \mathrm{PConf}_{b}(\mathcal{P}_2(d^*))}\# \Omega_{\omega,g_{[b]}}} \biggr| \\
            &< 16 \cdot \ell^4 \cdot D_{T,\ell,q} \cdot (n \log q)^{-2 m_{n,q} + 6 \log \ell + 1}.
        \end{split}
    \end{align*}
    Finally, we take an average of the inequality above over the set $\omega \in \widetilde{\Omega_{\overline{h}^*}}$ to obtain for each $J_{b+1}$:
    \begin{align} \label{eqn:symmetric-2-ind-2}
        \begin{split}
            & \hspace{15pt} \biggl| \sum_{\omega \in \widetilde{\Omega_{\overline{h}^*}}} \frac{\sum_{\widetilde{(g_{[b+1]})}_a \in \mathrm{PConf}_{b+1}(\mathcal{P}_2(d^*))} \# \{\omega^{[b+1]} \in \Omega_{\omega, g_{[b+1]}}: \mathrm{rk}(\omega^{[b+1]}) - \mathrm{rk}(\omega) = J_{b+1} - J_0\}}{ \# \widetilde{\Omega_{\overline{h}^*}} \cdot  \sum_{\widetilde{(g_{[b+1]})}_a \in \mathrm{PConf}_{b+1}(\mathcal{P}_2(d^*))} \# \Omega_{\omega, g_{[b+1]}}} \\
            & \hspace{30pt} - \sum_{\omega \in \widetilde{\Omega_{\overline{h}^*}}} \; \sum_{\widetilde{(g_{[b]})}_a \in \mathrm{PConf}_{b}(\mathcal{P}_2(d^*))} \sum_{\omega^{[b]} \in \Omega_{\omega,g_{[b]}}} \frac{c_{2,J_{b+1} - J_{b}, \omega^{[b]}}}{ \# \widetilde{\Omega_{\overline{h}^*}} \cdot \sum_{\widetilde{(g_{[b]})_a} \in \mathrm{PConf}_{b}(\mathcal{P}_2(d^*))}\# \Omega_{\omega,g_{[b]}}} \biggr| \\
            &< 16 \cdot \ell^4 \cdot D_{T,\ell,q} \cdot (n \log q)^{-2 m_{n,q} + 6 \log \ell + 1}.
        \end{split}
    \end{align}
    We denote by $\delta_{b}: \mathbb{Z}_{\geq 0} \to [0,1]$ the probability distribution
    \begin{equation}
        \delta_{b}(J) := \frac{\#\{f \in \pi_{d^*,j^*,2}^{-1}(\phi_{d^*,j^*,2}^{-1}(h)) : \mathrm{rk}(\omega_f^{[b]}) = J\}}{\# \pi_{d^*,j^*,2}^{-1}(\phi_{d^*,j^*,2}^{-1}(h))}.
    \end{equation}
    We note that by the induction hypothesis,
    \begin{equation}
        \sup_{J \in \mathbb{Z}_{\geq 0}} \left| \delta_b(J) - (\mathcal{W}_b \cdots \mathcal{W}_1 \delta_h)(J) \right| < b \cdot B_{T,\ell,q} \cdot (n \log q)^{-2m_{n,q} + 6 \log \ell + 1}.
    \end{equation}
    For each $\omega^{[b]} \in \Omega_{\omega, g_{[b]}}$ given some $\omega \in \widetilde{\Omega_{\overline{h}^*}}$ and $g_{[b]}$, we denote by $\tilde{\delta}_{\omega^{[b]}}: \mathbb{Z}_{\geq 0} \to [0,1]$ a function such that
    \begin{equation}
        \tilde{\delta}_{\omega^{[b]}}(J) = \begin{cases} 
        \frac{1}{\# \widetilde{\Omega_{\overline{h}^*}} \cdot \# \mathrm{PConf}_b(\mathcal{P}_2(d^*)) \cdot \# \Omega_{\omega, g_{[b]}}} &\text{ if } J = \mathrm{rk}(\omega^{[b]}) \\
        0 &\text{ otherwise}.
        \end{cases}
    \end{equation}
    Note that the following function is a probability distribution over $\mathbb{Z}_{\geq 0}$:
    \begin{equation}
        \delta_b^+ := \sum_{\omega \in \widetilde{\Omega_{\overline{h}^*}}} \; \sum_{\widetilde{(g_{[b]})}_a \in \mathrm{PConf}_b(\mathcal{P}_2(d^*))} \; \sum_{\omega^{[b]} \in \Omega_{\omega, g_{[b]}}} \tilde{\delta}_{\omega^{[b]}}.
    \end{equation}
    Then (\ref{eqn:symmetric-2-global}) can be rephrased as
    \begin{align}
    \begin{split} \label{eqn:symmetric-2-ind-rephrase}
        \sup_{J \in \mathbb{Z}_{\geq 0}} \left|\delta_b(J) - \delta_b^+(J) \right| < 3 \cdot \hat{C}_{T,\ell} \cdot (n \log q)^{-2 m_{n,q} + 4 \log \ell}.
    \end{split}
    \end{align}
    For each $\omega^{[b]} \in \Omega_{\omega, g_{[b]}}$ given some $\omega \in \widetilde{\Omega_{\overline{h}^*}}$ and $g_{[b]}$, we also denote by $\mathcal{W}_{\omega^{[b]}}$ the auxiliary Markov operator whose transition probabilities from $J_b := \mathrm{rk}(\omega^{[b]})$ to $J_{b+1}$ are given by constants $c_{2, J_{b+1} - J_b, \omega^{[b]}}$. Then the second term from (\ref{eqn:symmetric-2-ind-2}) can be written as
    \begin{align*}
        & \hspace{15pt} \sum_{\omega \in \widetilde{\Omega_{\overline{h}^*}}} \; \sum_{\widetilde{(g_{[b]})}_a \in \mathrm{PConf}_{b}(\mathcal{P}_2(d^*))} \sum_{\omega^{[b]} \in \Omega_{\omega,g_{[b]}}} \frac{c_{2,J_{b+1} - J_{b}, \omega^{[b]}}}{ \# \widetilde{\Omega_{\overline{h}^*}} \cdot \sum_{\widetilde{(g_{[b]})_a} \in \mathrm{PConf}_{b}(\mathcal{P}_2(d^*))}\# \Omega_{\omega,g_{[b]}}} \\
        &= \sum_{\omega \in \widetilde{\Omega_{\overline{h}^*}}} \; \sum_{\widetilde{(g_{[b]})}_a \in \mathrm{PConf}_b(\mathcal{P}_2(d^*))} \; \sum_{\omega^{[b]} \in \Omega_{\omega, g_{[b]}}} (\mathcal{W}_{\omega^{[b]}} \tilde{\delta}_{\omega^{[b]}})(J_{b+1}).
    \end{align*}
    Using the finiteness of $\widetilde{\Omega_{\overline{h}^*}}$ and the fact that $\delta_b^+$ satisfies the condition that $\delta_b^+(J) \neq 0$ over a finite set of integers $J$, we can construct an auxiliary Markov operator $\mathcal{W}_{b+1}$ over $\mathbb{Z}_{\geq 0}$ which satisfies the equation
    \begin{equation}
        (\mathcal{W}_{b+1} \delta_b^+)(J) = \sum_{\omega \in \widetilde{\Omega_{\overline{h}^*}}} \; \sum_{\widetilde{(g_{[b]})}_a \in \mathrm{PConf}_b(\mathcal{P}_2(d^*))} \; \sum_{\omega^{[b]} \in \Omega_{\omega, g_{[b]}}} (\mathcal{W}_{\omega^{[b]}} \tilde{\delta}_{\omega^{[b]}})(J).
    \end{equation}
    In particular, the operator $\mathcal{W}_{b+1}$ is obtained by taking averages of transition probabilities of $\mathcal{W}_{\omega^{[b]}}$. Applying the operator $\mathcal{W}_{b+1}$ to (\ref{eqn:symmetric-2-ind-rephrase}), we obtain
    \begin{equation}
        \sup_{J \in \mathbb{Z}_{\geq 0}} |\mathcal{W}_{b+1} \delta_b(J) - \mathcal{W}_{b+1} \delta_b^+(J)| < 15 \cdot \hat{C}_{T,\ell} \cdot (n \log q)^{-2m_{n,q} + 4 \log \ell}.
    \end{equation}
    We note that the constant changes from $3$ to $15$ because $\mathcal{W}_{b+1}$ allows only at $5$ states $J-2, J-1, \cdots, J+2$ to move to state $J$.
    
    Substituting the above relation to (\ref{eqn:symmetric-2-ind-2}) and using (\ref{eqn:symmetric-2-global}) gives
    \begin{align}
    \begin{split}
        & \hspace{15pt} \biggl| \frac{\sum_{J_0 = -\infty}^{\infty} \#\left\{(f_a, \widetilde{(g_a)}_a) \in \pi_{d^*,j^*,2}^{-1}(\phi^{-1}_{d^*,j^*,2}(h)): \mathrm{rk}(\omega_f^{[b+1]}) = J_{b}, \; \mathrm{rk}(\omega) = J_0 \right\}}{\# \pi_{d^*,j^*,2}^{-1} (\phi_{d^*,j^*,2}^{-1}(h))} - \mathcal{W}_{b+1} \delta_b(J) \biggr| \\
        &< (16 \cdot \ell^4 \cdot D_{T,\ell,q} + 18 \cdot \hat{C}_{T,\ell})\cdot (n \log q)^{-2 m_{n,q} + 6 \log \ell + 1}.
    \end{split}
    \end{align}
    We then use the induction hypothesis on $\delta_b$ to obtain (\ref{eqn:symmetric-2-induction}) for $b+1$. 
    \begin{align}
    \begin{split}
        & \hspace{15pt} \biggl| \frac{\sum_{J_0 = -\infty}^\infty \# \{f \in \pi_{d^*,j^*,2}^{-1}(\phi_{d^*,j^*,2}^{-1}(h)) : \mathrm{rk}(\omega_f^{[b+1]}) = J_{b+1}, \mathrm{rk}(\omega_f) = J_0\} }{\# \pi_{d^*,j^*,2}^{-1}(\phi_{d^*,j^*,2}^{-1}(h))} \\
        & \hspace{30pt} - (\mathcal{W}_{b} \cdots \mathcal{W}_{1} \delta_h)(J_b) \biggr| < (b+1) \cdot B_{T,\ell,q} \cdot (n \log q)^{-2m_{n,q} + 6 \log \ell + 1}.
    \end{split}
    \end{align}
    The upper bound on $\mathbb{E}[\ell^{\mathcal{W}_{b+1} \delta_h}]$ can be computed from setting $s_{1, \omega^{[b]}} = 1$ and $s_{2, \omega^{[b]}} = s_{3, \omega^{[b]}} = 0$ for all $\omega^{[b]} \in \Omega_{\omega, g_{[b]}}$ given any $\omega \in \widetilde{\Omega_{\overline{h}^*}}$ and $\widetilde{(g_{[b]})}_a \in \mathrm{PConf}_b(\mathcal{P}_2(d^*))$. Such constraints on the variables gives transition probabilities for $\mathcal{V}$, which gives $\mathbb{E}[\ell^{\mathcal{W}_{b+1} \delta_h}] \leq \mathbb{E}[\ell^{\mathcal{V} \delta_h}]$.

    The statement of the proposition then follows from repeating the induction step up until $b = \overline{\lambda}$.
\end{proof}

With both propositions obtained, we prove the following analogue of \cite[Proposition 6.11]{park2022prime} which computes the distribution of dimensions of global Selmer structures $\mathrm{Sel}(T, \alpha : \overline{\chi}_f)$ over the set of polynomials $\hat{F}_{(n,N),(w,w')}(\mathbb{F}_q)$.

\begin{proposition} \label{prop:major_contribution-symmetric}
    Fix integers $n> N$, $w > w'$ such that $w < 2m_{n,q}$, $w' = (1-\epsilon) w$ for some small enough $0 < \epsilon < 1$, and $m_{n,q} > \max(e^{e^e}, \sum_{v \in \Sigma_T} \deg v, 6 \log \ell + 2)$.

    Let $T$ be a symmetric admissible $\mathrm{Gal}(\overline{K}/K)$-module with a twisting data $\alpha := (\alpha_v)_v$. Then there exist fixed constants $\tilde{B}_{T, \ell, q}$ depending only on $T, \ell, q$ and $\gamma_T$ depending only on $T$ such that
    \begin{align}
    \begin{split}
        & \hspace{15pt} \left| \frac{\# \{f \in \hat{F}_{(n,N),(w,w')}(\mathbb{F}_q) : \dim_{\mathbb{F}_\ell} \mathrm{Sel}(T, \alpha : \overline{\chi}_f) = r\}}{\# \hat{F}_{(n,N),(w,w')}(\mathbb{F}_q)} - \prod_{k=1}^\infty \frac{1}{1 + \ell^{-k}} \prod_{i=1}^j \frac{1}{\ell^i-1} \right| \\
        &< \left(2 m_{n,q} \cdot B_{T,\ell,q} \cdot (n \log q)^{- m_{n,q}} + \ell^{2 \epsilon w  + 4} \cdot D_{T,\ell,q} \cdot c_T \cdot \gamma_T^{(\delta(\mathcal{P}_1) - \epsilon)(1-\epsilon)w} \right) \\
        &\hspace{15pt} + (1 - \delta(\mathcal{P}_1))^{(1-\delta(\mathcal{P}_1) + \epsilon)(1-\epsilon)w},
    \end{split}    
    \end{align}
    where $D_{T,\ell,q}$ is the constant defined in Proposition \ref{prop:local-Selmer}.
\end{proposition}
\begin{proof}
    We proceed as in the proof of Proposition \ref{prop:major_contribution-alternating} and \cite[Proposition 6.11]{park2022prime}, but the part where we use two types of Markov operators $M_\ell$ and $\mathcal{V}$ is different. As in the proof of Proposition \ref{prop:major_contribution-alternating}, the global Selmer group $\mathrm{Sel}(T, \alpha : \overline{\chi}_f)$ can be written as a local Selmer group $\mathrm{Sel}(T, \alpha : (\chi_{f,v})_{v \in \Sigma_T(\overline{f})})$. 

    Recall that
    \begin{align*}
        \hat{F}_{(n,N),(w,w')}(\mathbb{F}_q) &= \bigsqcup_{\lambda \in \Lambda_{N,w'}^{la}} \bigsqcup_{\eta \in \Lambda_{n-N,w-w'}^{for}} F_{(n,N),(w,w')}^{(\lambda,\eta)}(\mathbb{F}_q) \\
        &= \bigsqcup_{\lambda \in \Lambda_{N,w'}^{la}} \bigsqcup_{\eta \in \Lambda_{n-N,w-w'}^{for}} \prod_{i,j,k} \mathrm{Conf}_{\lambda_{i,j,k}}(\mathcal{P}_k(i)) \times \prod_{\hat{i},\hat{j},\hat{k}} \mathrm{Conf}_{\eta_{\hat{i},\hat{j},\hat{k}}} (\mathcal{P}_{\hat{k}}(\hat{i})).
    \end{align*}
    Given a fixed $\lambda \in \Lambda_{N,w'}^{la}$ and $\eta \in \Lambda_{n-N,w-w'}^{for}$, we consider two projection maps
    \begin{align}
        \begin{split}
            \Phi_{01} &: \prod_{i,j,k} \mathrm{Conf}_{\lambda_{i,j,k}}(\mathcal{P}_k(i)) \times \prod_{\hat{i},\hat{j},\hat{k}} \mathrm{Conf}_{\eta_{\hat{i},\hat{j},\hat{k}}} (\mathcal{P}_{\hat{k}}(\hat{i})) \to \prod_{i,j} \mathrm{Conf}_{\lambda_{i,j,2}}(\mathcal{P}_2(i)) \times \prod_{\hat{i},\hat{j},\hat{k}} \mathrm{Conf}_{\eta_{\hat{i},\hat{j},\hat{k}}} (\mathcal{P}_{\hat{k}}(\hat{i})) \\
            \Phi_{2} &: \prod_{i,j} \mathrm{Conf}_{\lambda_{i,j,2}}(\mathcal{P}_2(i)) \times \prod_{\hat{i},\hat{j},\hat{k}} \mathrm{Conf}_{\eta_{\hat{i},\hat{j},\hat{k}}} (\mathcal{P}_{\hat{k}}(\hat{i})) \to  \prod_{\hat{i},\hat{j},\hat{k}} \mathrm{Conf}_{\eta_{\hat{i},\hat{j},\hat{k}}} (\mathcal{P}_{\hat{k}}(\hat{i}))
        \end{split}
    \end{align}
    In particular, $\Phi_{01}$ forgets all irreducible factors of $f^*$ in $\mathcal{P}_0 \cup \mathcal{P}_1$, whereas $\Phi_{2}$ forgets all those of $f^*$ in $\mathcal{P}_2$.

    Fix $\lambda \in \Lambda_{N,w'}^{la}$ and $\eta \in \Lambda_{n-N,w-w'}^{for}$. Fix a polynomial $h$ that admits $\eta$. Denote by $\delta_h: \mathbb{Z}_{\geq 0} \to [0,1]$ the probability distribution
    \begin{equation}
        \delta_h(J) := \frac{\#\{f \in \Phi_{01}^{-1}(\Phi_2^{-1}(h)) : \mathrm{rk}((\chi_{f,v})_{v \in \Sigma_T(\overline{h})}) = J\}}{\# \Phi_{01}^{-1}(\Phi_2^{-1}(h))}.
    \end{equation}
    Our first goal is to obtain a governing Markov operator that models the variation of dimensions of local Selmer structures over the preimage $\Phi_{01}^{-1}(\Phi_2^{-1}(h))$. We denote by
    \begin{equation}
        \left\{ \mathcal{W}_{1}^{[i,j,2]}, \mathcal{W}_{2}^{[i,j,2]}, \cdots, \mathcal{W}_{\lambda_{i,j,2}}^{[i,j,2]} \right\}_{\substack{i,j \\ \lambda_{i,j,2} \neq 0}}
    \end{equation}
    the set of auxiliary Markov operators constructed in Proposition \ref{prop:P2_contribution} which models the distribution of dimensions of local Selmer structures over the preimage of $\phi_{i,j,2}^{-1}$ for which $\lambda_{i,j,2} \neq 0$. We define two constants pertaining to the number of irreducible factors of $f^*$:
    \begin{equation}
        w_{\lambda, 1} := \sum_{i,j} \lambda_{i,j,1}, \hspace{15pt} w_{\lambda, 2} := \sum_{i,j} \lambda_{i,j,2}.
    \end{equation}
    We define the following operator
    \begin{equation}
        \mathcal{W}_{(\lambda,\eta)} := \prod_{\substack{i,j \\ \lambda_{i,j,2} \neq 0}} \prod_{b=1}^{\lambda_{i,j,2}} \mathcal{W}_b^{[i,j,2]}.
    \end{equation}
    Using the upper bounds on exponential moments of $\mathcal{W}_{b}^{[i,j,2]}$, we also obtain that
    \begin{equation} \label{eqn:W_lambdaeta_bounds}
        \mathbb{E}[\ell^{\mathcal{W}_{(\lambda,\eta)} \delta_h}] \leq \mathbb{E}[\ell^{\mathcal{V}^{w_{\lambda,2}} \delta_h}].
    \end{equation}
    Then using Proposition \ref{prop:P01_contribution} and Proposition \ref{prop:P2_contribution} by a total of $\omega(\overline{f}^*)$ many iterations, we obtain that for any $m_{n,q} > \max\{\sum_{v \in \Sigma_T} \deg v, 3 \log \ell \}$
    \begin{align}
        \begin{split}
            & \hspace{15pt} \left| \frac{\#\{f \in \Phi_{01}^{-1}(\Phi_2^{-1}(h)) : \dim_{\mathbb{F}_\ell} \mathrm{Sel}(T, \alpha : \overline{\chi}_f) = J\}}{\# \Phi_{01}^{-1}(\Phi_2^{-1}(h))} - \left( \mathcal{M}_\ell^{w_{\lambda,1}} \mathcal{W}_{(\lambda,\eta)} \delta_h \right)(J) \right| \\
            &< \omega(f^*) \cdot B_{T,\ell,q} \cdot (n \log q)^{-2 m_{n,q} + 6 \log \ell + 1}.
        \end{split}
    \end{align}
    Using the upper bound from (\ref{eqn:W_lambdaeta_bounds}) and Theorem \ref{thm:Markov-symmetric}, there exist fixed constants $c_T > 0$ and $0 < \gamma_T < 1$ such that
    \begin{align}
        \begin{split}
            \sup_{J \in \mathbb{Z}_{\geq 0}} \left| (M_\ell^{w_{\lambda,1}} \mathcal{W}_{(\lambda,\eta)} \delta_h)(J) - \prod_{k=1}^\infty \frac{1}{1+ \ell^{-k}} \cdot \prod_{i=1}^J \frac{1}{\ell^i-1} \right| &< c_T \cdot \gamma_T^{w_{\lambda,1}} \cdot \mathbb{E}[\ell^{\mathcal{W}_{(\lambda,\eta)} \delta_h}] \\
            &\leq c_T \cdot \gamma_T^{w_{\lambda,1}} \cdot \mathbb{E}[\ell^{\mathcal{V}^{w_{\lambda,2}} \delta_h}].
        \end{split}
    \end{align}
    Hence we obtain for every $J \in \mathbb{Z}_{\geq 0}$,
    \begin{align}
        \begin{split}
            & \hspace{15pt} \left| \frac{\#\{f \in \Phi_{01}^{-1}(\Phi_2^{-1}(h)) : \dim_{\mathbb{F}_\ell} \mathrm{Sel}(T, \alpha : \overline{\chi}_f) = J\}}{\# \Phi_{01}^{-1}(\Phi_2^{-1}(h))} - \prod_{k=1}^\infty \frac{1}{1+ \ell^{-k}} \cdot \prod_{i=1}^J \frac{1}{\ell^i-1} \right| \\
            &< \left(\omega(f^*) \cdot B_{T,\ell,q} \cdot (n \log q)^{-2 m_{n,q} + 6 \log \ell + 1} + c_T \cdot \gamma_T^{w_{\lambda,1}} \cdot \mathbb{E}[\ell^{\mathcal{V}^{w_{\lambda,2}} \delta_h}] \right).
        \end{split}
    \end{align}
    We take averages of the above expression over elements in $\prod_{\hat{i},\hat{j},\hat{k}}\mathrm{Conf}_{\eta_{\hat{i},\hat{j},\hat{k}}}(\mathcal{P}_{\hat{k}}(\hat{i}))$ by varying $h$. Recall that by definition we have $\omega(h) = w - w'$. Let $D_{w,T}$ be the following positive integer which is the upper bound on $\max_{\chi' \in \Omega_{\overline{h}^*}} \mathrm{rk}(\chi')$ for every such $h$:
    \begin{equation}
        D_{w,T} := 2(w - w') + \max_{\chi \in \Omega_T} \mathrm{rk}(\chi) = 2 \epsilon w + \max_{\chi \in \Omega_T} \mathrm{rk}(\chi)
    \end{equation}
    Denote by $\delta_{D_{w,T}}: \mathbb{Z}_{\geq 0} \to [0,1]$ the probability distribution where $\delta_{D_{w,T}}(D_{w,T}) = 1$ and $0$ otherwise. Then for every $J \in \mathbb{Z}_{\geq 0}$, we have
    \begin{align}
        \begin{split}
            & \hspace{15pt} \left| \frac{\#\{f \in F_{(n,N),(w,w')}^{(\lambda,\eta)}(\mathbb{F}_q) : \dim_{\mathbb{F}_\ell} \mathrm{Sel}(T, \alpha : \overline{\chi}_f) = J\}}{\# F_{(n,N),(w,w')}^{(\lambda,\eta)}(\mathbb{F}_q)} - \prod_{k=1}^\infty \frac{1}{1+ \ell^{-k}} \cdot \prod_{i=1}^J \frac{1}{\ell^i-1} \right| \\
            &< \left((w-w') \cdot B_{T,\ell,q} \cdot (n \log q)^{-2 m_{n,q} + 6 \log \ell + 1} + c_T \cdot \gamma_T^{w_{\lambda,1}} \cdot \mathbb{E}[\ell^{\mathcal{V}^{w_{\lambda,2}} \delta_{D_{w,T}}}] \right).
        \end{split}
    \end{align}
    We take averages of the above expression over $\lambda \in \Lambda_{N,w'}^{la}$ and $\eta \in \Lambda_{n-N,w-w'}^{for}$. This is where we use Chebotarev density theorem over the Galois extension $K(T)/K$ to control the terms $\gamma_T^{w_{\lambda,1}}$ and $\mathbb{E}[\ell^{\mathcal{V}^{w_{\lambda,2}} \delta_{D_{w,T}}}]$. We first use Chebotarev density theorem over the subset of collections of splitting partitions $\{(\lambda_{i,j,k},i,j,k)\}_{i,j,k} \in \Lambda_{N,w'}^{la}$ such that the collection of degrees $i$ at which $\lambda_{i,j,k} \neq 0$ for some $j,k$ is a fixed set. Over such subsets, the following densities can be computed:
    \begin{itemize}
        \item For all but a proportion of at most $1 - (1-\delta(\mathcal{P}_1))^{(1-\delta(\mathcal{P}_1) + \epsilon)w'}$ many elements in the subset, we have $w_{\lambda,1} \geq (\delta(\mathcal{P}_1) - \epsilon) w'$.
        \item For all but a proportion of at most $1 - (1-\delta(\mathcal{P}_2))^{w' - D_{w,T}}$ many elements in the subset, we have $w_{\lambda,2} \geq D_{w,T}$. Such a lower bound on $w_{\lambda,2}$ implies $\mathbb{E}[\ell^{\mathcal{V}^{w_{\lambda,2}} \delta_{D_{w,T}}}] \leq \ell^4$. To see this, we use the definition of the operator $\mathcal{V}$ to obtain that whenever $D_{w,T} \geq 4$,
        \begin{equation*}
            \mathbb{E}[\ell^{\mathcal{V} \delta_{D_{w,T}}}] \leq \ell^{D_{w,T} - 1} = \mathbb{E}[\ell^{\delta_{D_{w,T}-1}}].
        \end{equation*}
        Furthermore, if $D_{w,T} \leq 3$, then
        \begin{equation*}
            \mathbb{E}[\ell^{\mathcal{V} \delta_{D_{w,T}}}] \leq \ell^4.
        \end{equation*}
        We then use induction on the number of iterations of $\mathcal{V}$ to obtain the desired upper bound on the exponential moments.
        
        Even if $w_{\lambda,2} < D_{w,T}$, we have
        \begin{equation*}
            c_T \cdot \gamma_T^{w_{\lambda,1}} \cdot \mathbb{E}[\ell^{\delta_{D_{w,T}}}] = c_T \cdot \gamma_T^{w_{\lambda,1}} \cdot \ell^{2 \epsilon w} \cdot D_{T,\ell,q},
        \end{equation*}
        where we recall that $D_{T,\ell,q}$ was defined as $\ell^{\max_{\chi \in \Omega_T} \mathrm{rk}(\chi)}$. We hence have a trivial bound that holds for any value of $w_{\lambda,2}$:
        \begin{equation*}
            c_T \cdot \gamma_T^{w_{\lambda,1}} \cdot \mathbb{E}[\mathcal{V}^{w_{\lambda,2}}\ell^{\delta_{D_{w,T}}}] \leq c_T \cdot \gamma_T^{w_{\lambda,1}} \cdot \ell^{2 \epsilon w + 4} \cdot D_{T,\ell,q}.
        \end{equation*}
    \end{itemize}
    With these densities computed, we take disjoint unions of such subsets by varying configurations of collections of degrees $i$ to obtain identical densities over the set $\hat{F}_{(n,N),(w,w')}(\mathbb{F}_q)$.
    With these terms in place, we hence obtain the desired statement:
    \begin{align*}
        \begin{split}
            & \hspace{15pt} \left| \frac{\#\{f \in \hat{F}_{(n,N),(w,w')}(\mathbb{F}_q) : \dim_{\mathbb{F}_\ell} \mathrm{Sel}(T, \alpha : \overline{\chi}_f) = J\}}{\# \hat{F}_{(n,N),(w,w')}(\mathbb{F}_q)} - \prod_{k=1}^\infty \frac{1}{1+ \ell^{-k}} \cdot \prod_{i=1}^J \frac{1}{\ell^i-1} \right| \\
            &< \left(2 m_{n,q} \cdot B_{T,\ell,q} \cdot (n \log q)^{-2 m_{n,q} + 6 \log \ell + 1} + \ell^{2 \epsilon w  + 4} \cdot D_{T,\ell,q} \cdot c_T \cdot \gamma_T^{(\delta(\mathcal{P}_1) - \epsilon)(1-\epsilon)w} \right) \\
            &\hspace{15pt} + (1 - \delta(\mathcal{P}_1))^{(1-\delta(\mathcal{P}_1) + \epsilon)(1-\epsilon)w}.
        \end{split}
    \end{align*}
\end{proof}

Theorem \ref{thm:Galois_module-symmetric} hence follows by combining Proposition \ref{proposition:fan_approximation} and Proposition \ref{prop:major_contribution-symmetric}, and optimizing the rate of convergence by choosing $\epsilon := (\log \log m_{n,q})^{-1}$. As in the proof of Theorem \ref{thm:Galois_module-alternating}, we perform a line-by-line adaptation of \cite[Proof of Theorem 1.2, pp.3316--3317]{park2022prime}. The implied constant $B_{sym}$ is given by
\begin{equation}
    B_{sym} := 12 \cdot (B_{T,\ell,q} + c_T) \cdot D_{T,\ell,q}.
\end{equation}
The first two terms of the exponent $\alpha(T)$ are obtained from Proposition \ref{proposition:fan_approximation}, whereas the last two terms are obtained from Proposition \ref{prop:major_contribution-symmetric}.

%% file: Application.tex
\section{Applications to twists of some superelliptic curves}
\label{sec:superelliptic}

We recall that $C$ is a non-isotrivial superelliptic curve over $K$ whose Weierstrass equation is
\begin{equation}
    C: y^\ell = F(x),
\end{equation}
such that $F$ is a degree $3$ polynomial over $\mathbb{F}_q(t)$. We note that $C$ is a curve of genus $\ell-1$.

Because we assumed that $\mu_\ell \subset K$, the isogeny $1-\zeta_\ell: \mathrm{Jac}(C) \to \mathrm{Jac}(C)$ induced from the map $\zeta_\ell: C \to C$ defined as $(x,y) \mapsto (x,\zeta_\ell y)$ is a well-defined endomorphism of $\mathrm{Jac}(C)$ over $K$. Hence, we let our $\Gal(\overline{K}/K)$-module of interest to be
\begin{equation}
    T := \mathrm{Jac}(C)[1-\zeta_\ell].
\end{equation}

Given a polynomial $f \in \mathbb{F}_q[t]$, we recall that $C_f$ is a cyclic order $\ell$ twist of $C$ over $K$ whose Weierstrass equation can be written as
\begin{equation}
    C_f: f y^\ell = F(x).
\end{equation}
As an application of Theorem \ref{thm:Galois_module-symmetric}, we compute the distribution of $1-\zeta_\ell$ Selmer groups of $\mathrm{Jac}(C_f)$ as $f$ ranges over the set of polynomials of degree $n$.
\begin{theorem}
\label{thm:superelliptic}
    Fix a prime number $\ell \neq 2,3$. Let $K = \mathbb{F}_q(t)$ be a global function field whose characteristic $p$ is coprime to $2,3$, and $q \equiv 1 \text{ mod } \ell$. Let $C: y^\ell = F(x)$ be a non-isotrivial superelliptic curve over $K$. Suppose the following four conditions hold.
    \begin{itemize}
        \item $C$ is not defined over $K^p$.
        \item The constant field of $K(\mathrm{Jac}(C)[1-\zeta_\ell])$ is equal to $\mathbb{F}_q$.
        \item $\Gal(K(\mathrm{Jac}(C)[1-\zeta_\ell])/K) \cong S_3$,
        \item $\mathrm{Jac}(C)$ has a place $v$ of totally split multiplicative reduction.
    \end{itemize}
    Denote by $\gamma_\ell$ the geometric rate of convergence of the Markov operator $M_\ell$ from Theorem \ref{thm:Markov-symmetric}. 
    Then for any small enough $\delta > 0$, there exist sufficiently large $n$ and a fixed constant $B(C,q) > 0$ independent of $n$ such that
    \begin{align*}
        \left| \frac{\# \{f \in F_n(\mathbb{F}_q) \; | \; \dim_{\mathbb{F}_\ell} \Sel_{1-\zeta_\ell}(\mathrm{Jac}(C_f)/K) = r \}}{\#F_n(\mathbb{F}_q)} - \prod_{k=1}^\infty \frac{1}{1+\ell^{-k}} \prod_{i=1}^r \frac{1}{\ell^i-1} \right| < \frac{B(C,q)}{n^{\alpha(\ell) - \delta}},
    \end{align*}
    where
    \begin{equation*}
        \alpha(\ell) := \sup_{0 < \rho < 1} \left( \min \left( \rho \log \rho + 1 - \rho, \; \; -\frac{\rho}{2} \cdot \log \gamma_\ell, \; \; \frac{\rho}{2} \log 2 \right) \right).
    \end{equation*}
\end{theorem}

Using Theorem \ref{thm:superelliptic}, we prove the following theorem, namely the uniform upper bound on the moments of the number of $K$-rational solutions of superelliptic curves $C_f$ in the twist family $\mathcal{C}_n$.
\begin{theorem} \label{thm:main}
    Assume conditions from Theorem \ref{thm:superelliptic}. Given a positive integer $m$, we denote by $S(\ell, p)$, $\mathbb{E}(\ell, p, m)$, and $\mathbb{P}(\ell, m)$ the following constants:
    \begin{align*}
        S(\ell, p) &:= 3^{\ell - 1} \cdot p^{3\ell - 3} \cdot (8\ell - 10) \cdot (\ell - 1)!, \\
        \mathbb{E}(\ell, p, m) &:= \prod_{k=1}^\infty \frac{1}{1+\ell^{-k}} \cdot \sum_{r = 0}^\infty \prod_{i=1}^r \frac{p^{(\ell-1)m}}{\ell^i-1},\\
        \mathbb{P}(\ell, m) &:= \prod_{k=1}^\infty \frac{1}{1+\ell^{-k}} \cdot \sum_{0 \leq r \leq m} \prod_{i=1}^r \frac{1}{\ell^i-1}.
    \end{align*}
    Then for each positive integer $m$, the following two inequalities hold:
    \begin{align}
        \limsup_{n \to \infty} \frac{\sum_{\substack{f \in F_n(\mathbb{F}_q)}} (\# C_f(K))^m}{\#F_n(\mathbb{F}_q)} &\leq \mathbb{E}(\ell, p, m) \cdot S(\ell,p)^m \\
        \liminf_{n \to \infty} \frac{\#\{f \in F_n(\mathbb{F}_q) \; | \; \#C_f(K) \leq p^{m(\ell - 1)} \cdot S(\ell,p)\}}{\# F_n(\mathbb{F}_q)} &\geq \mathbb{P}(\ell, m)
    \end{align}
\end{theorem}
\begin{proof}[Proof of Theorem \ref{thm:main}]
    By adapting the proof of Proposition 3.5 of \cite{Sch98} to global function fields, we have 
    \begin{equation}
        \mathrm{rank}_\mathbb{Z}(\mathrm{Jac}(C_f)(K)) \leq (\ell-1) \cdot \dim_{\mathbb{F}_\ell} \Sel_{1-\zeta_\ell}(\mathrm{Jac}(C_f)/K).
    \end{equation}
    Because $C$ is not defined over $K^p$, it follows that for any $f \in F_n(\mathbb{F}_q)$, $C_f$ is also not defined over $K^p$. This is because we have an isomorphism $C \cong C_f$ over $K(\sqrt[\ell]{f})$. If it were the case that $C_f$ is defined over $(K(\sqrt[\ell]{f}))^p = K^p(\sqrt[\ell]{f})$, then $C$ is also defined over $K^p(\sqrt[\ell]{f})$. But because the model of $C$ is already defined over $K$ and $\ell$ is coprime to $p$, we obtain that $C$ is in fact defined over $K_p$, a contradiction.
    
    Hence for any $f \in F_n(\mathbb{F}_q)$, we can use the main theorem of \cite{BV96} and the fact that the genus of $C_f$ is equal to $\ell - 1$ to obtain
    \begin{align}
    \begin{split}
        \# C_f(K) &\leq \# \frac{\mathrm{Jac}(C_f)(K)}{p \cdot \mathrm{Jac}(C_f)(K)} \cdot (3p)^{\ell-1} \cdot (8\ell - 10) \cdot (\ell-1)! \\
        &\leq p^{(\ell - 1) \dim_{\mathbb{F}_\ell} \Sel_{1-\zeta_\ell}(\mathrm{Jac}(C_f)/K)} \cdot 3^{\ell - 1} \cdot p^{3\ell - 3} \cdot (8 \ell - 10) \cdot (\ell - 1)!.
    \end{split}
    \end{align}
    Applying Theorem \ref{thm:superelliptic} gives the statement of the main theorem. The fact that $\mathbb{E}(\ell, p, m) < \infty$ for each $m \geq 1$ follows from the fact that the quantity $\prod_{i=1}^r \frac{1}{\ell^i-1}$ decays quadratic exponentially, i.e. at a rate of $O(\ell^{-r(r+1)/2})$, as $r$ grows arbitrarily large.
\end{proof}
The upper bound computed in Theorem \ref{thm:main} is far from impressive. Nevertheless, we provide in Table \ref{tab:my_label} some approximate values of $\log_p \mathbb{E}(\ell,p,1)$, $\mathbb{P}(\ell,0)$, and $\mathbb{P}(\ell,1)$, where $p, \ell \in \{5,7,11,13,17\}$. For the probabilities $\mathbb{P}(\ell,0)$ and $\mathbb{P}(\ell,1)$, we compute  their asymptotic values when $\ell$ grows arbitrarily large, whose implied constants of the error terms are positive and independent of the choice of $C$. In particular, at least $99\%$ of curves $C_f \in \mathcal{C}_n$ satisfies the inequality $\#C_f(K) \leq p^{(\ell-1)} \cdot S(\ell,p)$ as $n$ grows arbitrarily large.
\begin{table}[ht]
    \centering
    \begin{tabular}{|c||c|c|c|c|c|c|}
    \hline
    & $\ell = 5$ & $\ell = 7$ & $\ell = 11$ & $\ell = 13$ & $\ell = 17$ & $\ell \to \infty$ \\
    \hline
    \hline
    $\log_5 \mathbb{E}(\ell,5,1)$ & $\times$ & $12.41584$ & $29.05016$ & $39.65991$ & $65.18045$ & $-$ \\
    \hline
    $\log_7 \mathbb{E}(\ell,7,1)$ & $8.14741$ & $\times$ & $35.98102$ & $49.02079$ & $80.30054$ & $-$\\
    \hline
    $\log_{11} \mathbb{E}(\ell,11,1)$ & $10.30474$ & $19.53559$ & $\times$ & $61.63361$ & $100.64901$ & $-$\\
    \hline
    $\log_{13} \mathbb{E}(\ell,13,1)$ & $11.11003$ & $21.05772$ & $48.79132$ & $\times$ & $108.17569$ & $-$\\
    \hline
    $\log_{17} \mathbb{E}(\ell,17,1)$ & $12.40935$ & $23.50785$ & $54.35576$ & $73.80376$ & $\times$ & $-$\\
    \hline
    \hline
    $\mathbb{P}(\ell,0)$ & $0.79334$ & $0.85459$ & $0.90840$ & $0.92265$ & $0.94098$ & $1 - O(1/\ell)$ \\
    \hline
    $\mathbb{P}(\ell,1)$ & $0.99167$ & $0.99702$ & $0.99924$ & $0.99954$ & $0.99980$ & $1 - O(1/\ell^3)$ \\
    \hline
    \end{tabular}
    \vspace{0.5pt}
    \caption{A table of approximate values of $\log \mathbb{E}(\ell,p,1)$, $\mathbb{P}(\ell,0))$, and $\mathbb{P}(\ell,1)$ for some first $5$ values of $\ell$ and $p$ up to 5 decimals. For the probabilities, we also include the asymptotic values of these constants as $\ell$ grows arbitrarily large with error terms. We omit the cases where $\ell = p$ for which our main theorem is not applicable.}
    \label{tab:my_label}
\end{table}

Before we prove Theorem \ref{thm:superelliptic}, we state and prove a few lemmas.
\begin{lemma} \label{lemma:T_superelliptic}
The following properties hold for the superelliptic curve $C: y^\ell = F(x)$ over $K$ such that $F$ is a degree $3$ polynomial over $\mathbb{F}_q(t)$.
\begin{itemize}
    \item The torsion subgroup $\mathrm{Jac}(C)[1-\zeta_\ell]$ is a 2-dimensional $\mathbb{F}_\ell$ vector space generated by a set of divisors $\{(P_2 - P_1), (P_3 - P_1)\}$, where $P_i$'s are three Weierstrass points of $C$.
    \item For any polynomial $f \in \mathbb{F}_q[t]$, there is a $\Gal(\overline{K}/K)$-equivariant isomorphism
    \begin{equation}
        \mathrm{Jac}(C_f)[1-\zeta_\ell] \cong \mathrm{Jac}(C)[1-\zeta_\ell]
    \end{equation}
\end{itemize}
\end{lemma}
\begin{proof}
    The proof follows from \cite[Proposition 1.7]{Yu16} and \cite[Proposition 3.2]{Sch98}. The details of the proof for the second statement of the lemma can be obtained as follows. Denote by $\alpha_1, \alpha_2,\alpha_3$ the roots of the cubic polynomial $F$ in $\overline{K}$. Then we can rewrite the Weierstrass equation for $C$ as
    \begin{equation*}
        y^\ell = (x-\alpha_1)(x-\alpha_2)(x-\alpha_3).
    \end{equation*}
    As a convention, we denote by $P_1 = (\alpha_1,0)$, $P_2 = (\alpha_2,0)$, and $P_3 = (\alpha_3,0)$.
    Let $b$ be an integer such that $b\ell \equiv 1 \text{ mod } 3$. We observe that multiplying both sides of the Weierstrass equation for $C_f$ by $f^{b\ell - 1}$ and applying change of coordinates $(x^{\frac{b\ell-1}{3}},y^b) \mapsto (X,Y)$ gives
    \begin{equation*}
        Y^\ell = (X - f^{\frac{b\ell - 1}{3}} \alpha_1)(X - f^{\frac{b\ell - 1}{3}} \alpha_2)(X - f^{\frac{b\ell - 1}{3}} \alpha_3)
    \end{equation*}
    Because $f^{\frac{b\ell-1}{3}} \in \mathbb{F}_q[t]$, it follows that the bijection defined as
    \begin{align*}
        \mathrm{Jac}(C_f)[1-\zeta_\ell] &\to \mathrm{Jac}(C)[1-\zeta_\ell], \\
        ((f^{\frac{b\ell - 1}{3}} \alpha_2,0) - (f^{\frac{b\ell - 1}{3}} \alpha_1,0)) &\mapsto ((\alpha_2,0) - (\alpha_1,0)), \\
        ((f^{\frac{b\ell - 1}{3}} \alpha_2,0) - (f^{\frac{b\ell - 1}{3}} \alpha_1,0)) &\mapsto ((\alpha_3,0) - (\alpha_1,0)),
    \end{align*}
    is a $\Gal(\overline{K}/K)$-equivariant isomorphism.
\end{proof}

\begin{lemma} \label{lemma:T_conditions}
    Suppose the superelliptic curve $C$ satisfies the four conditions from Theorem \ref{thm:superelliptic}. Then $\mathrm{Jac}(C)[1-\zeta_\ell]$ is a symmetric admissible $\Gal(\overline{K}/K)$-module.
\end{lemma}
\begin{proof}
    Lemma \ref{lemma:T_superelliptic} implies statement (1) from Condition \ref{condition:admissible}. The conditions from Theorem \ref{thm:superelliptic} that $C$ is non-isotrivial, is not defined over $K^p$, and $K$ contains $\mu_\ell$ imply statement (2) of Condition \ref{condition:admissible}.

    We now show that the second condition from Theorem \ref{thm:superelliptic} implies statements (3),(4), and (5) of Condition \ref{condition:admissible}. We represent $(P_2 - P_1)$ and $(P_3 - P_1)$ as column vectors with entries in $\mathbb{F}_\ell$:
    \begin{equation*}
        (P_2 - P_1) \mapsto \begin{pmatrix} 1 \\ 0 \end{pmatrix}, \hspace{4pt} (P_3 - P_1) \mapsto \begin{pmatrix} 0 \\ 1 \end{pmatrix}.
    \end{equation*}
    We recall that the absolute Galois group $\Gal(\overline{K}/K)$ acts on $\mathrm{Jac}(C)[1-\zeta_\ell]$ by mapping the divisor $[P_i - P_j]$ to $[\sigma P_i - \sigma P_j]$ for any $\sigma \in \Gal(\overline{K}/K)$. The elements of $\Gal(K(\mathrm{Jac}(C)[1-\zeta_\ell])/K) \cong S_3$ can be represented by elements in $\GL_2(\mathbb{F}_\ell)$ as:
    \begin{align} \label{eqn:superelliptic_matrix-rep}
    \begin{split}
        () \mapsto \begin{pmatrix} 1 & 0 \\ 0 & 1 \end{pmatrix}, \hspace{4pt} (1 2) &\mapsto \begin{pmatrix} -1 & -1 \\ 0 & 1 \end{pmatrix}, \hspace{4pt} (1 3) \mapsto \begin{pmatrix} 1 & 0 \\ -1 & -1 \end{pmatrix}, \\
        (2 3) \mapsto \begin{pmatrix} 0 & 1 \\ 1 & 0 \end{pmatrix}, \hspace{4pt} (1 2 3) &\mapsto \begin{pmatrix} -1 & -1 \\ 1 & 0 \end{pmatrix}, \hspace{4pt} (1 3 2) \mapsto \begin{pmatrix} 0 & 1 \\ -1 & -1 \end{pmatrix}.
    \end{split}
    \end{align}
    The group acts on the divisors via matrix multiplication. Statement (3) follows from the fact that $\ell \geq 5$ and there are no proper subspaces of $\mathbb{F}_\ell^{\oplus 2}$ that is invariant with respect to multiplication by $6$ of the matrices above. Statement (4) can be checked from the fact that the only matrix in $\text{M}_{2 \times 2}(\mathbb{F}_\ell)$ such that commutes with all the $6$ matrices above are scalar multiples of the identity matrix. Statement (5) follows from the fact that the order of the group $S_3$ is coprime to $\ell$.

    The third condition from Theorem \ref{thm:superelliptic} implies that at the place $v$ the abelian variety $\mathrm{Jac}(C)$ is isomorphic to $\mathbb{G}_m^{\ell-1} / \prod_{i=1}^{\ell-1} q_i^\mathbb{Z}$ for some elements $q_1, q_2, \cdots, q_{\ell - 1}$ of $\mathbb{G}_m$ with non-trivial valuations with respect to $v$. It follows that $\mathrm{Jac}(C)[\ell](\overline{K}_v) \cap \mu_{\ell^\infty} = \mu_\ell$, from which Statement (6) follows. This concludes the proof that $\mathrm{Jac}(C)[1-\zeta_\ell]$ is an admissible $\mathrm{Gal}(\overline{K}/K)$-module.

    To show that $\mathrm{Jac}(C)[1-\zeta_\ell]$ is symmetric, it suffices to show that the Weil pairing
    \begin{equation*}
        \mathrm{Jac}(C)[1-\zeta_\ell] \times \mathrm{Jac}(C)[1-\zeta_\ell] \to \mu_\ell
    \end{equation*}
    is a non-degenerate symmetric pairing. This follows from using \cite[Lemmas 10.1, 10.2]{PS97} where we choose the field $k = \mathbb{F}_q(t)$ that contains primitive $\ell$-th roots of unity.
\end{proof}

We now prove Theorem \ref{thm:superelliptic}. By Lemma \ref{lemma:T_superelliptic}, we have
\begin{equation*}
    \Sel_{1-\zeta_\ell}(\mathrm{Jac}(C_f)/K) \subset H^1_\et(K, \mathrm{Jac}(C)[1-\zeta_\ell]).
\end{equation*}
Let $\alpha$ be a twisting data defined as
\begin{equation}
    \alpha(\overline{\chi}_f) := \begin{cases}
        \mathrm{Image} \left( \delta_v: \frac{\mathrm{Jac}(C_f)(K_v)}{(1-\zeta_\ell)(\mathrm{Jac}(C_f)(K_v))} \to H^1_{\et}(K_v, \mathrm{Jac}(C)[1-\zeta_\ell]) \right) &\text{ if } v \in \Sigma \text{ or } v \mid f, \\
        H^1_{\et}(\Oh_{K_v}, \mathrm{Jac}(C)[1-\zeta_\ell]) &\text{ otherwise }.
    \end{cases}
\end{equation}
We observe that
\begin{equation}
    \Sel_{1-\zeta_\ell}(\mathrm{Jac}(C_f)/K) = \Sel(\mathrm{Jac}(C)[1-\zeta_\ell], \alpha : \overline{\chi}_f).
\end{equation}
By Lemma \ref{lemma:T_conditions}, the torsion subgroup $\mathrm{Jac}(C)[1-\zeta_\ell]$ is an admissible $\Gal(\overline{K}/K)$-module. To apply Theorem \ref{thm:Galois_module-symmetric}, we need to check that $\delta(\mathcal{P}_1) \neq 0$ for the module $\mathrm{Jac}(C)[1-\zeta_\ell]$. For this we recall the matrix representations of elements in $\mathrm{Gal}(K(\mathrm{Jac}(C)[1-\zeta_\ell])/K)$, as shown in (\ref{eqn:superelliptic_matrix-rep}). Using the proof of Lemma \ref{lemma:H_ram}, we obtain
\begin{equation*}
        v \in \mathcal{P}_i \iff \dim_{\mathbb{F}_\ell} \frac{\mathrm{Jac}(C)[1-\zeta_\ell]}{(\mathrm{Frob}_v - 1) \mathrm{Jac}(C)[1-\zeta_\ell]} = i.
\end{equation*}
Using (\ref{eqn:superelliptic_matrix-rep}), we obtain the following equivalence relation:
\begin{equation}
    v \in \mathcal{P}_i \iff \begin{cases}
        \mathrm{Frob}_v \in \{(123), (132)\} &\text{ if } i = 0,\\
        \mathrm{Frob}_v \in \{(12), (13), (23) \} &\text{ if } i = 1, \\
        \mathrm{Frob}_v \in \{()\} &\text{ if } i = 2.
    \end{cases}
\end{equation}
We then use Chebotarev density theorem for $S_3$ Galois extensions to prove $\delta(\mathcal{P}_0) = \frac{1}{3}$, $\delta(\mathcal{P}_1) = \frac{1}{2}$, and $\delta(\mathcal{P}_2) = \frac{1}{6}$. Hence, we use Theorem \ref{thm:Galois_module-symmetric} to prove Theorem \ref{thm:superelliptic}.

\begin{remark} \label{remark:PR}
    We compare the statement of Theorem \ref{thm:superelliptic} with the predicted distribution from Poonen-Rains heuristics \cite{PR12}. Because $\mathrm{Jac}(C)[1-\zeta_\ell]$ is a direct summand of the permutation $\mathbb{F}_\ell$ vector space arising from the set of roots of $F$, Proposition 3.3.(b) and Theorem 4.14 of \cite{PR12} implies that $\Sel_{1-\zeta_\ell}(\mathrm{Jac}(C_f)/K)$ is an intersection of the images of the two homomorphisms
    \begin{equation}
        \begin{tikzcd}
            & H^1_\et(K,\mathrm{Jac}(C)[1-\zeta_\ell]) \arrow{d}{} \\
\prod_{v \text{ place of } K} \frac{\mathrm{Jac}(C_f)(K_v)}{(1-\zeta_\ell)\mathrm{Jac}(C_f)(K_v)} \arrow{r}{} & \prod_{v \text{ place of } K} H^1_\et(K_v, \mathrm{Jac}(C)[1-\zeta_\ell]),
        \end{tikzcd}
    \end{equation}
    both of which are maximal isotropic subspaces of the codomain of two homomorphisms with respect to the pairing $\{\mu_{T,v}\}$. However, unlike the Poonen-Rains heuristics which assume that $\mu_{T,v}$ is a symmetric pairing, our setup states that $\mu_{T,v}$ is an alternating pairing. 
    The distribution of intersections of two random maximal isotropic subspaces with respect to alternating pairing appears in the paper by Cais, Ellenberg, and Zureick-Brown \cite[Proposition 3.5]{Cais_Ellenberg_Zureick-Brown_2013}. Though the distribution was obtained for a different arithmetic application regarding random Dieudonn\'e modules, it predicts the distribution of the $1-\zeta_\ell$ Selmer groups of $\mathrm{Jac}(C_d)$ to be
    \begin{equation}
        \mathbb{P}[\dim_{\mathbb{F}_\ell} \mathrm{Sel}_{1-\zeta_\ell}(\mathrm{Jac}(C_d)) = r] =_? \prod_{k=1}^\infty \frac{1}{1 + \ell^{-k}} \cdot \prod_{i=1}^r \frac{1}{\ell^i - 1}.
    \end{equation}
    And indeed, Theorem \ref{thm:Galois_module-symmetric} achieves this distribution as a stationary distribution of the Markov operator $M_\ell$ from Theorem \ref{thm:Markov-symmetric}.

    One may ask whether the assumptions of the Poonen-Rains heuristics are still applicable to our setup. These are:
    \begin{itemize}
        \item One of the images of the two homomorphisms is fixed.
        \item The other image of the two homomorphisms varies uniformly at random over the set of maximal isotropic subspaces of $\prod_{v \text{ place of } K} H^1_\et(K_v, \mathrm{Jac}(C)[1-\zeta_\ell])$ with respect to $\mu_{T,v}$.
    \end{itemize}
    The Galois-equivariant isomorphism obtained from Lemma \ref{lemma:T_superelliptic} implies that the image of the restriction map (or the vertical map from the commutative diagram) satisfies the first condition. However, the description of the twisting data $\alpha$ defining the Selmer group $\Sel_{1-\zeta_\ell}(\mathrm{Jac}(C_f)/K)$ does not satisfy the second condition. This is because given any place $v \in \mathcal{P}_2$, the heuristic assumes that the images of the local Kummer maps (or the horizontal map from the commutative diagram) vary uniformly at random over all maximal isotropic subspaces of $H^1_\et(K_v, \mathrm{Jac}(C)[1-\zeta_\ell])$, whereas the proof of Theorem \ref{thm:superelliptic} shows that the images of the local Kummer maps range over only $\ell + 1$ of the maximal isotropic subspaces of $H^1_\et(K_v, \mathrm{Jac}(C)[1-\zeta_\ell])$ for each place $v \in \mathcal{P}_2$.
\end{remark}